\documentclass{amsproc}
\usepackage{overpic, graphicx, mathrsfs, amssymb, amsmath}

    \newtheoremstyle{upright}%
        {8pt plus2pt minus4pt}%
        {8pt plus2pt minus4pt}%
        {\upshape}%
        {}%
        {\bfseries\scshape}%
        {\bf{.}}%
        {1em}%
        {}%

\newtheorem{theorem}{Theorem}[section]
\newtheorem{example}[theorem]{Example}
\newtheorem{lemma}[theorem]{Lemma}
\newtheorem{proposition}[theorem]{Proposition}
\newtheorem{corollary}[theorem]{Corollary}
\newtheorem{conjecture}[theorem]{Conjecture}

\newtheorem{statement}[theorem]{Statement}
\newtheorem{notation}[theorem]{Notation}

\theoremstyle{upright}
\newtheorem{definition}[theorem]{Definition}
\newtheorem{remark}[theorem]{Remark}
\newcommand{\yooi}[1]{{(y_0^0)}_{#1}}
\newcommand{\xooi}[1]{{(x_0^0)}_{#1}}
\newcommand{\xoo}{x_0^0}
\newcommand{\yoo}{y_0^0}
\newcommand{\xio}[1]{x_{#1}^0}
\newcommand{\yio}[1]{y_{#1}^0}
\newcommand{\xii}[2]{x_{#1}^{#2}}
\newcommand{\yii}[2]{y_{#1}^{#2}}

\newcommand{\ppfi}[1]{pp\mathscr{F}_{#1}}
\newcommand{\pfix}[2]{p\mathscr{F}_{#1}^{#2}}

\newcommand{\pfx}[1]{p\mathscr{F}^{#1}}

\newcommand{\cantor}{\mathscr{C}}
\newcommand{\omegaks}{\Omega(KS)}
\newcommand{\omegaksi}[1]{\Omega(KS_{#1})}
\newcommand{\compseq}{\{\mathscr{O}_i(x_i^0,\theta_i^0)\}_{i=0}^\infty}
\newcommand{\compseqang}[1]{\{\mathscr{O}_i(x_i^0,#1)\}_{i=0}^\infty}
\newcommand{\compseqangx}[2]{\{\mathscr{O}_i(#2_i^0,#1)\}_{i=0}^\infty}
\newcommand{\compseqangi}[2]{\{\mathscr{O}_i(x_i^0,#1)\}_{i=#2}^\infty}
\newcommand{\orbitksi}[1]{\mathscr{O}_{#1}(x_{#1}^0,\theta_{#1}^0)}
\newcommand{\orbitksiang}[2]{\mathscr{O}_{#1}(x_{#1}^0,#2)}

\newcommand{\orbitksix}[2]{\mathscr{O}_{#1}(#2_{#1}^0,\pi/3)}
\newcommand{\orbitksixang}[3]{\mathscr{O}_{#1}(#2_{#1}^0,#3)}
\newcommand{\orbitksixNOi}[2]{\mathscr{O}_{#1}(#2_{#1}^0,\pi/3)}
\newcommand{\orbitksixNOiang}[3]{\mathscr{O}_{#1}({#2}_{#1}^{0},#3)}
\newcommand{\orbitks}{\underleftarrow{\mathscr{O}(x^0,\pi/3)}}

\newcommand{\FL}{\mathscr{L}}
\newcommand{\sumpieces}[1]{\sum_{i=2}^{#1}\chi[\xooi{i}]\#\orbitksiang{i-1}{\pi/3}\frac{\FL}{3^i}}

 \newcommand{\hstr}[1][3]{\rule{#1ex}{0ex}}  %horizontal strut
\numberwithin{equation}{section} 
\begin{document}
\title{Families of Periodic Orbits of the Koch Snowflake Fractal Billiard}

\author{Michel L. Lapidus}
\address{Department of Mathematics, University of California, Riverside, CA 92521--0135, USA}

\email{lapidus@math.ucr.edu}
\thanks{The research of the first author (MLL) was supported in part by NSF Grant DMS-0707524.}

\author{Robert G. Niemeyer}
\address{Department of Mathematics, University of California, Riverside, CA 92521--0135, USA}

\email{niemeyer@math.ucr.edu}

%    General info
\date{June 12, 2011.}


\begin{abstract}
The Koch snowflake $KS$ is a nondifferentiable curve.  Hence, any attempt to define reflection in the boundary may seem like an exercise in futility.  In this paper, for each integer $n\geq 0$, we describe the periodic orbits of the prefractal billiard $\omegaksi{n}$ (the $n$th inner rational polygonal approximation of the Koch snowflake billiard).  Moreover, we use this information in order to define and describe a particular collection of periodic orbits of the Koch snowflake billiard $\omegaks$.  

In the finite case, an orbit of $\omegaksi{n}$ can be reduced to its Poincar\'e section, which simply amounts to a finite collection of points in the boundary of the prefractal billiard.  We show that, for each $n\geq 0$, the collection of directions for which the billiard flow on $\omegaksi{n}$ is closed is exactly the collection of directions for which the billiard flow on $\omegaksi{0}$ is closed.  Such a result relies on the fact that the corresponding flat surface $\mathcal{S}(KS_n)$ ($n\geq 1$) is shown to be a branched cover of the flat surface $\mathcal{S}(KS_0)$, the hexagonal torus. Extending this result, we define what we call a compatible sequence of periodic orbits.  Focusing on the direction given by an initial angle of $\pi/3$, we define 1) a compatible sequence of piecewise Fagnano orbits, 2) an eventually constant compatible sequence of orbits and 3) a compatible sequence of generalized piecewise Fagnano orbits.

In the case of the infinite (fractal) billiard table, we will describe what we call \textit{stabilizing periodic orbits} of the Koch snowflake billiard $\omegaks$.  An eventually constant compatible sequence of periodic orbits is comprised (for all but finitely many) of $\cantor$-orbits (or what we also call \textit{stabilizing periodic orbits}). We show that the trivial limit of an eventually constant compatible sequence of periodic orbits is, in fact, a periodic orbit of $\omegaks$.  In a sense, we show that it is possible to define billiard dynamics on a Cantor set.

%(Although, finitely many of the orbits in an eventually constant compatible sequence, which are found only at the beginning of the sequence, are not necessarily stabilizing periodic orbits.)

In addition, we will discuss the geometric and topological properties of what we call the \textit{footprint of a piecewise Fagnano orbit}.  We will show that the inverse limit of the footprints of orbits of the prefractal approximations (or, the Poincar\'e sections of the respective orbits) exists in a specific situation and provide a plausibility argument as to why such an inverse limit of footprints should constitute the footprint of a well-defined periodic orbit of $\omegaks$.  Using, in particular, known results for the inverse limit of a sequence of finite spaces, we deduce that the \textit{footprint} (i.e., the intersection of the orbit with the boundary) of a \textit{piecewise Fagnano orbit} is a topological Cantor set and even, a self-similar Cantor set.
%To demonstrate the existence of what we will call \textit{stabilizing periodic} orbits and footprints of piecewise Fagnano orbits, we proceed through the following steps: 1) show an equality of the sets of periodic directions in $\omegaksi{0}$ and $\omegaksi{n}$, for all $n\geq 0$, 2) describe orbits with an initial direction of $\pi/3$ in $\omegaksi{n}$, $n\geq 0$, entirely in terms of the ternary representation of the initial basepoint $\xoo$ of an orbit $\orbitksi{0}$ of the equilateral triangle billiard $\omegaksi{0}$, and 3) demonstrate the existence of stabilizing periodic orbits (or, sometimes referred to as $\mathscr{C}$-orbits) of the Koch snowflake billiard $\omegaks$ by showing that a particular compatible sequence of periodic orbits is eventually constant. 

%We also give a characterization and a description of the periodic orbits of the Koch snowflake billiard, where such orbits have an initial direction of $\pi/3$, and finally, discuss the geometric and topological properties of such periodic orbits.  
%In addition, we provide computer simulations in support of the existence of other periodic orbits of $\omegaks$, those being what we call \textit{piecewise Fagnano} and \textit{generalized piecewise Fagnano} orbits of $\omegaks$.  
% and the union of all footprints from such orbits constitutes a subset of the \textit{natural Cantor set} of the Koch snowflake curve $KS$.

We allude to a possible characterization of orbits with an initial direction of $\pi/3$.  That is, we provide support for a complete description of periodic orbits in the direction of $\pi/3$.  Such a characterization would allow one to describe an orbit with an initial direction of $\pi/3$ of the Koch snowflake billiard as either a piecewise Fagnano orbit, a stabilizing orbit or a generalized piecewise Fagnano orbit. We then close the paper by discussing several outstanding open problems and conjectures about the Koch snowflake billiard $\omegaks$, the associated `\textit{fractal flat surface}' and possible connections with the associated fractal drum $\mathcal{D}(KS)$ via fractal analogues of Gutzwiller-like trace formulae.

%In particular, we are able to describe such orbits entirely in terms of the ternary expansion of the initial basepoint $\xoo\in I$ of the first orbit $\orbitksi{0}$ of the equilateral triangle billiard $\omegaksi{0}$.  This gives rise to classification of periodic orbits of the Koch snowflake billiard $\omegaks$, where we stipulate that this classification assumes a fixed initial direction of $\pi/3$ for every orbit being considered.  We will show that such orbits have finite length.  Moreover, a certain subcollection of periodic orbits, which is comprised of what we call \textit{piecewise Fagnano orbits} and \textit{generalized piecewise Fagnano orbits} (see Figures \ref{fig:piecewiseFagnanoOrbits} and \ref{fig:fagnanos012}), constitute 1) the \textit{natural Cantor set} of the Koch snowflake $KS$ and 2) an individual orbit from this particular subcollection of orbits constitutes a topological Cantor set in its own right.

These problems and conjectures have natural counterparts for other fractal billiards.  In the long-term, the present work may help lay the foundations for a general theory of fractal billiards.  
\end{abstract}

\maketitle

%\newpage

\section{Introduction}
\label{sec:introduction}
\begin{figure}[htbp]
	\centering
		\includegraphics{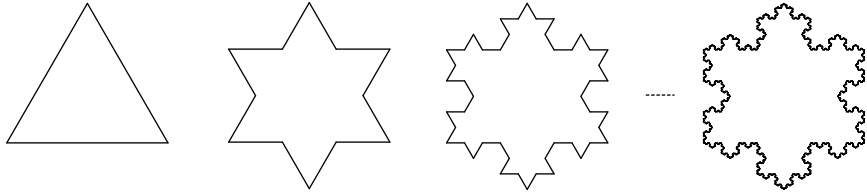}
	\caption{The Koch snowflake curve $KS$ and its prefractal approximations $KS_n$, for $n=0,1,2,...,6$.}
	\label{fig:kochconstruction75}
\end{figure}

\begin{figure}
	\centering
		\includegraphics{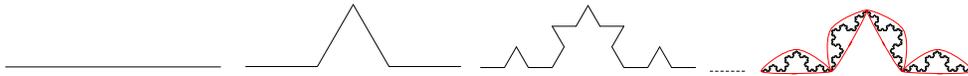}
	\caption{The Koch curve $KC$ and its prefractal approximations $KC_n$, for $n=0,1,2,...,6$.  The Koch curve is \textit{self-similar}, meaning that there are scaled (and rotated/reflected) copies of the Koch curve found as subsets of $KC$.  In this figure, we circle four copies of $KC$ scaled by $1/3$.  In general, one can find $4^n$ copies of the Koch curve scaled by $1/3^n$ as subsets of $KC$ such that the disjoint union (disjoint except at the endpoints of each copy) comprises the whole Koch curve.  This is, in fact, the essence of self-similarity. By abuse of language, we say that the Koch snowflake curve itself is ``self-similar''; see Figure \ref{fig:3kochcurvesLabeled}.}
	\label{fig:kochcurveconstruction}
\end{figure}

\begin{figure}
	\centering
		\includegraphics{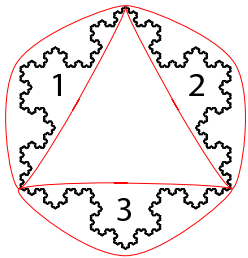}
	\caption{The self-similarity of the Koch snowflake curve $KS$.  Shown here is the Koch snowflake curve $KS$, viewed as the union of three isometric, abutting copies of the Koch curve.}
	\label{fig:3kochcurvesLabeled}
\end{figure}

The Koch snowflake curve, as depicted in Figure \ref{fig:kochconstruction75}, is a fractal.  In particular, it is the union of three \textit{self-similar} Koch curves, with the Koch curve being a continuous, nowhere differentiable curve with infinite length (see Figures \ref{fig:kochcurveconstruction} and \ref{fig:3kochcurvesLabeled}).  Consequently, any attempt to construct a line tangent to the Koch snowflake curve may seem like an exercise in futility.  This poses a unique problem for defining the trajectory of a \textit{billiard ball} \label{idx:snowflakebilliard}(i.e., a pointmass traversing the interior of the planar region bounded by the Koch snowflake $KS$).  Specifically, when this pointmass collides with the boundary $KS$ with unit speed, the absence of a well-defined tangent results in multiple choices for the angle of reflection, meaning there is a priori no well-defined angle of reflection. In \cite{LaNie1}, we provided experimental evidence in support of the existence of certain periodic orbits of the Koch snowflake billiard $\omegaks$.  Moreover, in \cite{LaNie1}, we stated several conjectures about the existence of a well-defined billiard $\omegaks$ and the dynamical equivalence between the conjectured billiard flow on $\omegaks$ and the associated geodesic flow on the proposed corresponding `\textit{fractal flat surface}'.  One of the main objectives of the present paper is to investigate what one means by \textit{reflection} in the snowflake boundary and to establish the existence and describe the topological and geometric properties of particular families of periodic orbits (and/or of their footprints) of the Koch snowflake billiard $\omegaks$.

In short, the point of view adopted here is to define certain ``periodic orbits'' of $\omegaks$ as suitable (inverse) limits of certain ``compatible sequences'' of periodic orbits of its (inner) rational polygonal billiard approximations $\omegaksi{n}$.  Using this definition and a study (conducted in \S3) of the periodic orbits of the $n$th prefractal billiard approximation \label{idx:omegaksi}$\omegaksi{n}$, for each fixed $n\geq 0$, we characterize and describe (in terms of the ternary expansion of their initial basepoint) the periodic orbits with an initial direction making an angle of $\pi/3$ with the horizontal in $\omegaksi{0}$.

More specifically, we are able to construct what we call the \textit{footprint} of the \textit{primary piecewise Fagnano} and \textit{piecewise Fagnano} orbits, and \textit{stabilizing periodic orbits} of $\omegaks$ (all in the direction of $\pi/3$).  As of now, the only family of well-defined orbits of the Koch snowflake billiard $\omegaks$ is the family consisting of what we call \textit{stabilizing orbits}. That which we propose to be a piecewise Fagnano orbit of $\omegaks$ has a footprint $\mathcal{F}(\xoo)$ that is the inverse limit of footprints (of piecewise Fagnano orbits) of the prefractal approximations.  A footprint of a prefractal approximation amounts to the Poincar\'e section of the billiard map describing the billiard flow on the corresponding phase space.  While we say ``piecewise Fagnano orbit,'' we are making an abuse of language in that we do not mean to imply that an orbit actually exists, but that whatever the orbit truly \textit{is}, it has a footprint $\mathcal{F}(\xoo)$.  Furthermore, even less is known about what we have called the \textit{generalized piecewise Fagnano orbits}.  Again, there really is no orbit to speak of, nor is there any footprint to speak of.  We discuss all of these `orbits' in \S 4 and \S5 with differing degrees of rigor, sometimes only providing a plausibility argument as to why a given orbit \textit{should} have a particular property.

This paper draws upon various subjects in mathematics.  So as to accommodate a diverse audience of readers, we make a considerable effort in developing the necessary background material in \S2.  In particular, we give a thorough description of the billiard flow associated with a billiard table $\Omega(B)$ with (piecewise) smooth boundary $B$. We also recall the notion of a \textit{flat surface} and how one can construct a flat surface from a \textit{rational polygonal billiard table} (that is, a planar billiard table whose boundary is a polygon with interior angles that are rational multiples of $\pi$).  In this context, a flat surface is a mathematical device used to rigorously describe the billiard flow on $\Omega(B)$ in terms of the geodesic flow on the surface. In addition, we recall the definition of \textit{inverse limit} and explain how the Cantor set $\cantor$ can be viewed as the inverse limit of an \textit{inverse limit sequence} of its prefractal approximations, denoted by $\cantor_n$, with the index $n$ corresponding to the approximation with $2^n$ many points. %segments of length $1/3^n$.

Similarly, the snowflake curve $KS$ can be viewed as the inverse limit of its prefractal approximations $KS_n$.  Here, $KS_n$ is the $n$th (inner) polygonal approximation to $KS$, and hence defines a rational polygonal billiard table $\omegaksi{n}$.  Since the theory of rational polygonal billiards is very well developed (see, e.g., [GaStVo,Gu1,GuJu1--2,HuSc,KaHa2,KaZe,Mas,MasTa,Ve1--3,Vo,Zo]), it is then natural to define the dynamics on the fractal ``billiard table'' $\omegaks$ in terms of the dynamics on its prefractal approximations $\omegaksi{n}$.  As a result, much of the focus of this paper will be to first obtain a good understanding of the periodic orbits of $\omegaksi{n}$, and then to provide a plausibility argument as to why we can view piecewise Fagnano orbits on $\omegaks$ as suitable (inverse) limits of appropriate sequences of piecewise Fagnano orbits on $\omegaksi{n}$, for $n\geq 0$.

It is our intention that those familiar with the topic of mathematical billiards and not fractal geometry find readily accessible the basic notions of self-similarity and that of an iterated function system (IFS).  We briefly describe these notions by means of a simple example in \S\ref{subsec:compSeqPF} and refer the reader to various references for further details (e.g., \cite{Ba,Ed,Fa}).  So as to accommodate readers from the physical sciences, we also attempt to explain the necessary concepts from topology and geometry.  As such, the interested reader will find references to \cite{Ma} for further details on \textit{covering spaces}, \cite{McL} for category theory, \cite{Bo} for general topology and \cite{HoYo} for the specialized topic of inverse and direct limits in the context of the category of topological spaces.

In \S3, we begin the discussion of our results.  We prove that the \textit{prefractal flat surface} $\mathcal{S}(KS_n)$ associated with the Koch snowflake prefractal billiard $\omegaksi{n}$ is a \textit{branched cover}\footnote{We briefly discuss the definition of \textit{branched cover} in \S3.} of the \textit{hexagonal torus} $\mathcal{S}(KS_0)$.  We then use this result to show that initial directions of periodic orbits in the billiard $\omegaksi{0}$ are exactly the initial directions of periodic orbits in  the billiard $\omegaksi{n}$, and vice-versa.

This key fact serves as the foundation for \S4.  It is there that we develop much of the machinery necessary to describe what we call \textit{compatible sequences of periodic orbits}.  Further discussion about the nature of particular points in the unit interval $I$ \label{idx:unitInterval}(which we always view as the base of the equilateral triangle $KS_0 := \Delta$) gives rise to specific \textit{compatible sequences}.  We note that \S4 is very dense, serving as a strong foundation for \S5.  While the addressing system used  in \S5 (and introduced and used in \cite{LaPa}) may seem to make some of the tools developed in \S4 redundant, in principle, many of the proofs of the results in \S5 demonstrate the interconnectedness of the two sections, making the implicit (and even explicit) dependence of \S5 on \S4 readily apparent. 

Therefore, we model the structure of \S5 on that of \S4, so as to allow the reader to draw parallels more quickly and to see how the results in \S4 influence our later study of the periodic orbits of the Koch snowflake billiard.  In \S 4, we show that directions for which an orbit $\orbitksi{n}$ of $\omegaksi{n}$ are periodic are exactly the same for which an orbit $\orbitksi{0}$ of $\omegaksi{0}$ is periodic.  This aids us in constructing what we call a  \textit{compatible sequence of closed orbits}.  We focus our investigation on orbits with an initial direction of $\pi/3$.  As such, we describe what we call a \textit{compatible sequence of piecewise Fagnano orbits}, an \textit{eventually constant compatible sequence of orbits} and \textit{compatible sequence of generalized piecewise Fagnano orbits}.  The period and length of \textit{piecewise Fagnano orbits}, $\cantor$\textit{-orbits} and \textit{generalized piecewise Fagnano orbits} of $\omegaksi{n}$ are given in terms of the ternary representation of the initial basepoint of the initial orbit of the respective compatible sequence of periodic orbits.  In \S \ref{subsubsec:topologicalAndGeometricPropertiesOfFootprint} and \S\ref{subsec:stabilizingOrbitsKS}, we describe the topological and geometric properties of what we call the \textit{footprint} of a piecewise Fagnano orbit and \textit{stabilizing orbits} (or \textit{$\cantor$-orbits}), respectively.  In \S\ref{subsec:plausibilityArgument}, we then provide a plausibility argument for the existence of what we call a \textit{piecewise Fagnano orbit} of the Koch snowflake billiard $\omegaks$.  Finally, we close \S5 by conjecturing the existence of \textit{generalized piecewise Fagnano orbits} of the Koch snowflake billiard $\omegaks$.

%is a \textit{periodic orbit} of the Koch snowflake billiard $\omegaks$.  In particular, we argue that it is reasonable to expect that the inverse limit of a compatible sequence of piecewise Fagnano orbits is indeed an orbit of the Koch snowflake billiard $\omegaks$.  Moreover, we argue that it should then be possible to describe such orbits entirely in terms of the ternary expansion of the initial basepoint $\xoo\in I$ of a particular orbit $\orbitksi{0}$ of the equilateral triangle billiard $\omegaksi{0}$ (in fact, the first orbit in what we will call a \textit{compatible sequence of periodic orbits}).  This gives rise to a description of periodic orbits of the Koch snowflake billiard $\omegaks$, where we stipulate that this description assumes a fixed initial direction of $\pi/3$ for every orbit being considered.  We will show that such orbits have finite length.  Moreover, a certain subcollection of these periodic orbits, which is comprised of what we call \textit{piecewise Fagnano orbits} (see Figure \ref{fig:fagnanos012}), that 1) have footprints that, collectively, constitute a subset of the \textit{natural Cantor set} of the Koch snowflake $KS$ and 2) the footprint of an individual orbit from this particular subcollection of orbits constitutes a topological Cantor set in its own right and is a self-similar fractal.

Considering the fact that the field of ``fractal billiards'' is still in its infancy, we provide many open questions and conjectures in \S6.  In particular, we stress that $\omegaks$ does not constitute a well-defined mathematical billiard, in the sense that we have not provided a well-defined phase space, let alone a geodesic flow on such a phase space.  Such a mathematical object has yet to be precisely defined, but the work we have laid out in \S5 and the remarks made in \S6 indicate a possible path for constructing such a phase space and geodesic flow.  

In addition to determining the nature of such a geodesic flow and whether it could be dynamically equivalent to the billiard flow, we ask questions regarding the \textit{ergodic} nature of the conjectured geodesic and billiard flows on the hypothesized `fractal flat surface' and the corresponding fractal billiard.  We then state an open question asking whether or not it is possible to construct analogs of various classical trace formulas (e.g., the Gutzwiller trace formula) which can be used, in the classical billiard case, as a tool for connecting the length spectrum of particular billiards and the eigenvalue spectrum of the Laplacian defined on their associated fractal drums.

In the long-term, we hope that the present (preliminary) study of the Koch snowflake billiard will help lay the foundations for a general theory of fractal billiards.

% Let $X$ be a topological space. An \textit{Iterated Function System} is a collection of maps $\{\phi_i\}^{n}_{i=1}$ where $\Phi:=\bigcup_i^n\phi_i$ is a map defined on the space of all compact subsets of $X$.  If $\phi_i$ is a contraction mapping, then so is $\Phi$.  Hutchinson's theorem then states that there is a unique fixed point attractor of this IFS.  We then call such an attractor a \textit{self-similar} set.  Specifically, a self-similar set $S = \Phi(S) = \bigcup_{i=1}^n \phi_i(S)$. 

\subsection{Index of notation}

At this point, we would like to provide the reader with an index of notation.  Such an index will list the notation, a brief explanation of the notation and where this notation is first used or defined.

\begin{figure}[h!]
\begin{center}
\begin{tabular}{|c|p{9 cm} c|}%{|l|lc|}
\hline
$I$, $[0,1]$ & The unit interval.  By convention,  $I$ is the base of the equilateral triangle $\Delta=KS_0$, and is identified with $[0,1]$ & \pageref{idx:unitInterval}\\
$KS_n$ & The $n$th Koch snowflake prefractal curve, $n\geq 0$ &  	\pageref{fig:kochconstruction75}\\
$KS$ & The Koch snowflake fractal curve & \pageref{fig:kochconstruction75} \\
$\omegaksi{n}$ & The $n$th Koch snowflake prefractal rational billiard & \pageref{idx:omegaksi}\\
$\omegaks$ & The Koch snowflake fractal billiard & \pageref{idx:snowflakebilliard}\\
$f_B$ & The billiard map associated with the billiard $\Omega(B)$& \pageref{idx:billiardMapB}\\
$f_n$ & The billiard map associated with the billiard $\omegaksi{n}$& \pageref{idx:billiardMapKSn}\\
$\mathcal{S}(R)$ & The flat surface corresponding to the rational billiard $\Omega(R)$ &\pageref{idx:rationalFlatSurface}\\
$\cantor$ & The ternary Cantor set& \pageref{exp:ternaryCantorSet}\\
$\mathcal{S}(KS_n)$ & The $n$th Koch snowflake prefractal flat surface &\pageref{idx:SKSnFlatSurface}\\
$s_{n,k}$ & A side of $\omegaksi{n}$& \pageref{rmk:SnkXnkThetank}\\
$\xio{n}$ & An initial basepoint of an orbit of $\omegaksi{n}$& \pageref{rmk:SnkXnkThetank}\\
$\xii{n}{k_n}$ & A basepoint of an orbit of $\omegaksi{n}$ & \pageref{rmk:SnkXnkThetank}\\
$\orbitksi{n}$ & An orbit of $\omegaksi{n}$ with initial condition $(\xio{n}, \theta_{n}^{0})$ & \pageref{rmk:SnkXnkThetank}\\
$g_{n,k}$ & A ghost of a side $s_{n,k}$ of $\omegaksi{n}$ ($1\leq k\leq 3\cdot 4^n$)& \pageref{def:ghost}\\
$C_{n,k}$ & A cell of $\omegaksi{n}$ corresponding to a side $s_{n-1,k}$ of $\omegaksi{n-1}$ & \pageref{def:cell}\\
$\pfix{n}{\yio{n}}$ & A piecewise Fagnano orbit of $\omegaksi{n}$&\pageref{rmk:labelingOfpfis}\\
$M(\cantor)$ & A specific collection of points in the unit interval $I$  &\pageref{eqn:MCantor}\\
%$\gpfix{n}{\yio{n}}$ & A generalized piecewise Fagnano orbit of $\omegaksi{n}$ &$\cdots$&  \\
$0.u_1u_2..._3$ & A base-3 expansion of a number in the unit interval $I$ & \pageref{idx:baseThree}\\
$\cantor_{n,k}$; $\cantor_{n,k}'$ &The ternary Cantor set of a side $s_{n,k}$ of $\omegaksi{n}$; the non-ternary points of $\cantor_{n,k}$  & \pageref{idx:cnk}\\
$l,c,r$ & An alternate alphabet used to give a ternary representation of an element of $I$ & \pageref{idx:lcr}\\
\hline
\end{tabular}
\end{center}
\end{figure}

\newpage

\setcounter{tocdepth}{2}
\tableofcontents

\section{Background}

\subsection{Mathematical Billiards and the Billiard Map}
\label{subsec:mathBillandBillMap}

Under ideal conditions, we know that a point mass having a perfect elastic collision with a $C^1$ surface (or curve) will reflect at an angle which is equal to the incoming angle, both measured relative to the normal at the point of collision.

%The Law of Reflection states that when a pointmass collides with a sufficiently smooth boundary, the angle of reflection will be equal to the angle of incidence.  

Consider a compact region $\Omega(B)$ in the plane with connected boundary $B$.  Then, $\Omega(B)$ is called a \textit{planar billiard} when $B$ is smooth enough to allow the Law of Reflection to hold, off of a set of measure zero.  When $B$ is a nontrivial connected polygon in $\mathbb{R}^2$, $\Omega(B)$ is called a \textit{polygonal billiard}, and the collection of vertices forms a finite set, which is a set of zero measure (when we take our measure to be the Hausdorff measure or simply, the arc-length measure on $\Omega(B)$). For the reader's easy reference, we provide a formal definition of \textit{rational billiard} below.

\begin{definition}[Rational polygon and rational billiard]
If $B$ is a nontrivial connected polygon such that for each interior angle $\theta_i$ of $B$ there are relatively prime integers $p_i > 0$ and $q_i>0$ such that $\theta_i = \frac{p_i}{q_i} \pi$, then we call $B$ a \textit{rational polygon} and $\Omega(B)$ a \textit{rational billiard}.
\label{def:ratBilliard}
\end{definition}

\begin{figure}
\begin{center}
\includegraphics[scale = .5]{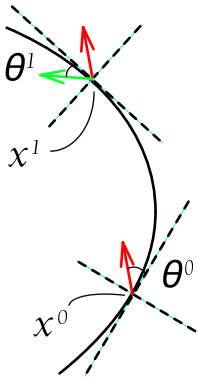}
\end{center}
\caption{A billiard ball traverses the interior of a billiard and collides with the boundary.  The velocity vector is pointed outward at the point of collision. The presence of a well-defined tangent at this point provides for the existence of a normal to the tangent and the recovery of the Law of Reflection (i.e., Snell's Law).  The resulting direction of flow is found by either reflecting the vector through the tangent or by reflecting the incidence vector through the normal and reversing the direction of the vector.  We use the former method throughout this paper. A rigorous discussion of the Law of Reflection in this context is given in \cite{Sm}.}
%{The billiard ball hits the boundary and the vector describing the position and velocity is not pointing outwards.  Since the boundary is sufficiently smooth in a neighborhood of this point, we know that there exists a well-defined tangent.  We reflect this vector through the tangent to indicate the direction (velocity) of the billiard ball after collision in the boundary.  Moreover, by reflecting the billiard ball in the boundary, we form an equivalence class of vectors.} 
\label{fig:billiardMap}
\end{figure}

The Law of Reflection essentially amounts to reflecting the incoming vector through the normal vector and then reversing the direction. 	Instead, we adhere to the convention in the field of mathematical billiards according to which the vector describing the position and velocity of the billiard ball (which amounts to the position and angle, since we are assuming unit speed) be reflected in the tangent to the point of incidence.  Then we can rigorously reformulate the Law of Reflection as follows: the vector describing the motion is the reflection of the incoming vector through the tangent at the point of collision.  Moreover, we can identify these two vectors and form an equivalence class of vectors in the unit tangent bundle corresponding to the billiard table $\Omega(B)$.  (See Figure \ref{fig:billiardMap} and \cite{Sm} for a detailed discussion of this equivalence relation on the unit tangent bundle $\Omega(B)\times S^1$.)

For the remainder of this section, we suppose that $\Omega(B)$ is a polygonal billiard with connected boundary.\footnote{In fact, in the remainder of the paper, except when $B$ is the Koch snowflake, $B$ will be a nontrivial connected polygon.  As such, the billiard flow given by the billiard map $f_B:(\Omega(B)\times S^1)/\sim \,\,\to (\Omega(B)\times S^1)/\sim $ will have finitely many singularities, whereas $\Omega(KS)$ has, a priori, uncountably many singularities.} Denote by $S^1$ the unit circle, which we may consider to represent all the possible directions (or angles) in which a billiard ball may initially move.  In practice, one restricts one's attention to $(B\times S^1)/\sim$ when discussing the phase space of the billiard $\Omega(B)$.  To clearly understand how one forms equivalence classes from elements of $B\times S^1$, we let $(x,\theta),(y,\gamma)\in B\times S^1$ and say that $(x,\theta)\sim (y,\gamma)$ if and only if $x=y$ and one of the following is true:

\begin{enumerate}
\item{$x=y$ is not a vertex of the boundary $B$ and $\theta = \gamma$,}
\item{$x=y$ is not a vertex of the boundary $B$, but $x=y$ is a point on a segment $s_i$ of the polygon $B$ and $\theta = r_i(\gamma)$, where $r_i$ denotes reflection in the segment $s_i$,}
\item{If $x=y$ is a vertex of $B$, then we identify $(x,\theta)$ with $(y,g(\gamma))$ for every $g$ in the group generated by reflections in the two adjacent sides having $x$ (or $y$) as a common vertex.}
\end{enumerate}

\label{idx:billiardMapB} The \textit{billiard map} (sometimes also called the \textit{Poincar\'e map}) is denoted by $f_B$ and is determined in the following way.  Let $[(x^0,\theta^0)] \in (B\times S^1)/\sim$ be the equivalence class of $(x^0,\theta^0)$, with $x^0\in B$ and $\theta^0\in S^1$ an inward pointing direction (unit vector) based at $x^0$.  Then $f_B([x^0,\theta^0])$ identifies the first point of collision in the boundary and the outward pointing vector $\widetilde{\theta^0}$ that is parallel with the unit vector in the direction $\theta^0$.  Since $f_B:(B\times S^1)/\sim\,\,\longrightarrow (B\times S^1)/\sim$, we have that $(x^1,\widetilde{\theta}^0)\in [(x^1,\theta^1)]=f_B([x^0,\theta^0])$.  As expected, the outward pointing vector is identified with the inward pointing vector and one then takes as the direction of travel the inward pointing vector $\theta^1$ based at $x^1$.\footnote{The representative element is always the one for which the vector describing the direction after collision is inward pointing.}   

The \textit{billiard map} $f_B$  describes the \textit{discrete billiard flow} on the phase space $(\Omega(B)\times S^1)/\sim$.  It is related as follows to the \textit{continuous billiard flow} $f_B^t: \mathbb{R}\times (\Omega(B)\times S^1)/\sim\,\,\to (\Omega(B)\times S^1)/\sim$.  When we restrict our attention to $B\times S^1$,  the flow $f_B^t$ on $(B\times S^1)/\sim$ is discrete.  Specifically, one denotes the time at which the billiard ball collides with the boundary $B$ by $t_i$, $i\in \mathbb{Z}$.  The discrete billiard flow is then given by $f_B^{t_i}:\mathbb{Z}\times (B\times S^1)/\sim\,\,\to (B\times S^1)/\sim$.  Such a flow is also called the \textit{Poincar\'e section} of the continuous flow $f_B^t$; see \cite{KaHa1} for a more complete discussion of the billiard flow.  In the sequel, instead of denoting the discrete flow by $f_B^{t_i}$, we use the standard notation $f_B^i$.  In general, for every integer $k\geq 0$, we have $[x^k,\theta^k] = f_B^k([x^0,\theta^0])$, where $f^k_B$ (defined above) can also be viewed as the $k$th iterate of the billiard map $f_B$.

\begin{remark}
\label{idx:billiardMapKSn}
So as not to introduce cumbersome notation, when we begin discussing the billiard map $f_{KS_n}$ corresponding to the $n$th prefractal billiard $\omegaksi{n}$, we will simply write $f_{KS_n}$ as $f_n$.  Moreover, when discussing the billiard flow on $(\omegaksi{n}\times S^1)/\sim$, the $k$th point in an orbit $[(x^k,\theta^k)]\in (\omegaksi{n}\times S^1)/\sim$ will instead be denoted by $[(\xii{n}{k_n},\theta_n^{k_n})]$, so as to be clear as to which space such a point belongs.
\end{remark}

%The \textit{billiard map} $f:\Omega(B)\times S^1/\sim \,\,\rightarrow \Omega(B)\times S^1/\sim$ maps $[(x,\theta)]\in \Omega(B)\times S^1/\sim$ to the appropriate element $[(x',\theta')]$, where 1) $\theta'$ is inward pointing and 2) $(x',\theta')$ we take to be the representative element. (Here, for example $[(x,\theta)]$ denotes the equivalence class of $(x,\theta)\in \Omega(B)\times S^1$.)

In the event that a basepoint $x^1$ of $f_B(x^0,\theta^0)$ is a corner of $\Omega(B)$ (that is, a vertex of the polygonal boundary $B$) and $\theta^0$ was a direction for which the billiard flow would be periodic, then the resulting closed orbit is said to be \textit{singular}.  In addition, since $\theta^0$ is a direction for which the resulting orbit is periodic, there exists a positive integer $k$ such that the basepoint $x^{-k}$ of $f_B^{-k}(x^0,\theta^0)$ is a corner of $\Omega(B)$. (Here, $f_B^{-k}$ denotes the $k$th inverse iterate of $f_B$.)  The path then traced out by the billiard ball connecting $x^1$ and $x^{-k}$ is called a \textit{saddle connection}.  

\begin{remark}
Within the subject of mathematical billiards, there appears to be a slight abuse of language.  One may refer to the orbit of a billiard as the path traced out by the billiard ball or as the collection of incidence points in the boundary.  In the latter case, such a set of points is referred to as the Poincar\'e section of the billiard map $f_B$.  When we want to be clear as to which concept we are referring, we will specifically write `the path corresponding to the orbit' or `the footprint (or Poincar\'e section) of the corresponding orbit', respectively.
\end{remark}

\begin{figure}
\begin{center}
\includegraphics{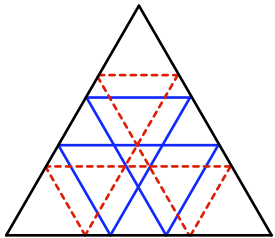}
\end{center}
\caption{These two orbits are equivalent.  The initial angle is $\pi/3$ and the segments in the path remain parallel after each subsequent reflection.  It should be noted, however, that an orbit with an initial basepoint at the midpoint of the base of the equilateral triangle and an initial angle of $\pi/3$ results in what is called the \textit{Fagnano orbit}, which has period $3$ and a path with a length that is half of the length of either orbit shown in this figure. This does not mean that we are making an exception and saying that this Fagnano orbit \textit{should} be equivalent to the other two shown, but merely that one must be careful when making the claim that all orbits with the same direction are necessarily equivalent.}
\label{fig:equivalentOrbits}
\end{figure}

When considering \textit{periodic orbits} (which are, by definition, non-singular closed orbits), one point of interest is the corresponding path traced out by the billiard ball and the length of such a path.  Such a path necessarily has finite length since the orbit has finite period.  If two periodic orbits $\mathscr{O}(x,\theta)$ and $\mathscr{O}(y,\phi)$ have corresponding paths $\mathscr{P}(x,\theta)$ and $\mathscr{P}(y,\phi)$ that have equal lengths, equal periods and, $x$ and $y$ lie on the same segment of $B$ such that $\mathscr{P}(x,\theta)$ and $\mathscr{P}(y,\phi)$ remain parallel from one basepoint to the next, then we say that these two orbits are equivalent.  From this, we see that equivalence of orbits does not depend at all on the choice of the initial basepoint of the given periodic orbits;\footnote{Except for a few special cases, one of which is discussed in the caption of Figure \ref{fig:equivalentOrbits}.} see Figure \ref{fig:equivalentOrbits}.

\subsection{Flat Structures and Flat Surfaces}
\label{subsec:flatStructuresandFlatSurfaces}
In this section, we deal only with flat surfaces constructed from rational billiards; see Definition \ref{def:ratBilliard} in \S\ref{subsec:mathBillandBillMap} above.  By showing that there is an intimate connection between what is called a flat surface and a billiard $\Omega(B)$, we will eventually be able to demonstrate that reflection in certain corners can be defined and, for certain rational billiards, the nature of an orbit of a special class of rational billiards does not depend at all on the initial point, but only on the initial direction.\footnote{In \S\ref{sec:openProblemsConj}, we discuss what are called \textit{uniquely ergodic flows} and the \textit{Veech dichotomy}.  Understanding such concepts is not necessary for \S\S3--5; so we defer explanation of those concepts until \S\ref{sec:openProblemsConj} and refer the reader to [\textbf{Ve1--3}] for an immediate and detailed description of the Veech dichotomy and to \cite{MasTa,HuSc,Zo} for a pedagogical discussion of uniquely ergodic, periodic and aperiodic billiard flows.}   To such end, we discuss the notion of a flat structure.

A flat structure on a connected surface $M$ is an atlas $(U_\alpha,\phi_\alpha)$ with a finite number of singularities such that, away from these singularities, each `coordinate changing' map

\begin{equation}
\phi_\alpha\circ \phi^{-1}_\beta : \phi_\beta(U_\alpha\cap U_\beta) \to \phi_\alpha(U_\alpha\cap U_\beta)
\end{equation}

\noindent is a translation in $\mathbb{R}^2$.  Specifically, a flat structure satisfies the following definition.

\begin{definition}[Flat structure]
\label{def:flatStructure}
Let $M$ be a compact, connected, orientable surface.  A \textit{flat structure} on $M$ is an atlas $\omega$, consisting of charts of the form $(U_\alpha,\phi_\alpha)$, where $U_\alpha$ is a domain (i.e., a connected open set) in $M$ and $\phi_\alpha$ is a homeomorphism from $U_\alpha$ to a domain in $\mathbb{R}^2$, such that the following conditions hold:

\begin{enumerate}
\item{the domains $U_\alpha$ cover the whole surface $M$ except for finitely many points $z_1,z_2,...,z_k$, called \textit{singular points};}
\item{all coordinate changing functions are shifts (i.e., translations) in $\mathbb{R}^2$;}
\item{the atlas $\omega$ is maximal with respect to properties $(1)$ and $(2)$;}
\item{for each singular point $z_i$, there is a positive integer $m_i$, a punctured neighborhood $\dot{U}_i$ not containing other singular points, and a map $\psi_i$ from this neighborhood to a punctured neighborhood $\dot{V}_i$ of a point in $\mathbb{R}^2$ that is a shift in the local coordinates from $\omega$, and is such that each point in $\dot{V}_i$ has exactly $m_i$ preimages under $\psi_i$.}
\end{enumerate}

\end{definition}

\begin{definition}[Flat surface]
\label{def:flatSurface}
We say that a connected, compact surface equipped with a flat structure is a \textit{flat surface}.  
\end{definition}

\begin{remark}
Calling a connected, compact ($2$-dimensional) manifold with flat structure a flat surface is somewhat of an abuse of language, but it enables us to be briefer and refer to related notions with greater ease.  Note that in the literature on billiards and dynamical systems, the terminology and definitions pertaining to this topic are not completely uniform; see, for example, [\textbf{GaStVo,Gu1,GuJu1--2,HuSc,KaHa2,Mas,MasTa, Ve1--3,Vo,Zo}].  We have adopted the above definition for clarity and the reader's convenience.    
\end{remark}

The singular points of a flat surface $\mathcal{S}$ are called \textit{singularities of the flow}.  There are two types of singularities in a flat surface: removable and nonremovable.  They are called such because it may or may not be possible to define the flow at these points.  We are primarily interested in how the geodesic flow behaves at these points.  We now turn to a discussion of why these singularities can be termed ``removable'' or ``nonremovable'' and how to discern between the two types.  Consider $z_i$ in the set of singularities of a flat structure.  This point has what is called a conic angle.  A \textit{conic angle} is the number of radians required to form a closed circle about $z_i$.  

We state the following definitions, for completeness.

\begin{definition}[Removable conic singularity]
A singularity of a flat structure $(U_\alpha,\phi_\alpha)$ is a \textit{removable singularity} of the flow if the conic angle about such a point is $2\pi$.  
\end{definition}

\begin{definition}[Nonremovable conic singularity]
A singularity of a flat structure $(U_\alpha,\phi_\alpha)$ is a \textit{nonremovable singularity} of the flow if the conic angle about such a point is $2\pi c$, for some integer $c>1$.  
\end{definition}

With regards to the geodesic flow on a flat surface $\mathcal{S}$, removable singularities of the flow pose no problem and the flow lines through such points can continue unimpeded.  However, when a singularity is nonremovable, the flow cannot continue unimpeded.  Moreover, when $z_i$ is a removable singularity of a flat structure $\omega$, there exists a flat structure $\tilde{\omega}$ on $\mathcal{S}$ such that $z_i$ is not a singular point of $\tilde{\omega}$.  Technically speaking, when $z_i$ is a removable singularity of the flow on a surface with a flat structure, there exist neighborhoods $U_\beta$ and $U_\gamma$ of $z_i$ and maps $\phi_\beta$ and $\phi_\gamma$ such that the associated transition map $\phi_\beta \circ \phi_\gamma^{-1}$ is a shift in the local coordinates at $z_i$.  Consequently, the geodesic flow may be continuously extended at removable singularities of the flat surface.

\label{idx:rationalFlatSurface}We now discuss how to construct a flat surface from a rational billiard.  Consider a rational polygon billiard $\Omega(P)$ with $m$ sides and interior angles $\pi(u_i/v_i)$ at each vertex $z_i$, for $1\leq i\leq m$. Here, $u_i$ and $v_i$ are relatively prime positive integers. Then, a (laborious) calculation shows that for some $j\leq m$, $v_j=\text{lcm} (v_i)_{i=1,i\neq j}^m$. Consequently, the linear portions of the planar symmetries generated by reflection in the sides of the polygonal billiard $\Omega(P)$ generate the dihedral group $D_N$, where $N:=\text{lcm} \{v_i\}_{i=1}^m$.  Here, by definition, the \textit{dihedral group $D_N$} denotes the group of symmetries of the regular $N$-gon.  So as to be clear, we mention that the notation $D_N$ does not refer to the (wrong) fact that such a finite group has $N$ elements.\footnote{Actually, $D_N$ has $2N$ elements, and the standard group theory notation for the dihedral group is then $D_{2N}$, since the cardinality of the group is often given more importance, from the perspective of group theory.}  We next consider $\Omega(P)\times D_N$ (equipped with the product topology).  We want to glue `sides' of $\Omega(P)\times  D_N$ together and construct a natural atlas on the resulting surface $M$ so that $M$ becomes a flat surface.  

To such end, let $p_1=p_2$ be a point on a side $s_a$ of $\Omega(P)$, $r_a$ be the linear portion of the reflection determined by reflecting $\Omega(P)$ in the side $a$, and $(p_1,r_1),(p_2,r_2)\in P\times D_N$.  Then, by definition, $(p_1,r_1)\sim (p_2,r_2)$ if and only if

\begin{enumerate} 
\item{$(p_1,r_1) = (p_2,r_2)$, or}
\item{$r_a = r_1^{-1}r_2$, or}
\item{$p_1$ and $p_2$ are the same vertex of $\Omega(P)$ having adjacent sides $s_a$ and $s_b$, with $r_1^{-1}r_2$ belonging to the subgroup generated by $r_a$ and $r_b$.}
\end{enumerate}

It is easy to check that $\sim$ is an equivalence relation on $\Omega(P)\times D_N$.  More work is required to show that $M:= (\Omega(P)\times D_N)/\sim$ is a compact, connected, orientable surface.  As a result of the identification, the points of $M$ that correspond to the vertices of $\Omega(P)$ constitute (removable or  nonremovable) conic singularities of this surface.  Heuristically, $\Omega(P)\times D_N$ can be represented as $\{r_iP\}_{i=1}^{2N}$, in which case it is easy to see what points are made equivalent under the action of $\sim$.

Denote by $\Omega(P)^\circ$ the interior of the billiard $\Omega(P)$.  To construct a flat structure on $M$, let $U_i=\Omega(P)^\circ\times \{r_i\}$.  Then $\{U_i,\phi_i\}_{i=1}^{2N}$ can be naturally extended to constitute a flat structure on $M$, in the sense of Definition \ref{def:flatStructure}.  Hence, $M$ is a flat surface, in the sense of Definition \ref{def:flatSurface}. The map $\phi_i\circ\phi_j^{-1} : \phi_j(U_i\cap U_j)\to \phi_i(U_i\cap U_j)$ is a translation in the local coordinates of a point $z\in U_i\cap U_j$; i.e., $\phi_i\circ\phi_j^{-1}(z) = z+d$, where the constant $d$ is independent of the choice of $i$ and $j$.\footnote{A priori, the choice of $d$ describing the translation in the local coordinates \textit{does} depend on $i$ and $j$.  However, given the fact that one constructs the flat surface by identifying parallel and opposite sides of the polygon $B$, for a fixed direction $\theta$, one can describe a parallel line field in the direction $\theta := \text{arctan} \frac{\alpha}{\beta}$, with $d =:\alpha+\sqrt{-1}\beta$.}

\begin{example}[The equilateral triangle billiard $\omegaksi{0}$]
We consider the equilateral triangle $\Delta:=KS_0$,  an important example in the literature (see, e.g., \cite{BaxUm}) and an even more important example when considering the Koch snowflake prefractal billiard. Since the interior angles $\pi p_i/q_i$ are all the same and equal to $\pi/3$, we have

\begin{equation}
N=\text{\emph{lcm}}(q_i)_{i=1}^3 = \text{\emph{lcm}}(3,3,3) = 3.
\end{equation}

\noindent Consequently, the associated surface is given by $\mathcal{S}(KS_0):=(\Delta \times  D_3)/\sim$.  Moreover, there is a flat structure on this surface and since all the singularities are removable conic singularities, such a structure can be extended to include the singular points and the resulting surface is topologically equivalent to a torus. (Really, the associated flat surface $\mathcal{S}(KS_0)$ is the hexagonal torus, which is topologically equivalent to the standard torus; see Figure \ref{fig:hexagonalTorus} in \S3.)  
\end{example}

We close this discussion by recalling an important fact about the geodesic flow on a flat surface $\mathcal{S}(B)$ constructed from a rational billiard $\Omega(B)$ and the billiard flow on that rational billiard.  These two flows can be shown to be (dynamically) equivalent under the action of the group $D_N$ associated with the construction of the corresponding flat surface.  Heuristically, one may view the corresponding equivalence as follows: a billiard flow line may be straightened to a geodesic flow line and additionally, a geodesic flow line may be collapsed into a billiard flow line.  In order to be more technically correct, we further explain the details of the equivalence of these two flows.  If we consider the geodesic flow on the surface, the quotient of the phase space by the group of symmetries associated with the construction of $\mathcal{S}(B)$ has the effect of collapsing the space to a space with a quotient flow that is isomorphic to the billiard flow.  On the other hand, any given billiard orbit may be straightened by making successive reflections, via the action of $D_N$, in the identified sides of the flat surface, therefore producing a straight-line flow line on the flat surface $\mathcal{S}(B)$.

\subsection{Inverse limit sequence and inverse limit}
\label{subsec:invLimSeqInvLim}

Inverse limits of various topological (algebraic, or geometric) objects will play an important role in this paper.  Hence, since such a notion may not be familiar to all readers, it may be helpful to recall some basic facts pertaining to this subject.  Further information can be found, for example, in \cite{HoYo} and, in a more general context, in \cite{Bo,McL}. 

We discuss the inverse limit in the context of the category Set.  The objects are sets and the morphisms are set maps.\footnote{Discussing the inverse limit in the context of Set is purely a formality that allows us to speak in more concrete terms and utilize important existences and uniqueness properties of the inverse limit in the category Set. See \cite{Bo,McL} for a general discussion.} Consider a partially ordered collection of sets $\{X_i\}_{i=1}^\infty$, with each $X_i$ equipped (for each integer $n\geq 1$) with a map $F_n:\prod_{i=1}^\infty X_i \to X_n$ defined by $F_n((x_i)_{i=1}^\infty)=x_n$. If for every $n$ and every $m\leq n$, where $m,n\in\mathbb{N}^{*}:=\{1,2,...\}$, we define the map $F_{nm}:X_n\to X_m$ by $F_{nm}(x_n) = x_m$, then the set

\begin{equation}
\varprojlim X_i := \left\{(x_i)_{i=1}^\infty \in \prod_{i=1}^\infty X_i \vert F_{nm}(x_n)=x_m\text{, for all } m\leq n\right\}
\end{equation}

\noindent exists and is unique; it is called the \textit{inverse limit} of the \textit{inverse limit sequence} $\{X_i,F_i\}_{i=1}^\infty$. Furthermore, the maps $F_{nm}$ are called the \textit{transition maps} of the inverse limit system. 

Naturally, if we work instead within the category of Topological Spaces, then all the maps involved should be morphisms of that same category (and hence, here, continuous maps). 

We next give an example of the inverse limit of a particular collection of sets (or rather, of topological spaces).

\begin{example}[The ternary Cantor set $\cantor$]\label{exp:ternaryCantorSet}
Recall that a ternary number is a number expressed in terms of a base-$3$ number system, say, in terms of the characters (or symbols) $\{0,1,2\}$.\footnote{Any three symbols suffice.  One can, and we do so in \S\ref{sec:orbitsOfOmegaKS}, represent elements of the Cantor set in terms of the characters $\{l,c,r\}$, where such characters stand for \textit{left}, \textit{center} and \textit{right}, respectively.  In general, any element in the unit interval $I$ can be represented by a finite or an infinite sequence expressed in terms of an alphabet consisting of the characters $l,c,r$.}  For example, 1102 in base-3 is actually the number 38 in the base-10 number system.  When one constructs the Cantor set $\cantor$, one may do so by removing middle thirds (open intervals) of successive approximations.\footnote{This is not the only way to construct the ternary Cantor set, but most, if not all, methods amount essentially to the same process.}  Let $n\geq 1$.  In removing the middle third from an interval of length $1/3^{n-1}$, we are essentially producing a left third and a right third interval of length $1/3^n$.  As such, we can label the left third as 0 and the right third as 2.  One quickly sees that the Cantor set $\cantor$ is a collection of infinite words written entirely in terms of 0's and 2's.  There are no ternary numbers in $\cantor$ whose address contains 1, because that would mean that we did not remove a middle third from some interval in the construction process.\footnote{For the endpoints of the deleted intervals (also called `ternary points' and necessarily of the form $p/3^q$, for $p,q$ nonnegative integers with $p\leq 3^q$ and $p$ not divisible by $3$, when we restrict our attention to the unit interval $I$), the resulting address is not unique and may contain 1's, although only finitely many 1's.  For example, one may represent $1/3$ by the finite ternary expansion $0.1$ or by the infinite ternary expansion $0.0\bar{2}$, where the overbar indicates that $2$ is repeated ad infinitum.  We adhere to the convention that such ternary numbers are always represented by infinite expansions given in terms of only $0$'s and $2$'s.} More importantly, our labeling system described above enables us to determine particular ternary numbers.

Let $\mathscr{C}_n$ be the $n$th prefractal approximation of the Cantor set.  In the context of the inverse limit construction and addressing system above, $\mathscr{C}_n$ is the collection of $2^n$ points, each having an address given by a finite ternary expansion of length $n$ and each never containing the character $1$.  On the other hand, in the context of the geometric construction detailed in the previous paragraph, $\cantor_n$ consists of $2^n$ compact intervals of length $1/3^n$ (i.e., of `scale $n$').  Furthermore, each address in $\mathscr{C}_n$ may be thought of as an address of a particular segment that remains after removing $2^{n-1}$ segments from the unit interval $I$; see Figure \ref{fig:cantor}.\footnote{Such a construction is called \textit{construction by tremas}, where one removes segments ad infinitum, thereby producing the fractal set.  At each stage, middle thirds are removed from the remaining segments, thereby producing the Cantor set in the limit.}  We take the more geometric interpretation as the definition of $\mathscr{C}_n$; in that case, the sequence $\{\cantor_n\}_{n=1}^\infty$ of prefractal approximations converges to the Cantor set $\cantor$: $\cantor_n\to\cantor$ as $n\to\infty$, in the sense of the Hausdorff metric.  \emph{(}Also, more simply, $\cantor_n$ is monotonically decreasing and $\cantor = \bigcap_{n=1}^\infty \cantor_n$.\emph{)} We will next discuss another way in which $\cantor$ can be viewed as the `limit' of $\{\cantor_n\}_{n=1}^\infty$.

If for positive integers $m\leq n$, we now define the transition map $\tau_{nm}:\cantor_n\to\cantor_m$ as the truncation map that truncates addresses of segments in the prefractal approximations $\mathscr{C}_n$ to addresses of segments in the prefractal approximation $\cantor_m$ by simply removing the last $n-m$ characters from the address, then we form an inverse limit of the prefractal approximations $\cantor_n$ that is exactly the Cantor set $\cantor$:

\begin{equation}
\cantor =\varprojlim \cantor_i = \left\{(c_i)_{i=1}^\infty \in \prod_{i=1}^\infty \cantor_i \vert \tau_{nm}(c_n)=c_m\text{ for all } m\leq n\right\}.
\end{equation}

\begin{figure}
\begin{center}
\includegraphics[scale=.5]{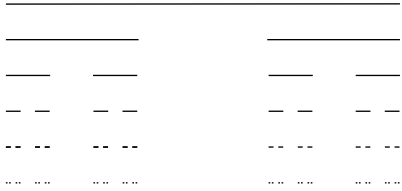}
\end{center}
\caption{The geometric construction of the Cantor set.  One removes middle thirds from remaining segments until all that remains are the points of the Cantor set. Benoit Mandelbrot is credited with the coining of the illustrative name \textit{Cantor Dust}, a name that is quite fitting as it appears that all that remains in the limit is `dust.' Shown here in this figure are the first six approximations of the Cantor set.}
\label{fig:cantor}
\end{figure}
\end{example}

More background information about inverse limits is provided, for example, in \cite[\S2--14 and \S2--15]{HoYo}, where the following well-known theorems can be found; see \cite{HoYo}, Theorems 2--95 and 2--97, along with Corollaries 2--98 and 2--99.  (All the topological spaces considered in Theorems \ref{thm:invLimitTotDisc}, \ref{thm:totCompDisconHomeo} and \ref{thm:twoTotDisconSpacesCompactHomeo} below are implicitly assumed to be metrizable.  Moreover, by the ``Cantor set'', we mean the classic ternary Cantor set discussed in Example \ref{exp:ternaryCantorSet} just above.)

\begin{theorem}
The inverse limit of finite sets is a compact and totally disconnected set.  Conversely, any such topological space is homeomorhpic to an inverse limit of finite sets.
\label{thm:invLimitTotDisc}
\end{theorem}

\begin{theorem}
Any compact and totally disconnected space is homeomorphic to a (closed) subset of the Cantor set.
\label{thm:totCompDisconHomeo}
\end{theorem}

\begin{theorem}
Any two totally disconnected and perfect\footnote{Recall that a subset of a topological space is called \textit{perfect} if it is closed and contains no isolated points.  Hence, a compact space is perfect if it has no isolated points or equivalently, if each of its points is a limit point.} compact spaces are homeomorphic to one another \emph{(}and hence also to the Cantor set\emph{)}.
\label{thm:twoTotDisconSpacesCompactHomeo}
\end{theorem}

\begin{remark}
In the literature on dynamical systems, it is common to use the term \textit{topological Cantor set} to refer to a totally disconnected and perfect compact space (i.e., to a metrizable space that is homeomorphic to the Cantor set).
\end{remark}

\begin{remark}
In \S5, we will show that what we will be referring to as the \textit{footprints} of the \textit{primary piecewise Fagnano orbit} and, in general, \textit{piecewise Fagnano orbits} of the Koch snowflake billiard $\Omega(KS)$ are, in fact, topological Cantor sets (see Theorem \ref{thm:orbitsAreCantorSets}).  Moreover, the footprint of the primary piecewise Fagnano orbit will be the analog of the classic ternary Cantor set, and as a subset of the Koch snowflake $KS$, the union of the footprints of every piecewise Fagnano orbit will constitute a subset of what we will call the \textit{elusive limit points} of the Koch snowflake $KS$.
\end{remark}

\section{The Flat Surface $\mathcal{S}(KS_n)$ as a Branched Cover of $\mathcal{S}(KS_0)$}
\label{sec:SKSnAsABranchedCover}

\begin{figure}
\begin{center}
\includegraphics[scale = .5]{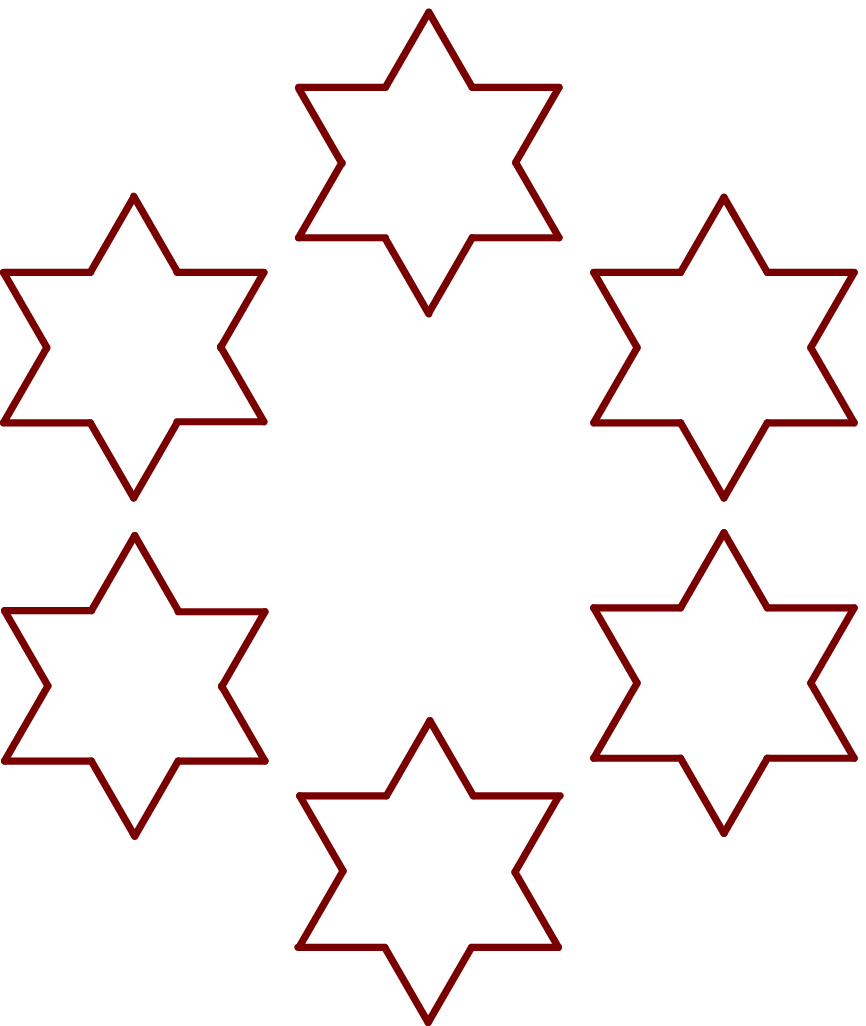}

\includegraphics[scale = .5]{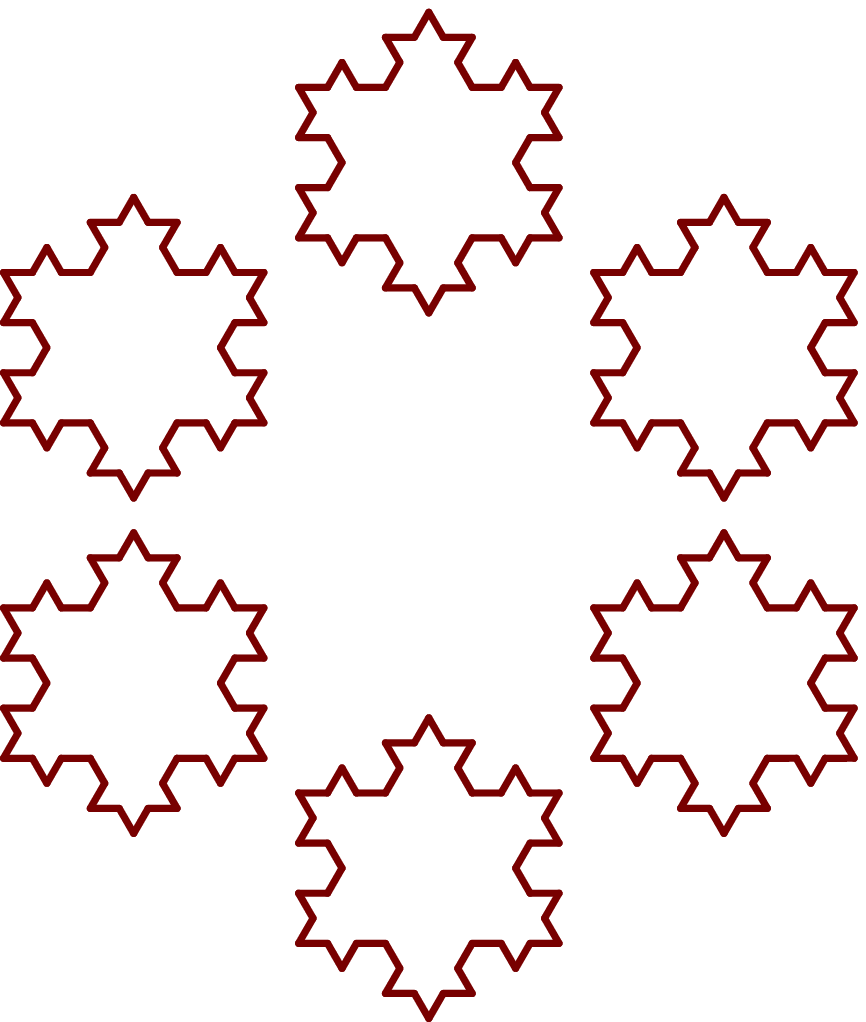}

\includegraphics[scale = .5]{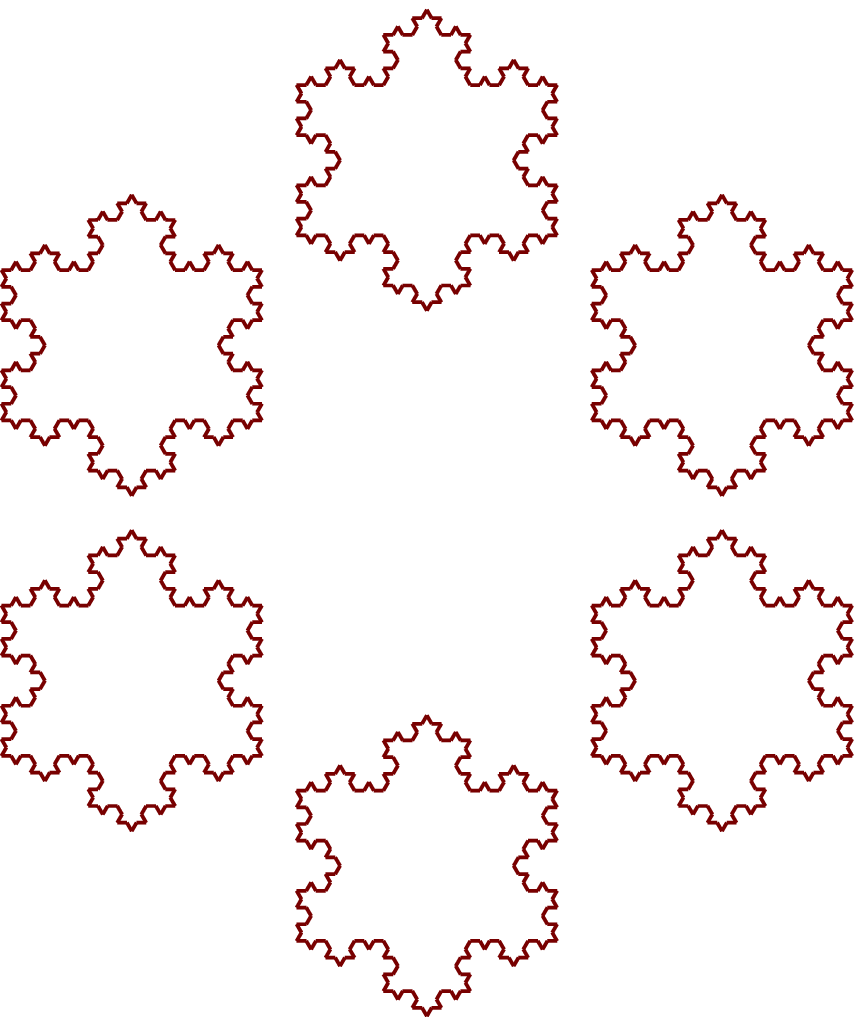}
\end{center}
\caption{The flat surfaces $\mathcal{S}(KS_i)$, $i=1,2,3$.  Note that the proper identification is not shown in the figures above.  Given the arrangement of the six copies of $KS_n$, one then identifies opposite and parallel sides to make the proper identification that results in a geodesic flow that is dynamically equivalent with the billiard flow on the associated billiards $\omegaksi{1},\omegaksi{2},\omegaksi{3}$.}
\label{fig:kochsurface}
\end{figure}

We denote the flat surface $M$ as constructed from a particular rational billiard $\Omega(P)$ by $\mathcal{S}(P)$.  \label{idx:SKSnFlatSurface}In particular, $\mathcal{S}(KS_n)$ is the flat surface associated with the prefractal billiard $\Omega(KS_n)$.  The flat surfaces $\mathcal{S}(KS_n)$, $n=1,2,3$, are given in Figure \ref{fig:kochsurface}.  For each billiard $\Omega(KS_n)$, the group of symmetries $D_N$, where $N=\text{lcm}\{v_i\}_{i=1}^{3\cdot 4^n}$ (that is, the second component in the product $\Omega(KS_n)\times D_N$) is the dyhedral group $D_3$, and thus is independent of $n$.  From this, we deduce that for any $n\geq 0$, there are six copies of the prefractal billiard table $\omegaksi{n}$ (with sides appropriately identified) used in the construction of the associated flat surface $\mathcal{S}(KS_n):=(\omegaksi{n}\times D_3)/\sim$; see Figure \ref{fig:kochsurface}.  We refer the reader back to \S\ref{subsec:flatStructuresandFlatSurfaces} for the discussion of flat surfaces and the associated conical singularities.

\begin{definition}[Covering map]
Let $E$ and $B$ be topological spaces. A \textit{covering map} $p:E\to B$ is a continuous surjective map such that each point $b\in B$ admits a neighborhood $U$ of $b$ for which the preimage $p^{-1}(U)$ is a disjoint collection of open sets in $E$, each of which is mapped homeomorphically  onto $U$ via $p$. One then says that $U$ is \textit{evenly covered} by $p$, and that $E$ is a \textit{covering space} for $B$; see, for example, \cite[Chap. 5]{Ma}.
\end{definition}

\begin{definition}[Branched (or ramified) cover]
Let $E$ and $B$ be topological spaces.  A continuous map $p:E\to B$ is a \textit{branched cover} of $B$ if for all but a finite number of points of $B$, $p$ is a covering map of $E$ onto $B$.  The set of points of $B$ that are not evenly covered by $p$ is called the \textit{branch locus} (or \textit{set of ramification points}).
\end{definition}

\begin{example}[The map $p:\mathbb{C}\to\mathbb{C}$]  The map $p:\mathbb{C}\to\mathbb{C}$, given by $p(z)=z^2$, is a branched covering of $\mathbb{C}$, with branch locus $\{0\}$.  Hence, it is certainly not a covering map.  On the other hand, $p:\mathbb{C}-\{0\}\to \mathbb{C}-\{0\}$, given by the same expression $p(z) = z^2$, is a covering map, since it is \textit{locally trivial}: indeed, each nonzero complex number $z$ in the target space has an open neighborhood $U$ such that $p$, restricted to $p^{-1}(U)$, is equivalent to the projection onto $U\times \{+,-\}$.
\end{example}

For the remainder of the paper, when we say that a regular polygon is of \textit{scale $n$}, we mean that the side length of the regular polygon is $1/3^n$.  For example, an equilateral triangle of scale $n$ is one for which the side length is $1/3^n$.  

Taking as inspiration the results and methods of Gutkin and Judge in \cite{GuJu1} and \cite{GuJu2}, we now show that for each $n\geq 1$, the flat surface $\mathcal{S}(KS_n)$ is a branched cover of the hexagonal torus $\mathcal{S}(KS_0)$; see Corollary \ref{cor:ksnBranchedCoverOfks0}. To such end, we establish several results culminating in the fact that $\mathcal{S}(KS_n)$ is tiled by equilateral triangles of scale $n$.

\begin{lemma}
Let $n\in \mathbb{N}$.  Then, for any positive integer $k\geq n$, $\mathcal{S}(KS_n)$ can be tiled by equilateral triangles of scale $k$.
\label{lem:kochTiledByEqui}
\end{lemma}

\begin{proof}
This follows from the construction of the Koch snowflake.  We note that each triangle of scale $n$, denoted $\Delta_n$, can be tiled by $9^{k-n}$ triangles of scale $k\geq n$; see Figure \ref{fig:equiTiledHex} for the case when $k=n+1$.  Note that $\mathcal{S}(KS_n)=(\Omega(KS_n)\times D_3)/\sim$ and that $\Omega(KS_n)$ is constructed from $\Omega(KS_{n-1})$ by gluing a copy of the billiard table $\omegaksi{n}$ of scale $n$ to every side $s_k$ at the middle third of $s_k$ and then removing the segment common to $s_k$ and $\Delta_{n+1}$. So, $\mathcal{S}(KS_n)$ can be tiled by equilateral triangles of scale $k$. 
\end{proof}

In the sequel, given a bounded set $A\subseteq \mathbb{R}^2$, we will write that ``$A$ \textit{can be tiled by} $H_n$'' in order to indicate that $A$ can be tiled by finitely many copies of hexagonal tiles $H_n$ of scale $n$. 

\begin{lemma}
Let $n\in \mathbb{N}$.  Then the hexagonal torus $\mathcal{S}(KS_0)$ can be tiled by $H_n$.
\label{lem:hextile}
\end{lemma}

\begin{proof}
Consider the hexagonal torus $\mathcal{S}(KS_0)$, as shown in Figure \ref{fig:hexagonalTorus}. In Figure \ref{fig:hexagonalTorusTiled}, we see that $\mathcal{S}(KS_0)$ can be tiled by nine hexagons.. As mentioned in the proof of Lemma \ref{lem:kochTiledByEqui}, $\Delta_n$ can be tiled by nine copies of $\Delta_{n+1}$ overlapping only at the edges.   Moreover, at the center of each $\Delta_n$ is a hexagon of scale $n+1$, $H_{n+1}$; see Figure \ref{fig:equiTiledHex}.

\begin{figure}
\begin{center}
\includegraphics[scale = .5]{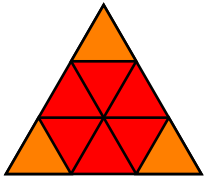}
\end{center}
\caption{We see that $\Delta_n$ is tiled by nine copies of $\Delta_{n+1}$, with a hexagonal tile $H_{n+1}$ in the center.} 
\label{fig:equiTiledHex}
\end{figure}

\begin{figure}
\begin{center}
\includegraphics[scale = .5]{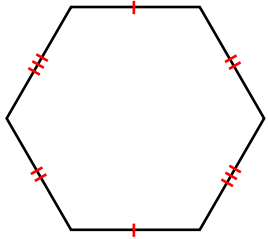}
\end{center}
\caption{The hexagonal torus $\mathcal{S}(KS_0)$. (As usual, similarly marked sides are identified.)  It should be noted that $\mathcal{S}(KS_0)$ is topologically (but not metrically) equivalent to the flat square torus.} 
\label{fig:hexagonalTorus}
\end{figure}

\begin{figure}
\begin{center}
\includegraphics[scale = .5]{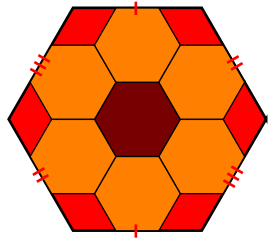}
\end{center}
\caption{We see that there are seven hexagonal tiles of scale $n=1$ tiling $\mathcal{S}(KS_0)$.  By means of the identification, we recover the other two hexagonal tiles, as claimed.} 
\label{fig:hexagonalTorusTiled}
\end{figure}

A hexagon $H_n$ contains six equilateral triangles $\Delta_n$.  As we have seen before, each triangle $\Delta_n$ has at its center a hexagon $H_{n+1}$; see Figure \ref{fig:equiTiledHex}.  These six triangles $\Delta_n$ are then all tiled so that for each $\Delta_n$, there is a $\Delta_{n+1}$ such that a total of six copies of $\Delta_{n+1}$ comprise an additional hexagonal tile $H_{n+1}$ placed at the center of $\Delta_n$. Consequently, $H_n$ is tiled as shown in Figure \ref{fig:hexagonTiled}.  Then, as shown in Figure \ref{fig:threeHexagonTiled}, three copies of $H_n$ arrange so that at each common vertex, there are three rhombic tiles comprising another additional copy of $H_{n+1}$.  Therefore, for any $n$, $\mathcal{S}(KS_0)$ is tiled by hexagons $H_n$.  
\end{proof}

\begin{proposition}
\label{prop:ConicSingAtCenter}
For every $n\in\mathbb{N}$, $\mathcal{S}(KS_n)$ can be tiled by $H_{n+1}$ in such a way that each conical singularity is at the center of a hexagonal tile.
\end{proposition}

\begin{proof}
We proceed by induction on $n$.  First, assume that $n=1$.  Then $\mathcal{S}(KS_1)$ can be tiled by $H_2$ so that every conic singularity of $\mathcal{S}(KS_1)$ is at the center of some hexagonal tile.  Now, given $n\geq 2$, suppose that $\mathcal{S}(KS_{n-1})$ is tiled by hexagons of scale $n$ so that every conic singularity is at the center of some hexagonal tile of scale $n$. 

We first note that the nature of a conic singularity in the surface dictates that such a singularity is actually common to four hexagonal tiles.  We now embed in  $\mathcal{S}(KS_n)$ the tiling of $\mathcal{S}(KS_{n-1})$ by $H_n$.  We see by way of the equivalence relation on $\omegaksi{n}\times D_n$ that for every triangular region $\Delta_n$ in $\mathcal{S}(KS_n)$ that is not in $\mathcal{S}(KS_{n-1})$, there are five other regions in $\mathcal{S}(KS_n)$ that comprise a complete hexagonal tile $H_n$, such that every side of the hexagonal tile is common to some other hexagonal tile in the tiling embedded from $\mathcal{S}(KS_{n-1})$.  Furthermore, by Lemma \ref{lem:hextile}, each hexagon can be tiled as shown in Figure \ref{fig:hexagonTiled}.

\begin{figure}
\begin{center}
\includegraphics[scale = .5]{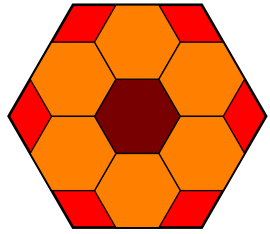}
\end{center}
\caption{Six triangles $\Delta_n$ tile $H_n$.  The hexagonal tile $H_n$ is tiled by seven tiles $H_{n+1}$ with six rhombic tiles.} 
\label{fig:hexagonTiled}
\end{figure}

Note that every hexagon $H_n$ is bordered by six other hexagons $H_n$.  Therefore, by the proof of Lemma \ref{lem:kochTiledByEqui}, every vertex of $H_n$ is at the center of a hexagon $H_{n+1}$.  Moreover, vertices of $H_n$ with edges common to a segment of the polygon $KS_n$ coincide with ternary points $\frac{p}{3^{n+1}}$, where $p=1$ $(\text{mod } 3)$ or $p=2$ $(\text{mod } 3)$.  

\begin{figure}
\begin{center}
\includegraphics[scale = .5]{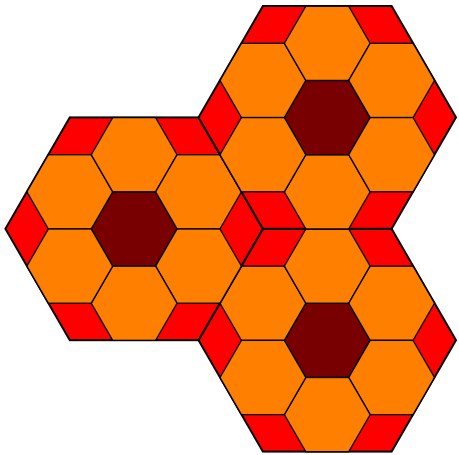}
\end{center}
\caption{We see that three hexagonal tiles $H_n$ tiled as in Figure \ref{fig:hexagonTiled} can be arranged so that a rhombic tile from each all combine to form another hexagonal tile $H_{n+1}$.} 
\label{fig:threeHexagonTiled}
\end{figure}

Consequently, vertices of any tile $H_{n+1}$ with a side common to $\Omega(KS_n)$ that are not at a conic singularity of $\mathcal{S}(KS_n)$ are a distance $\frac{p}{3^{n+1}}$ from a conic singularity of $\mathcal{S}(KS_n)$.  This is exactly the distance to the center of a hexagon tile $H_{n+1}$ along side of $\Delta_{n+1}$ tiling $H_{n+1}$. Hence, we conclude that every conic singularity is at the center of some hexagonal tile $H_{n+1}$, as desired.
\end{proof}

\begin{corollary}
\label{cor:ksnBranchedCoverOfks0}
For every $n\in \mathbb{N}$, the prefractal Koch snowflake flat surface $\mathcal{S}(KS_n)$ is a branched cover of the prefractal Koch snowflake flat surface $\mathcal{S}(KS_0)$.  
\label{prop:branchedCover}
\end{corollary}

\begin{proof}
The center point $x_0$ of the flat hexagonal torus $\mathcal{S}(KS_0)$ is a branched locus of the cover $\mathcal{S}(KS_n)$ when $\mathcal{S}(KS_n)$ is tiled by $H_{n+1}$ as described in Proposition \ref{prop:ConicSingAtCenter}.  This follows from the fact that every conic singularity is at the center of four hexagonal tiles.  Specifically, this means that this center point $z_0$ is not evenly covered by the quotient map $p_n:\mathcal{S}(KS_n)\rightarrow \mathcal{S}(KS_n)/(3^{n+1}H_{n+1})$.\footnote{The notation $3^{n+1}H_{n+1}$ indicates that we are scaling the hexagonal torus of scale $n+1$ by $3^{n+1}$, thereby producing $H_0$, the hexagonal torus.}  Any other point in $\mathcal{S}(KS_0)$ is evenly covered since every element in the fiber $p_n^{-1}(z)$, $z\neq z_0$, has a conic angle of $2\pi$. 
\end{proof}

\begin{figure}
\begin{center}
\includegraphics[scale=.5]{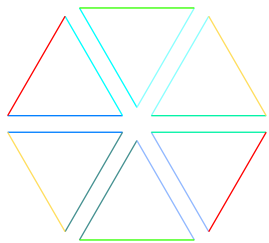}
\end{center}
\caption{Six equilateral triangles glued together correctly constitute a hexagonal torus.  We show the surface in an exploded view so as to emphasize the fact that there are six copies of the equilateral triangle embedded in the surface.  Prior to identifying properly, this figure constitutes the heuristic description of $\Omega(KS_0)\times D_3$.}
\label{fig:sixEquiSurface}
\end{figure}

\begin{figure}
\begin{center}
\includegraphics[scale=.5]{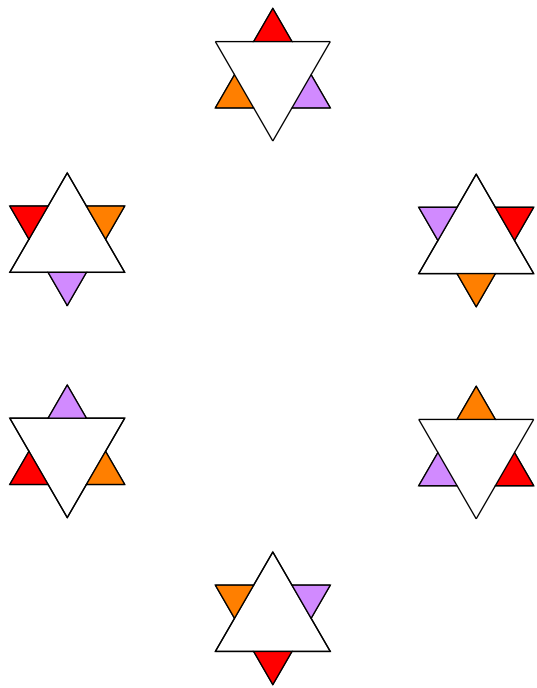}
\end{center}
\caption{Here we show the surface $\mathcal{S}(KS_1)$ in an exploded view.  Opposite sides that are parallel and opposite are identified.  We also see that the surface $\mathcal{S}(KS_1)$ can be viewed as gluing three tori to the existing torus in such a way that the resulting flat surface has genus $g=10$.}
\label{fig:threeToriInSurf}
\end{figure}

\begin{figure}
\begin{center}
\includegraphics[scale=.5]{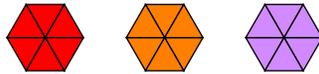}
\end{center}
\caption{We show here the three tori that are appended to the existing torus, rearranged to clearly see the three tori discussed in the caption of Figure \ref{fig:threeToriInSurf}.}
\label{fig:threeToriInSurfRearranged}
\end{figure}

We close this section by highlighting a possible pattern in the construction of the flat surfaces $\mathcal{S}(KS_n)$.  In Figure \ref{fig:threeToriInSurf}, we see that, under the proper re-identification, three equilateral triangle tiles appended to each copy of the equilateral triangle in the flat surface $\mathcal{S}(KS_0)$ (see Figure \ref{fig:sixEquiSurface}) results in the flat surface $\mathcal{S}(KS_1)$.  In Figure \ref{fig:threeToriInSurfRearranged}, the proper re-identification yields three tori that are interconnected in such a way that results in a flat surface with genus $g=10$.  Heuristically, one tears three holes in the hexagonal torus, then glues three additional tori to the existing flat surface in such a way that the proper surface results.  While this is admittedly very difficult to visualize, one can get the impression that each surface $\mathcal{S}(KS_n)$ results from $\mathcal{S}(KS_{n-1})$ by tearing and gluing to the newly opened holes $3\cdot 4^n$ many equilateral triangles of scale $n$ to $\mathcal{S}(KS_{n-1})$.

\section{Periodic orbits of $\mathbf{\omegaksi{n}}$ in the direction of $\pi/3$}
\label{sec:orbitsOfOmegaKSn}
%In CITE BAXTER/UMBLE, Baxter and Umble completely classify periodic orbits of an equilateral triangle billiard $\Omega(\Delta)$.  For our purposes, we let $\Omega(KS_0):=\Omega(\Delta)$.  In the last section we described how closed orbits in $\Omega(KS_n)$ are in exact correspondence with closed orbits in $\Omega(KS_0)$.  Also, we showed that directions in $S^1$ for which flows on $\mathcal{S}(KS_0)$ are uniquely ergodic are in exact correspondence with directions in $S^1$ for which flows on $\mathcal{S}(KS_n)$ are uniquely ergodic.  Consequently, the classification of periodic orbits in $\Omega(\Delta)$ holds for $\Omega(KS_n)$.  Such insights allow us to make the following definitions.

In the previous section, we saw that $\mathcal{S}(KS_n)$ is a branched cover of $\mathcal{S}(KS_0)$.  This implies that for every $n\in\mathbb{N}$, a closed geodesic on $\mathcal{S}(KS_n)$ will project down to a closed geodesic on the hexagonal torus $\mathcal{S}(KS_0)$ under the action of the covering map $p_n$ defined in the proof of Corollary \ref{prop:branchedCover}.  Also, a closed geodesic $\gamma$ on $\mathcal{S}(KS_0)$ lifts to a segment on $\mathcal{S}(KS_n)$ (that is not necessarily closed).  However, according to the discussion in \S\ref{sec:SKSnAsABranchedCover}, there exists a positive integer $k$ such that the lift of $\gamma^k$ is a closed geodesic on $\mathcal{S}(KS_n)$.\footnote{Recall that the notation $\gamma^k$ is meant to represent $k-1$ many concatenations of $\gamma$ with itself: $\gamma*\gamma*...*\gamma = \gamma^k$}  Since the geodesic flow on $\mathcal{S}(KS_n)$ is dynamically equivalent to the billiard flow on $(\omegaksi{n}\times S^1)/\sim$, it follows that a direction giving rise to a closed orbit in $\omegaksi{0}$ is a direction giving rise to a closed orbit in $\omegaksi{n}$ for every $n\geq0$, and vice-versa. We may put this more succinctly as periodic directions in $\mathcal{S}(KS_0)$ are exactly the periodic directions in $\mathcal{S}(KS_n)$, and vice-versa. We summarize the above discussion in the following theorem.

\begin{theorem}
\label{thm:periodicDirectionsAreSame}
The geodesic flow on $\mathcal{S}(KS_0)$ is closed if and only if for every $n\geq 0$, the geodesic flow on $\mathcal{S}(KS_n)$ is closed.  Moreover, the set of directions for which the discrete billiard flow $f_n$ on $\omegaksi{n}$ is closed\footnote{It should be noted that we are making a slight abuse of notation and language.  The billiard flow $f_{n}^t$ is a flow on the phase space ($\omegaksi{n}\times S^1)/\sim$ and $f_n$ is the billiard map.  Iterates of the billiard map then yield elements of the Poincar\'e section.  When we say the billiard flow is closed, we mean that the Poincar\'e section is finite and vice-versa.} \emph{(}i.e., regardless of the initial basepoint, a direction for which a geodesic will be closed\emph{)} is exactly the set of directions for which the billiard flow $f_0$ is closed on $\omegaksi{0}$.
\end{theorem}

\begin{remark}
If $\{e_1,e_2\}$ is a basis for $\mathbb{R}^2$, then a vector $z\in \mathbb{R}^2$  is called \textit{rational with respect to} $\{e_1,e_2\}$ if $z=ne_1+me_2$, $n/m\in \mathbb{Q}$ (that is, $n,m\in \mathbb{Z}$, $m\neq 0$). The plane can be tiled by $\omegaksi{0}$.  As proved in \cite{Gu2}, the collection of directions that gives rise to closed orbits of the equilateral triangle billiard $\omegaksi{0}$ is exactly the set of directions that are rational with respect to the basis $\{e_1,e_2\}= \{(1,0),(1/2,\sqrt{3}/2)\}$. By Theorem \ref{thm:periodicDirectionsAreSame}, for every $n\geq 0$, the same collection of rational directions (with respect to $\{e_1,e_2\}$) describes the directions for which the billiard flow on $\omegaksi{n}$ is closed. 
%We note that, in accordance with the general literature \cite{MasTa,Sm,Ta,Vo}, the collection of directions for which the billiard flow is closed in a rational billiard is at most countable.  In fact, the collection of directions for which the billiard flow is closed in the billiard $\omegaksi{0}$ is given by all directions $\theta$ such that $\tan(\theta)$ is rational with respect to the basis $\{(1,0),(1/2,\sqrt{3}/2)\}$.  By Theorem \ref{thm:periodicDirectionsAreSame}, such directions give rise to closed orbits in $\omegaksi{n}$, for every $n\geq 0$.
\end{remark}

In the sequel, when we say that an angle is measured with respect to a fixed coordinate system, we mean that every direction is measured with respect to the same coordinate system up to translation, but not rotation.  For example, a direction of $\pi/3$ measured with respect to a side of a prefractal billiard table may not be $\pi/3$ when measured with respect to some fixed coordinate system.  Rather, it may be $5\pi/3$; see Figure \ref{fig:thetapi3}.  In order to maintain a consistent measurement, we fix a coordinate system at the base of the equilateral triangle billiard $\Omega(KS_0)$; see Figure \ref{fig:fixedcoordinatesystem}.  

\begin{remark}
\label{rmk:SnkXnkThetank}
Let $n\geq 0$ and $k\leq 3\cdot 4^n$.  We denote by $s_{n,k}$ a side of the billiard $\omegaksi{n}$.  A basepoint of an iterate of the billiard map $f_n:=f_{KS_n}$ of $\omegaksi{n}$ is then denoted by $\xii{n}{k_n}$.  So as to be perfectly clear, we mention that the superscript $k_n$ is not related to the number of sides of $\omegaksi{n}$.  Rather, it is a notation that helps us differentiate between `$k$' iterates of the billiard map $f_n$ and $f_m$, $m\neq n$.  Hence, $f_n^{k_n}$ is the $k_n$th iterate of the billiard map $f_n$.  An initial basepoint of an orbit is then denoted by $\xio{n}$ and an initial direction is $\theta_n^0$.
\end{remark}
Consider an inward pointing vector with basepoint $\xio{n}$ and angle $\phi_n^0$ measured relative to the side of $\omegaksi{n}$ on which $\xio{n}$ resides (assuming $\xio{n}$ is not a vertex).  Since we always want to measure angles relative to a fixed coordinate system, we denote the measure of the angle $\phi_n^0$ relative to the fixed coordinate system by $\theta(\phi)_n^0$. Because of the unique geometry of the prefractal approximation, for every inward pointing vector forming an angle $\phi_n^0$ with a side $s_{n,k}$ on which the basepoint lies, there exists $k\in \mathbb{Z}$ such that $\theta(\phi)_n^0 = k\pi/3\pm \phi$.\footnote{Hence, the angle $\theta(\pi/3)$ is merely an integer multiple of $\pi/3$.}

We have the following definitions and theorems (see Figure \ref{fig:fagnanos012witharrows} for an illustration of Definitions \ref{def:compSeqOrbits} and \ref{def:compSeqInitConditions}):

\begin{figure}
\begin{center}
\includegraphics[scale=.5]{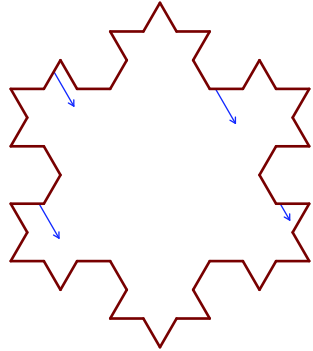}
\end{center}
\caption{While each of these angles can be measured as $\pi/3$ or $2\pi/3$, relative to the sides from which each one emanates, we instead measure relative to a fixed coordinate system given in Figure \ref{fig:fixedcoordinatesystem}, resulting in the measurement of each angle being equal to $5\pi/3$.}
\label{fig:thetapi3}
\end{figure}

\begin{figure}
\begin{center}
\includegraphics[scale=.5]{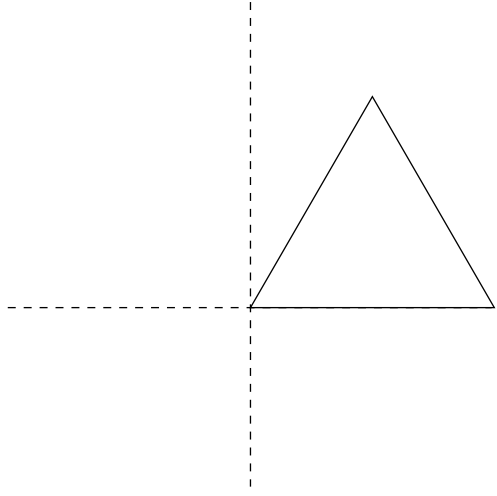}
\end{center}
\caption{The fixed coordinate system relative to which all angles are measured.}
\label{fig:fixedcoordinatesystem}
\end{figure}

\begin{definition}[Compatible sequence of orbits]
\label{def:compSeqOrbits}
Let $\orbitksi{0}$ be an orbit of the equilateral triangle billiard $\Omega(KS_0) =\Omega(\Delta)$.  Assume that $(x^0_n,\theta^0_n)$ is the initial condition of an orbit $\orbitksi{n}$ of $\Omega(KS_n)$ such that $\theta_n^0 = \theta_0^0$ (where each angle is measured relative to the fixed coordinate system) and $\xio{n}$ is collinear with $\xoo$ in the direction $\theta_n^0$ and there are no points of $KS_n$ between $\xio{n}$ and $\xoo$.  Then we say that $\orbitksi{n}$ is \textit{compatible with} $\orbitksi{0}$ and the resulting collection of orbits satisfying such a condition is called a \textit{compatible sequence of orbits}. We write the compatible sequence as $\compseq$.
\end{definition}

\begin{definition}[Compatible sequence of initial conditions]
\label{def:compSeqInitConditions}
Let $\theta_0^0$ be a fixed direction measured relative to the fixed coordinate system. Assume that $\xoo$ is a point on the equilateral triangle $KS_0$ such that $\theta_0^0$ based at $\xoo$ is inward pointing.  Further, assume that $\xio{n}$, $n\geq 1$, is a point on the boundary of $\omegaksi{n}$ such that $\xio{n}$ and $\xoo$ are collinear in the direction $\theta_0^0$ with no points of $KS_n$ common to the segment (interior to $\omegaksi{n}$) joining $\xoo$ and $\xio{n}$, when $\xoo$ is viewed as a point of $\omegaksi{n}$ (if $n=0$, then $\xio{n}$ is trivially collinear with $\xoo$).  Then we say that $\xio{n}$ and $\xoo$ are \textit{compatible} and that the sequence of compatible points $\{\xio{i}\}_{i=0}^\infty$ is a \textit{compatible sequence of initial basepoints}; furthermore, if for each $i\geq 0$, $\theta_{i}^0=\theta_{0}^0$ is an angle in $S^1$ (measured relative to the fixed coordinate system) that corresponds to an inward pointing vector at a basepoint $\xio{i}$, then  $\{(\xio{i},\theta_i^0)\}_{i=0}^\infty$ is called a \textit{compatible sequence of initial conditions}.
\end{definition}

\begin{figure}
\begin{center}
\includegraphics{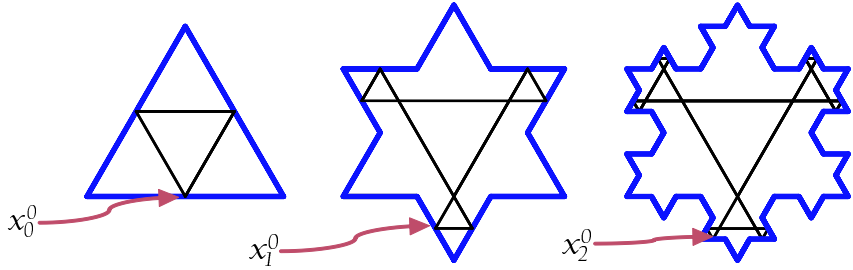}
\end{center}
\caption{A compatible sequence of orbits.  The first three orbits in a compatible sequence are illustrated, and the first three initial basepoints in the compatible sequence of initial basepoints are indicated by the arrows.  Here, the initial direction is $\pi/3$ and $\xoo = 1/2$.}
\label{fig:fagnanos012witharrows}
\end{figure}

\begin{theorem}
\label{thm:compOrbitiffCompInitCond}
A sequence of initial conditions $\{(\xio{i},\theta_{i}^0)\}_{i=0}^\infty$ is a compatible sequence of initial conditions if and only if the corresponding sequence of orbits $\compseq$ is a compatible sequence of orbits.
\end{theorem}

\begin{proof}
This follows easily from Definitions \ref{def:compSeqOrbits} and \ref{def:compSeqInitConditions}.
\end{proof}

\begin{theorem}
If $\orbitksi{0}$ is a closed orbit in $\Omega(KS_0)$, then every orbit in the corresponding compatible sequence of orbits $\{\orbitksiang{i}{\theta^0_i}\}_{i=0}^\infty$ is a closed orbit.  
\end{theorem}

\begin{proof}
Consider a closed orbit $\orbitksi{0}$ of $\omegaksi{0}$.  Then the period of $\orbitksi{0}$ is finite and the unfolding of the orbit in $\mathcal{S}(KS_0)$ is a closed geodesic (or, in the case of the Fagnano orbit, twice the unfolding of the Fagnano orbit is a closed geodesic in $\mathcal{S}(KS_0)$).  Then, by the discussion in \S3, the lift of the geodesic to $\mathcal{S}(KS_n)$ is some segment $\gamma$ such that for some positive integer $k$, $\gamma^k$ is a closed geodesic in $\mathcal{S}(KS_n)$.  Then, the corresponding billiard orbit $\orbitksi{n}$ on $\omegaksi{n}$ is also closed (by the dynamical equivalence between the billiard flow on $\omegaksi{n}$ and the geodesic flow on the flat surface $\mathcal{S}(KS_n)$).  

The orbit $\orbitksi{n}$ is then compatible with $\orbitksi{0}$.  Since $n$ was arbitrary, it follows that $\compseq$ is a compatible sequence of closed orbits.
\end{proof}

\begin{definition}[Compatible sequences of closed and periodic orbits]
\label{def:compSeqClosedPeriodicOrbits}
If every orbit in a compatible sequence of orbits is closed, then we say that it is a \textit{compatible sequence of closed orbits}.  If, in addition, no orbit in the compatible sequence is singular, then we call the sequence a \textit{compatible sequence of periodic orbits}.
\end{definition}

\begin{remark}
Let $\compseqangx{\theta(\phi)_{i}^0}{y}$ be a compatible sequence of orbits.  If we want to discuss a particular compatible sequence of orbits, we use the fact that $\theta(\phi)_i^0 = \theta(\phi)_0^0$ for all $i\geq 1$ and 1) write $\theta(\phi)_0^0$ as $\theta(\phi)$ and 2) $\orbitksixang{n}{y}{\theta(\phi)_n^0}$ as $\orbitksixang{n}{y}{\theta(\phi)}$. In the event that we are discussing a compatible sequence of orbits with $\theta(\phi)_n^0 = \theta(\pi/3)$ for every $n\geq 0$, we write $\orbitksixang{n}{y}{\theta(\phi)_n^0}$ as $\orbitksixang{n}{y}{\theta(\pi/3)}$.
\end{remark}

\begin{remark}
The following definitions, Definitions \ref{def:ghost}--\ref{def:ghostOfCell}, are also given in \cite[\S3 \& \S4]{LaNie1}, up to a few minor differences which have allowed for the convenient definition of \textit{piecewise Fagnano orbit}, which is given in Definition \ref{def:pfn}.
\end{remark}

\begin{figure}
\begin{tabular}{p{5 cm}  p{5 cm}}
\begin{center}\includegraphics{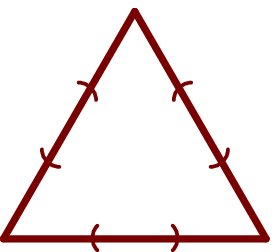}\end{center} & \begin{center}\includegraphics{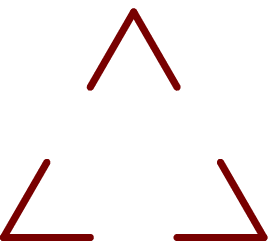}\end{center} \\ 
(a) The ghost set of $KS_0 = \Delta$, denoted by $G_0$. & (b) The elements of the ghost set $G_0$ are removed.\\
\begin{center}\includegraphics{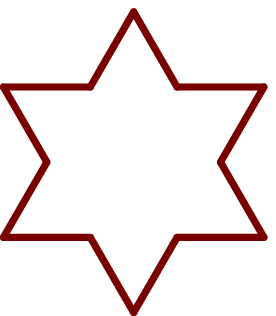}\end{center} & \begin{center}\includegraphics{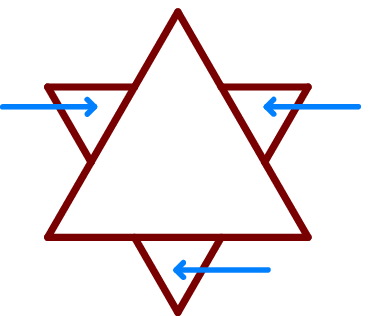}\end{center} \\
 (c) Out of every side there `sprouts' two segments, giving rise to $KS_1$.  & (d) $G_0\cup KS_1$.  The arrows indicate the cells $C_{1,k}$, $1\leq k \leq 3\cdot 4^0=3$, of $KS_1$.\\
\\
\end{tabular}
\caption{An illustration of Definitions \ref{def:ghost} and \ref{def:cell} in terms of $KS_0 = \Delta$ and $KS_1$.  The ghost of the segment $s_{0,k}$, denoted by $g_{0,k}$, is a middle-third segment of $s_{0,k}$ and is removed from $s_{0,k}$ so that we may construct the cell $C_{1,k}$ of $KS_{1}$.  Then $g_{0,k}$ is referred to as the ghost of the cell $C_{1,k}$ of $KS_{1}$.}
\label{fig:illustrateDefis}
\end{figure}

\begin{definition}[Ghosts of $KS_n$]
\label{def:ghost}
 Let $n\geq 0$ and $\{s_{n,k}\}_{k=1}^{3\cdot 4^n}$ be the collection of segments comprising the polygonal boundary $KS_n$ of the billiard $\Omega(KS_n)$.  Then, for $1\leq k \leq 3\cdot4^n$, the open  middle third of the side $s_{n,k}$ is denoted by $g_{n,k}$ and is called the \textit{ghost of the side} $s_{n,k}$.  Moreover, the collection $G_n = \{g_{n,k}\}_{k=1}^{3\cdot 4^n}$ is called the \textit{ghost set of} $KS_n$.  The segments $g_{n,k}$ are removed in order to generate $KS_{n+1}$; see Figures \ref{fig:illustrateDefis}(a)--(c).
\end{definition}

\begin{definition}[A cell $C_{n,k}$ of $\Omega(KS_n)$]
\label{def:cell}
Consider (the `set-theoretic difference') $\Omega(KS_n)\setminus \Omega(KS_{n-1})$.  The resulting triangular regions are then called \textit{cells of} $\Omega(KS_n)$. We denote a cell of $\Omega(KS_n)$ by $C_{n,k}$, where $k$ denotes the side of $\omegaksi{n-1}$ to which the cell was glued; see Figure \ref{fig:illustrateDefis}(d).\footnote{Hence, there are $3\cdot 4^{n-1}$ cells $C_{n,k}$ of $\omegaksi{n}$ and so $1\leq k\leq 3\cdot 4^{n-1}$, in both  Definitions \ref{def:cell} and \ref{def:ghostOfCell}.}
\end{definition}

\begin{definition}[Ghost of a cell $C_{n,k}$]
\label{def:ghostOfCell}
Let $n\geq 1$ and $1\leq k\leq 3\cdot 4^{n-1}$.  If $C_{n,k}$ is a cell of $\Omega(KS_n)$, then the ghost $g_{n-1,k}$ corresponding to the side $s_{n-1,k}$ that was removed in the construction of $\Omega(KS_n)$ is called the \textit{ghost of the cell} $C_{n,k}$. 
\end{definition}

\subsection{Compatible sequences of piecewise Fagnano orbits of $\mathbf{\omegaksi{n}}$}
\label{subsec:compSeqPF}
\begin{figure}
\begin{center} 
\includegraphics{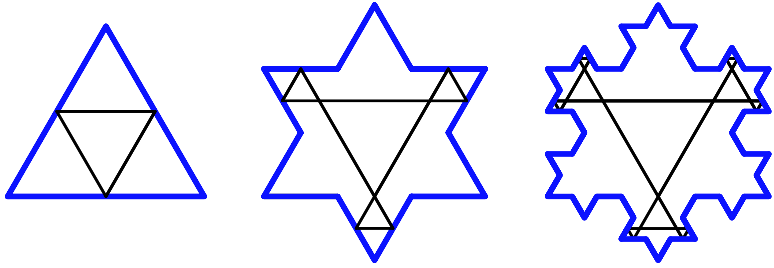}
\end{center}
\caption{The first orbit is the Fagnano orbit of the equilateral triangle billiard $\omegaksi{0}$.  The second orbit is the \textit{primary piecewise Fagnano} orbit of $\omegaksi{1}$ and the third orbit the primary piecewise Fagnano orbit of $\omegaksi{2}$.  These piecewise Fagnano orbits are called \textit{primary}, because the first element in the compatible sequence of orbits to which all of these orbits belong is \textit{the} Fagnano orbit of $\omegaksi{0}$; see also Definition \ref{def:ppfn}.}  
\label{fig:fagnanos012}
\end{figure}

\begin{figure}
\begin{center}
\includegraphics[scale=0.75]{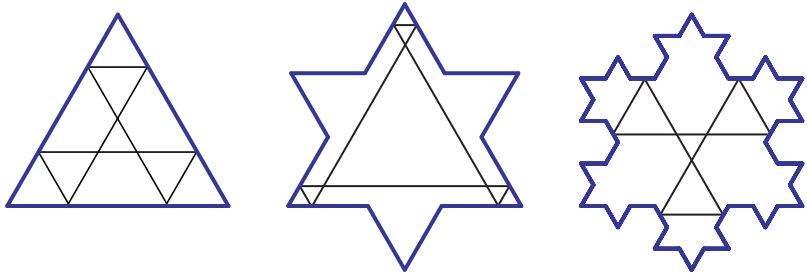}
\end{center}
\caption{Examples of orbits in the direction $\pi/3$.  The first orbit is a non Fagnano orbit of the equilateral triangle.  The second orbit is a non piecewise Fagnano orbit with an initial direction of $\pi/3$ (specifically, the second orbit will be an element is an \textit{eventually constant compatible sequence of periodic orbits}; see \S\ref{subsec:eventConstComSeq}).  The third orbit shown in $\omegaksi{2}$ is an orbit for which \textit{the next orbit} in the corresponding compatible sequence (to which the two orbits would belong) is a piecewise Fagnano orbit of $\omegaksi{3}$, since each basepoint of the orbit lies on a midpoint of some ghost of a cell $C_{3,k}$ of $\omegaksi{3}$.}
\label{fig:piecewiseFagnanoOrbits}
\end{figure}

In $\omegaksi{0}$, the direction $\theta=\pi/3$, measured relative to the fixed coordinate system described in Figure \ref{fig:fixedcoordinatesystem}, gives rise to two different orbits: The Fagnano orbit and a non-Fagnano orbit in the direction of $\pi/3$; see Figures \ref{fig:fagnanos012} and \ref{fig:piecewiseFagnanoOrbits}. The Fagnano orbit is the orbit in the direction $\pi/3$ starting at the midpoint of the base of $\Delta$.  This is, in fact, the shortest orbit of $\Omega(\Delta)$ (see \cite{BaxUm}).  Any other orbit in the direction $\pi/3$ is necessarily twice as long as the Fagnano orbit $\mathscr{F}$.  

We next discuss a generalization of the Fagnano orbit for the Koch snowflake prefractal approximations $\omegaksi{n}$.  

\begin{definition}[Piecewise Fagnano orbits of $\omegaksi{n}$]
If $\yio{n}\in KS_n$ is compatible with a midpoint $\yio{n-1}$ of a ghost $g_{n-1,k}$ of a cell $C_{n,k}$ of $\omegaksi{n}$ in the direction  $\theta(\pi/3)$ (which denotes an integer multiple of $\pi/3$ such that $\theta(\pi/3)$ is an inward pointing direction at $\yio{n}$), then the orbit $\orbitksixang{n}{x}{\theta(\pi/3)}$ is called a \textit{piecewise Fagnano orbit of $\omegaksi{n}$}. 
\label{def:pfn}
\end{definition}

A piecewise Fagnano orbit $\orbitksixang{n}{y}{\theta(\pi/3)}$ is named as such for the fact that one can view such an orbit as the result of appending scale $n$ copies of the Fagnano orbit of $\omegaksi{0}$ to every basepoint of an orbit $\orbitksixang{n-1}{y}{\theta(\pi/3)}$, where $\yio{n-1}$ is the midpoint of a ghost $g_{n-1,k}$ of a side $s_{n-1,k}$ of $\omegaksi{n}$ referred to in Definition \ref{def:pfn}; see Figures \ref{fig:fagnanos012} and \ref{fig:piecewiseFagnanoOrbits}.  We denote a piecewise Fagnano orbit of $\omegaksi{n}$ by $\pfix{n}{\yio{n}}$, where $\yio{n}$ is the initial basepoint of the orbit and $\theta(\pi/3)$ is the initial inward pointing vector.  When $\yio{n}$ is collinear with $\xoo = 1/2\in I$ in the direction $\theta(\pi/3)$ (or when, for some integer $k_n$, a basepoint $\yii{n}{k_n}$ of the orbit $\orbitksixang{n}{y}{\theta(\pi/3)}$ is collinear with $\xoo=1/2\in I$ in the direction $\theta'(\pi/3)$), then we write the orbit $\orbitksixang{n}{y}{\theta(\pi/3)}$ as $\ppfi{n}$ for the \textit{primary piecewise Fagnano orbit} of $\omegaksi{n}$; see Figure \ref{fig:fagnanos012}.

\begin{remark}
\label{rmk:labelingOfpfis}
Later, we will see that for every orbit $\orbitksixang{n}{y}{\theta(\pi/3)}$, there exists a unique element $\xoo\in I$ (the unit interval $[0,1]$ viewed as the base of $\Delta = KS_0$) such that $\xoo$ is compatible with a basepoint $y_n^{k_n}$ of the orbit $\orbitksixang{n}{y}{\theta(\pi/3)}$ in the direction $4\pi/3$.  (See Proposition \ref{prop:determiningPiecewiseFagnanoOrbits} and the discussion preceding it.) Consequently, we will eventually write $\orbitksixang{n}{y}{\theta(\pi/3)}$ as $\orbitksixang{n}{x}{\pi/3}$, since $\orbitksixang{n}{x}{\pi/3}$ determines, and is determined by, the orbit $\orbitksixang{n}{y}{\theta(\pi/3)}$.  Then, in the case when $\orbitksixang{n}{x}{\pi/3}$ is a piecewise Fagnano orbit, we denote the orbit by $\pfix{n}{\xoo}$, since $\xoo$ determines the point $\xio{n}$, and vice-versa.
\end{remark}

\begin{definition}[Primary piecewise Fagnano orbit of $\omegaksi{n}$]\label{def:ppfn}
Let $\pfix{n}{\xoo}$ be a piecewise Fagnano orbit of $\omegaksi{n}$.  If $\xoo = 1/2\in I$, then $\pfix{n}{\xoo}$ is called \textit{the primary piecewise Fagnano orbit} of $\omegaksi{n}$.  We denote the primary piecewise Fagnano of $\omegaksi{n}$ by $\ppfi{n}$.
\end{definition}

We define a compatible sequence of piecewise Fagnano orbits as follows.

\begin{definition}[Compatible sequence of piecewise Fagnano orbits]\label{def:compSeqPF}
The sequence of orbits $\compseqangx{\theta(\phi)}{y}$ is called a \textit{compatible sequence of piecewise Fagnano orbits} if $\compseqangx{\theta(\phi)}{y}$ is a compatible sequence of periodic orbits (in the sense of Definition \ref{def:compSeqClosedPeriodicOrbits}) and there exists $N\geq 0$ such that for every $n\geq N$, $\xio{n}$ is a midpoint of a side $s_{n,k}$ of a cell $C_{n,k'}$ of $\omegaksi{n}$ ($1\leq k'\leq 3\cdot 4^{n-1}$).  Alternately, we say that a compatible sequence of periodic orbits is a compatible sequence of piecewise Fagnano orbits if there exists $N\geq 1$ such that for every $n\geq N$, $\orbitksixang{n}{y}{\theta(\phi)}$ is a piecewise Fagnano orbit of $\omegaksi{n}$; see Figure \ref{fig:fagnanos012}.
\end{definition}

We note that it follows from Theorem \ref{thm:compOrbitiffCompInitCond} that, under the assumptions of Definition \ref{def:compSeqPF}, $\{(\yio{i},\theta(\phi)_i^0)\}_{i=0}^\infty$ is a compatible sequence of initial conditions, in the sense of Definition \ref{def:compSeqInitConditions}.  Moreover, it follows from Definition \ref{def:compSeqPF} that for every $i\geq 0$, there exists $k_i\in \mathbb{N}$ such that $\theta(\phi)_i^0 = k_i\pi/3$.

As alluded to above in Remark \ref{rmk:labelingOfpfis}, we want to be able to describe every piecewise Fagnano orbit $\pfix{n}{\yio{n}}$ of $\omegaksi{n}$ in such a way that $\pfix{n}{\yio{n}}$ is actually an element of some compatible sequence of piecewise Fagnano orbits determined by a particular orbit $\orbitksiang{0}{\pi/3}$ of $\omegaksi{0}$, with $\xoo$ an element of the base of $KS_0$.  In other words, if $\pfix{n}{\yio{n}}$ is an orbit of $\omegaksi{n}$, then we want to show that there exists $\xio{n}\in KS_n$ such that 1) $\xio{n}$ is a midpoint of a side $s_{n,k}$, 2) $\xio{n}$ is collinear with $\xoo$ and 3) $(\yio{n},\theta(\pi/3))\in \orbitksiang{n}{\pi/3}$.

We begin establishing such a connection by considering the following three contractive similarity transformations  $\phi_1,\phi_2,\phi_3:\mathbb{R}\to\mathbb{R}$:

\begin{eqnarray}
\notag\phi_1(x) &:=& \frac{1}{3}x,\\
\label{eqn:IFS}
\phi_2(x) &:=& \frac{1}{3}x+\frac{1}{3},\\
\notag\phi_3(x) &:=& \frac{1}{3}x+\frac{2}{3}.
\end{eqnarray}

Consider $\Phi : \mathbf{K}\rightarrow \mathbf{K}$, a map defined on the space $\mathbf{K}$ of all nonempty compact subsets of $\mathbb{R}$, as $\Phi(K):=\bigcup_{j=1}^3 \phi_j(K)$.  When $\mathbf{K}$ is equipped with the Hausdorff metric, $\mathbf{K}$ becomes a complete metric space. (See, e.g., \cite{Ba,Ed,Fa}.)  Since each $\phi_j$ is a contraction mapping, it follows from Hutchinson's Theorem \cite{Hut} that  $\Phi$ is a contraction mapping on $\mathbf{K}$ and by the Contraction Mapping Principle, $\Phi$ has a unique fixed point attractor in $\mathbf{K}$.  If $I=[0,1]$, we see that $\Phi(I) = \bigcup_{j=1}^3 \phi_j(I) = I$, meaning that the unique fixed point attractor of $\Phi$ is the unit interval $I$.   
%Then, the unique fixed point attractor of $\Phi$ is $[0,1]$, since $\Phi([0,1]) = \bigcup_{j=1}^3 \phi_j([0,1]) = [0,1/3]\cup[1/3,2/3]\cup[2/3,1]$.

\begin{figure}
\begin{center}
\includegraphics{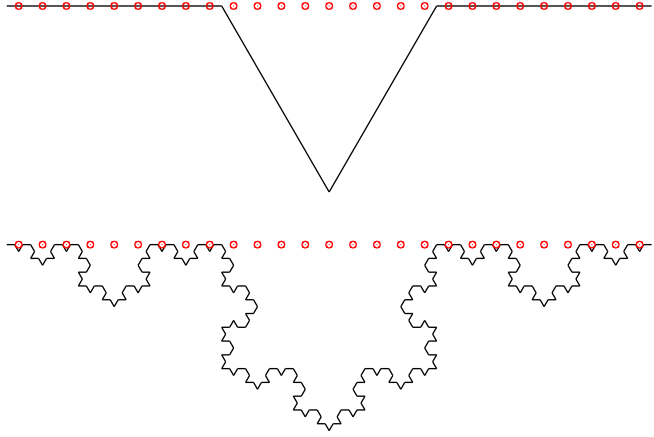}
\end{center}
\caption{In this figure, we demonstrate where exactly the points of $M(3)$ would be located in 1) $\omegaksi{0}$ and 2) $\omegaksi{3}$.}
\label{fig:cantorMidpts}
\end{figure}

For each integer $n\geq 0$, we define $M(n):=\Phi^n(\{1/2\})$, where $\Phi^n$ denotes the $n$th iterate of $\Phi$ and $\Phi^0:= \textit{Id}$; see Figure \ref{fig:cantorMidpts}.  Then it can be easily checked that

\begin{equation}
\label{eqn:MCantor}
M(\cantor):=\bigcup_{n=0}^\infty M(n)
\end{equation}

\noindent  is the collection of all elements of the unit interval $I$ with ternary expansions terminating in $1$'s. Note that each $M(n)$ is finite and hence, that $M(\cantor)$ is countably infinite.  Furthermore, $M(\cantor)$ is dense in the unit interval $I$.

The following lemma will help justify a key step in the proof of the next result (Proposition \ref{prop:ternaryRepOfMc}).

\begin{lemma}
\label{lem:inductionStepForPjiProof}
Let $N\geq 1$.  For every $1\leq i \leq N$, there exists $j_i\in\{1,2,3\}$ and $p_{j_i}\in \{0,1,2\}$ such that

\begin{align}
\phi_{j_N}\circ\phi_{j_{N-1}}\circ\cdots\circ\phi_{j_1}(1/2) &= \sum_{i=1}^N \frac{p_{j_i}}{3^{N-i+1}} +\frac{1}{2\cdot 3^{N}}.
\end{align}
\end{lemma}

\begin{proof}

We know that $\Phi^1(\{1/2\}) = \{1/6,1/2,5/6\}$.    More to the point, each number in $\Phi^1(\{1/2\})$ may be given by $\frac{p_{j_1}}{3^1}+\frac{1}{2\cdot 3^1}$ for a suitable choice of $j_1\in\{1,2,3\}$ and $p_{j_1}\in\{0,1,2\}$. That is, 

\begin{align}
\notag \phi_{j_1} (1/2) &= \left\{\begin{array}{cr}\frac{1}{6} = \frac{0}{3} + \frac{1}{2\cdot 3} &\text{ if } j_1=1, \,\, p_{j_1} = 0\\
\vspace{-2 mm}\\
											\frac{1}{2} = \frac{1}{3}+\frac{1}{2\cdot 3} &\text{ if } j_1 = 2,\,\, p_{j_1} = 1\\
\vspace{-2 mm}\\
											\frac{5}{6} = \frac{2}{3}+\frac{1}{2\cdot 3} &\text{ if } j_1 = 3,\,\, p_{j_1} =2.
\end{array}\right.
\end{align}

Next, let us proceed by induction. Let $N\geq 1$. Suppose that for each $i=1,...,N$, there exist $j_i\in\{1,2,3\}$ and $p_{j_i}\in \{0,1,2\}$ such that

\begin{align}
\notag \phi_{j_{N}}\circ\phi_{j_{N-1}}\circ...\circ\phi_{j_1}(1/2) &= \sum_{i=1}^N \frac{p_{j_i}}{3^{N-i+1}} + \frac{1}{2\cdot 3^{N}}.
\end{align}

\noindent Then

\begin{align}
\notag \phi_{j_{N+1}}(\phi_{j_{N}}\circ...\circ\phi_{j_1}(1/2)) &= \phi_{j_{N+1}}\left(\sum_{i=1}^N \frac{p_{j_i}}{3^{N-i+1}} + \frac{1}{2\cdot 3^{N}}\right)\\
														\notag	&= \left(\sum_{i=1}^N \frac{p_{j_i}}{3^{N-i+1}} + \frac{1}{2\cdot 3^{N}}\right)\frac{1}{3}+\frac{p_{j_{N+1}}}{3}\\
														\notag	&= \sum_{i=1}^{N+1} \frac{p_{j_i}}{3^{(N+1)-i+1}} + \frac{1}{2\cdot 3^{N+1}},
\end{align}

\noindent where $p_{j_{N+1}}\in\{0,1,2\}$.  

Therefore, for every $N\geq 1$, we have

\begin{align}
\notag \phi_{j_N}\circ \phi_{j_{N -1}}\circ ... \circ \phi_{j_1}(1/2)  &= \sum_{i=1}^{N} \frac{p_{j_i}}{3^{N-i+1}} +\frac{1}{2\cdot 3^{N}},
\end{align}

\noindent as desired.
\end{proof}

\begin{proposition}
\label{prop:ternaryRepOfMc}
Let $x\in I$.  Then $x\in M(\cantor)$ if and only if there exists a nonnegative integer $N$ such that for every $1\leq i\leq N$, there exist $j_i\in \{1,2,3\}$ and  $p_{j_i}\in\{0,1,2\}$ such that
\begin{equation}
x = \sum_{i=1}^N \frac{p_{j_i}}{3^{N-i+1}} +\frac{1}{2\cdot 3^{N}}.
\end{equation}
\label{prop:ternaryExpOfMcantor}
\end{proposition}

\begin{proof}
Let $x\in I$ such that $x\in M(\cantor)$. Then $x\in M(N)$ for some $N\geq 0$, where $N$ may be taken as the least such integer $N$. Since $M(N) = \Phi^N(\{1/2\})$, we have that $x\in \Phi^N(\{1/2\})$ and, by Lemma \ref{lem:inductionStepForPjiProof}, for every $1\leq i \leq N$, there exist $j_i\in\{1,2,3\}$ and $p_{j_i}\in \{0,1,2\}$ such that  

\begin{align}
x = \phi_{j_N}\circ \phi_{j_{N-1}}\circ ... \circ \phi_{j_1}(1/2) &= \sum_{i=1}^{N} \frac{p_{j_i}}{3^{N-i+1}} +\frac{1}{2\cdot 3^{N}}. 
\end{align}

%\noindent where for each $1\leq i\leq N$, $j_i$ is an element of $\{1,2,3\}$.  , we have that for each $i\leq N$ and $j_i$, there exists $p_{j_i}\in \{0,1,2\}$ such that

%\begin{align}
%x &= \phi_{j_N}\circ \phi_{j_{N-1}}\circ ... \circ \phi_{j_1}(1/2) \\
%	&= \sum_{i=1}^{N} \frac{p_{j_i}}{3^{N-i+1}} +\frac{1}{2\cdot 3^{N}}. 
%\end{align}

Conversely, suppose now that there exists $N\geq 1$ such that for each $1\leq i\leq N$, there exist $j_i\in \{1,2,3\}$ and  $p_{j_i}\in\{0,1,2\}$ such that

\begin{align}
\notag  x  &= \sum_{i=1}^{N} \frac{p_{j_i}}{3^{N-i+1}} +\frac{1}{2\cdot 3^{N}}. 
\end{align}

Then, we see that for $N=1$,

\begin{align}
\notag \sum_{i=1}^{1} \frac{p_{j_i}}{3^{1-i+1}} +\frac{1}{2\cdot 3^{1}} &= \frac{p_{j_1}}{3}+\frac{1}{2\cdot 3} = \phi_{j_1}(1/2),
\end{align}

\noindent for some $j_1\in\{1,2,3\}$ and $p_{j_1}\in\{0,1,2\}$.  Furthermore, if $N\geq 2$, we have successively:

\begin{align}
\notag x  &= \sum_{i=1}^{N} \frac{p_{j_i}}{3^{N-i+1}} +\frac{1}{2\cdot 3^{N}}\\
	\notag&= \left(\sum_{i=1}^{N-1} \frac{p_{j_i}}{3^{N-i}} +\frac{1}{2\cdot 3^{N-1}}\right)\frac{1}{3} + \frac{p_{j_N}}{3}\\
	\notag&= \left(\left(\sum_{i=1}^{N-2} \frac{p_{j_i}}{3^{N-i-1}} +\frac{1}{2\cdot 3^{N-2}}\right)\frac{1}{3}+\frac{p_{j_{N-1}}}{3}\right)\frac{1}{3} + \frac{p_{j_N}}{3}\\
	\notag&= ...\\
	\notag&= \left(\left(...\left(\left(\frac{p_{j_1}}{3}+\frac{1}{2\cdot 3}\right)\frac{1}{3}+\frac{p_{2,j}}{3}\right)...\right)\right)\frac{1}{3}+\frac{p_{j_N}}{3}\\
	\notag&= \phi_{j_N}\circ...\circ\phi_{j_1}(1/2).
\end{align}

\noindent This concludes the proof of Proposition \ref{prop:ternaryRepOfMc}.
\end{proof}

The number $\sum_{i=1}^N \frac{p_{j_i}}{3^{N-i+1}}$ is an endpoint of an interval $A\subset I$ of length $1/3^{N}$.  By the very nature of the Koch snowflake construction (via an iterated function system), there is a side $s_{N,k}$ of the prefractal approximation $\Omega(KS_{N})$ such that the interval $A$ is a translate (or a translate and rotation by $\pm \pi/3$) of this side $s_{N,k}$ in the direction $4\pi/3$.  
 
Let $\xoo\in M(N)$.  Then, by Proposition \ref{prop:ternaryExpOfMcantor}, there exists $\{p_{j_i}\}_{i=1}^N$, with $j_i\in\{1,2,3\}$ and $p_{i_j}\in\{0,1,2\}$ for each $i=1,...,N$, such that $\xoo = \sum_{i=1}^N \frac{p_{j_i}}{3^{N-i+1}} + \frac{1}{2\cdot 3^{N}}$, meaning, $\xoo$ is the midpoint of the interval in the base of the equilateral triangle for which $\sum_{i=1}^N \frac{p_{j_i}}{3^{N-i+1}}$ is an endpoint.  Therefore, in the direction $4\pi/3$, the endpoint of the line segment connecting $\xoo$ with the segment $s_{N,k}$ of the boundary $KS_{N}$ is a midpoint of the segment $s_{N,k}$.  Denote this point by $\xio{N}$.  Since there are no points of $\omegaksi{N}\setminus \{\xio{N}\}$ on the aforementioned line segment, we know that we may construct a compatible sequence of initial basepoints $\{\xio{i}\}_{i=0}^\infty$ such that $\compseqang{\pi/3}$ is a compatible sequence of piecewise Fagnano orbits.  

%By Remark \ref{rmk:labelingOfpfis} and the paragraph preceding it, we denote such a compatible sequence by $\{\pfix{i}{\xoo}\}_{i=0}^\infty$, where $\xoo\in M(\cantor)\setminus T$.

%This means that $\{x_i\}_{i=0}^{N}$ is a compatible sequence of initial conditions with $x_N$ a midpoint of $s_{N,k}$.

Conversely, we want to show that if $\yio{n}$ is a midpoint of a side $s_{n,k}$ of $\omegaksi{n}$, then, for some $N\geq 0$, $\orbitksixang{n}{y}{\theta(\pi/3)}$ determines an element $\xoo\in M(N) \subseteq M(\cantor)$, and ultimately, a compatible sequence of initial conditions that corresponds to a compatible sequence of piecewise Fagnano orbits $\compseqang{\pi/3}$.\footnote{Recall that the angle $\theta(\pi/3)$ is an integer multiple of $\pi/3$, measured with respect to the fixed coordinate system shown in Figure \ref{fig:fixedcoordinatesystem}.} Consider an inward pointing direction $\theta(\pi/3)$ at a midpoint $\yio{n}\in KS_n$ on a side $s_{n,k}$ of $\omegaksi{n}$.  The billiard ball must reflect off of at most one other side of a cell in which the side $s_{n,k}$ was a part of (this may not be a cell $C_{n,k}$, but rather a cell $C_{n',k'}$, with $n'\leq n$ and $k'\leq 3\cdot 4^{n'}$) before exiting the cell (by convention, a cell $C_{0,k}$ is the equilateral triangle); see Figure \ref{fig:onlyReflectTwice}.  

%A billiard ball entering into a cell $C_{n',k'}$ must exit $C_{n',k'}$ in a direction that would have been the same had the pointmass reflected off of the ghost of the cell, $g_{n',k'}$; see Figure \ref{fig:possibleEntriesExits}.  

%By means of the local and global symmetry of the prefractal table $\omegaksi{n}$, as shown in Figure \ref{fig:kochSymmetries} and discussed in the corresponding caption, one can deduce that there must be a point $x_n^{k_n}$ for which the billiard ball makes a collision and $x_n^{k_n}$ is a midpoint of a side, such that the billiard ball next traverses the interior of $\omegaksi{n}\cap\omegaksi{0}$.  We claim that upon doing so, the billiard ball must pass through some point of $M(N)$, $N\leq n$.  Indeed, if it did not, then $\yio{n}$ could not have been a midpoint of any side of $\omegaksi{N}$. This follows from the fact that $\yio{n}$, relative to the side on which it lies, is the point on the side with the value $1/(2\cdot 3^n)$, meaning that translating $x_n^{k_n}$  in the direction of $\pi/3$ results in identifying an element of $I$ that also has a ternary expansion terminating solely in 1's.  In fact, the element $\xoo$ in the unit interval collinear with $\xio{n}$ is equal to $\sum_{i=1}^N p_{j_i}/3^{N+1-i} + 1/(2\cdot 3^{N})$ for some finite collection $\{p_{j_i}\}_{i=1}^N$, which is an element of $M(N)$.   Such an orbit $\orbitksiang{n}{\theta(\pi/3)}$ contains an element $(\yio{n},\pi/3)$. Letting $y=\xio{n}$, we see that $\{\orbitksi{i}\}_{i=0}^\infty$ is a compatible sequence of periodic orbits.  

By means of the local and global symmetry of the prefractal billiard table $\omegaksi{n}$, as shown in Figure \ref{fig:kochSymmetries} and discussed in the corresponding caption.  One can deduce that there must be a basepoint $\yii{n}{k_n}$ of the orbit $\orbitksixang{n}{y}{\theta(\pi/3)}$ that corresponds to a midpoint and that, after colliding with the boundary $KS_n$, the billiard ball next traverses the interior of $\omegaksi{n}\cap\omegaksi{0}$ in the direction of $\pi/3$. We claim that upon doing so, the billiard ball must pass through some element of $M(\cantor)$.  Indeed, if it did not, then $\yio{n}$ could not have been a midpoint of any side.\footnote{Recall that $\yii{n}{k_n}$ can also be viewed as the image of $\yio{n}$ under the action of particular local and global symmetries of $\omegaksi{n}$.}  Specifically, $\yio{n}$ has a value $\frac{1}{2\cdot 3^n}$, relative to $s_{n,k}$, and $\yii{n}{k_n}$ has same the same value, but relative to some side $s_{n,k'}$, $1\leq k' \leq 3\cdot 4^n$.  Since $s_{n,k'}$ can be translated (or translated and rotated by $\pm\pi/3$) in the direction of $\pi/3$ so as to correspond with an interval of length $1/3^n$ in $I$, it follows that the billiard ball passes through an element of $M(\cantor)$; call this element $\xoo$.  Then, if we define $\xio{n}$ by $\xio{n} :=\yii{n}{k_n}$, the orbit $\orbitksiang{n}{\pi/3}$ contains $(\yio{n},\theta(\pi/3))$ and constitutes an element in a compatible sequence of piecewise Fagnano orbits.

We summarize our discussion in the following proposition:

\begin{figure}
\begin{center}
\includegraphics{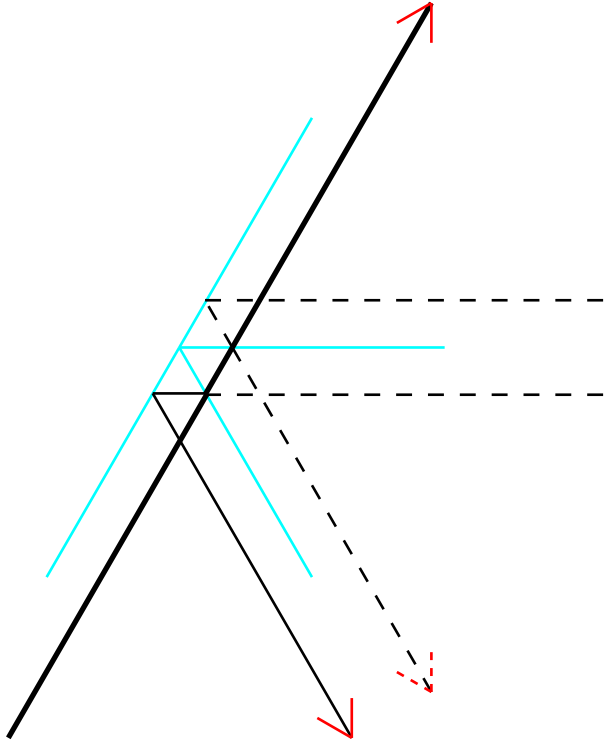}
\end{center}
\caption{Depicted in this figure is a classical tool for analyzing billiard orbits: \textit{unfolding the orbit}.  We consider a billiard orbit's trajectory entering a cell $C_n$ at an angle of $\pi/3$.  As expected, the billiard ball must exit the cell after two reflections and at an angle of $-\pi/3$.  We verify this by considering the unfolded trajectory in the plane and noticing that the trajectory continues on unimpeded as it passes through a reflected copy of the opening. (See Figure \ref{fig:possibleEntriesExits} for a generalization of this discussion to the case of two collisions with the boundary.)}
\label{fig:onlyReflectTwice}
\end{figure}

\begin{figure}
\begin{center}
\includegraphics{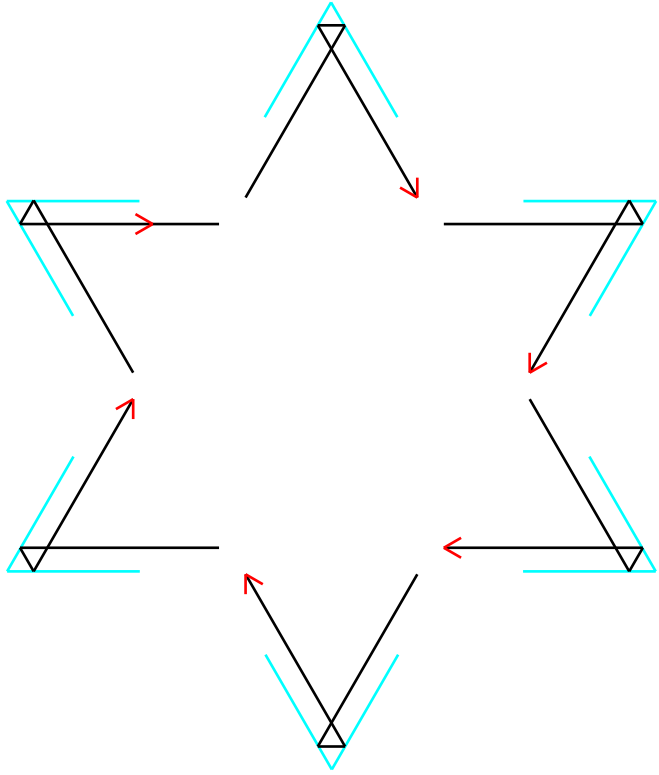}
\end{center}
\caption{This figure generalizes what we have seen in Figure \ref{fig:onlyReflectTwice}.  An orbit's trajectory that enters a cell $C_n$ parallel to a side of $C_n$ will exit after two collisions.}
\label{fig:possibleEntriesExits}
\end{figure}

\begin{figure}
\begin{center}
\includegraphics{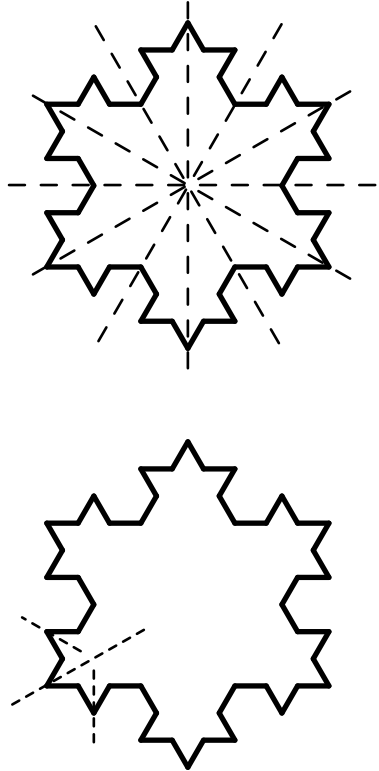}
\end{center}
\caption{For $n\geq 1$, the symmetry group of $KS_n$ is the dihedral group $D_6$ (see \S\ref{subsec:flatStructuresandFlatSurfaces}).  Local symmetry is seen at the level of a cell $C_n$ of $\Omega(KS_n)$.  The symmetry of interest is symmetry with respect to the angle bisector of the acute angle of the cell $C_n$.}
\label{fig:kochSymmetries}
\end{figure}

\begin{proposition}
\label{prop:determiningPiecewiseFagnanoOrbits}
Let $n \geq 0$ and $\orbitksixang{n}{y}{\theta(\pi/3)}$ be compatible with a piecewise Fagnano orbit $\pfix{N}{\yio{N}}$.  Then there exists a unique element $\xoo \in M(\cantor)$ and a compatible sequence of piecewise Fagnano orbits $\compseqang{\pi/3}$ such that  

\begin{align}
\orbitksixang{n}{y}{\theta(\pi/3)} \in\compseqang{\pi/3}
\end{align}

\noindent and
\begin{align}
\pfix{m}{\yio{m}} = \orbitksiang{m}{\pi/3}, \text{ for every }m\geq N. 
\end{align}
%%
%%\noindent

%Let $n\geq 0$ and let $\orbitksixang{n}{y}{\theta(\pi/3)}$ be an element of a compatible sequence of piecewise Fagnano orbits $\compseqang{\pi/3}$.  Then there exists a unique $\xoo\in M(\cantor)$ such that $\orbitksix{n}{x} = \orbitksix{n}{y}$ for all $n\geq 0$ and $\{\orbitksixang{i}{x}{\pi/3}\}_{i=0}^\infty$ forms (in the sense of Definition \ref{def:compSeqPF}) a compatible sequence of piecewise Fagnano orbits 

Conversely, for every $\xoo\in M(\cantor)$, there exists $n\geq 0$ and a piecewise Fagnano orbit $\pfix{n}{\xio{n}}$ of $\omegaksi{n}$ such that $\xio{n}$ is compatible with $\xoo$ and $\orbitksiang{0}{\pi/3}$ is compatible with $\pfix{n}{\xio{n}}$.  Moreover, there is a least integer $n\geq 0$ such that $\xio{n}$ is compatible with $\xoo$ and for every integer $k\geq n$, the compatible orbit $\orbitksi{k}$ is a piecewise Fagnano orbit $\pfix{n}{\xio{n}}$ determined by $\xoo$.
\end{proposition}

We then are in a position to say that $\xoo\in M(\cantor)$ determines a piecewise Fagnano orbit $\pfix{n}{\xio{n}}$ and, conversely, $\pfix{n}{\xio{n}}$ determines a unique element $\xoo\in M(\cantor)$, allowing us to write $\pfix{n}{\xio{n}}$ as $\pfix{n}{\xoo}$ without any ambiguity.

\subsection{Eventually constant compatible sequences of periodic orbits}
\label{subsec:eventConstComSeq}
%The ternary Cantor set $\cantor$ is the collection of all elements in the unit interval $I$ with ternary expansions expressed solely in terms of 0's and 2's.  Elements of the complement have ternary expansions containing at least one 1.  

\label{idx:baseThree}
Consider the ternary Cantor set $\cantor$.  It is the collection of all elements in the unit interval $I$ such that an element has a ternary expansion consisting only of $0$'s and $2$'s.  This is not to imply that an element of $\cantor$ will not have an expansion consisting of any number of $1$'s, but that such a number must have some ternary representation satisfying this rule.  For example, denoting a base-$3$ number in $I$ by $0.u_1u_2..._3$, we see that $1/3 = 0.1_{3}$ is in $\cantor$ because $1/3 = 0.0\overline{2}_3$, which has no $1$'s.  Similarly, $1/4 = 0.\overline{02}_3$ is also in $\cantor$.  However, $1/2 = 0.\overline{1}_3$ is not in $\cantor$, because $1/2$ does not have an equivalent representation that would satisfy the rule for being in $\cantor$.\footnote{Indeed, $1/2$ is in the deleted interval $(1/3,2/3)$ at the first stage of the construction of $\cantor$.}  The complement $\cantor^c=I\setminus \cantor$ (with respect to the unit interval $I$) can be described as the collection of elements in $I$ having a ternary expansion comprised of at least two $1$'s.  For example, $0.11_3$ is not in the Cantor set, because it lies in the open interval $(1/3,2/3)$, which is removed during the geometric construction process as shown in Figure \ref{fig:cantor}. (See Example \ref{exp:ternaryCantorSet} of \S\ref{subsec:invLimSeqInvLim}.)

Let $\cantor'$ be the collection of all elements of $\cantor$ with infinite ternary expansions that do not terminate solely in a single character (repeated ad infinitum). These elements are naturally represented as infinite sequences of 0's and 2's, with periodic expansions representing rational values and aperiodic expansions representing irrational values of $\cantor$. 

In any prefractal approximation $\cantor_n$ of $\cantor$, every $x\in \cantor'$ is contained in a (one-dimensional) connected open neighborhood $U_x$ such that $U_x\cap \mathscr{C}_n = U_x$.  Therefore, for every $n\geq 0$ and  $\xoo\in \cantor'\subseteq I$, $\orbitksiang{n}{\pi/3}$ is a well-defined periodic orbit of $\Omega(KS_n)$.\footnote{This means that it does not hit any of the vertices of $KS_n$.}  Moreover, the compatible sequence of periodic orbits $\{\orbitksiang{i}{\pi/3}\}_{i=0}^\infty$ is an eventually constant sequence of orbits, since $\xoo\in KS_n$, for every $n\geq 0$.

%We may extend this result.  Every side $s_{n,k}$ of $\Omega(KS_n)$, $1\leq k\leq 3\cdot 4^n$, contains a scaled (or scaled and rotated) copy of the ternary Cantor set; denote the points of this scaled Cantor set by $\cantor_{n,k}$ and the points of $\cantor_{n,k}$ which have infinite ternary expansions not terminating solely in 0's nor 2's by $\mathscr{C}_{n,k}'$.  We have already demonstrated (in \S\ref{subsec:compSeqPF}) that an orbit with an initial direction of $\theta(\pi/3)$ emanating from a midpoint of a side $s_{n,k}$ of $\omegaksi{n}$ is an element of a compatible sequence of piecewise Fagnano orbits. The same exact argument that allowed us to say that each piecewise Fagnano orbit $\pfix{n}{\yio{n}}$ of $\omegaksi{n}$ could be determined by $(\xoo,\pi/3)$, where $\xoo\in M(\cantor)$, also allows us to say that an orbit $\orbitksixang{n}{y}{\theta(\pi/3)}$ with $\yio{n}\in \mathscr{C}_{n,k}'$ can be shown to be an element of a compatible sequence of orbits entirely determined by some element $\xoo$ with a ternary expansion consisting of infinitely many $0$'s and $2$'s, where $M(\cantor) \subseteq I$ is the subset of the base of $KS_0$ (identified with $I$) defined by equation (\ref{eqn:MCantor}); see \S\ref{subsec:compSeqPF}.  

We may extend this result.  In order to do so, we first recall Proposition \ref{prop:determiningPiecewiseFagnanoOrbits}: $\xoo \in M(\cantor)$ determines and is determined by a piecewise Fagnano orbit of $\omegaksi{n}$, $n\geq 1$.  The proof of this statement (given in the discussion preceding the proposition) relies heavily on the local and global symmetry of the prefractal snowflake $\omegaksi{n}$.  We take advantage of such symmetries when constructing an orbit $\orbitksiang{n}{\pi/3}$  of $\omegaksi{n}$ with an initial basepoint $\xio{n}$ on the boundary $KS_n$ of $\omegaksi{n}$. 

\label{idx:cnkprime}\label{idx:cnk}Consider $\cantor_{n,k}\subseteq S_{n,k}$ a scaled copy of $\cantor$, where $n\geq 0$ and $1\leq k\leq 3\cdot 4^n$.  Then, define $\cantor_{n,k}'$ to be the collection of all points of $\cantor_{n,k}$ with a ternary expansion consisting of infinitely many $0$'s and $2$'s.  If $\yio{n}\in\cantor_{n,k}'$, then the orbit $\orbitksixang{n}{y}{\theta(\pi/3)}$ is a well-defined orbit of $\omegaksi{n}$.  In addition, there exists a basepoint $\yii{n}{k_n}$ of the orbit $\orbitksixang{n}{y}{\theta(\pi/3)}$ such that after colliding in the boundary at $\yii{n}{k_n}$, the billiard ball then traverses the interior of $\omegaksi{n}\cap\omegaksi{0}$.  In analogy with the case when $\yii{n}{k_n}$ was a midpoint of a side $s_{n,k}$, the billiard ball must now pass through a point $\xoo\in I $ with a ternary expansion consisting of finitely many $1$'s and infinitely many $0$'s and $2$'s. Moreover, this compatible sequence is eventually constant, meaning that there is a positive integer $N$ such that $\orbitksiang{n}{\pi/3}=\orbitksixNOi{N}{x}$, for all $n\geq N$.  We have therefore proved the following proposition.  %(See the beginning of \S4.3 below for the choice of such initial basepoints $\xio{n}$.)

\begin{proposition} \label{prop:constantCompSeq}
Let $N\geq 0$ and $k\leq 3\cdot 4^N$.  If $\yio{N} \in \mathscr{C}_{N,k}'$, then there exists $\xio{N}\in KS_N$ compatible with some element $\xoo\in I$ having a ternary representation consisting of finitely many $1$'s and infinitely many $0$'s and $2$'s and such that  $\{\orbitksiang{i}{\pi/3}\}_{0}^\infty$ is an eventually constant compatible sequence of periodic orbits with $\orbitksixang{N}{y}{\theta(\pi/3)}\in \{\orbitksiang{i}{\pi/3}\}_{0}^\infty$.
\end{proposition}

\begin{remark}
\label{rmk:eventuallyConstCompSeq}
If there exists a positive integer $N\geq 0$ such that $\{\orbitksixang{i}{x}{\pi/3}\}_{i\geq N}$  is a constant compatible sequence of periodic orbits, then $\{\orbitksixang{i}{x}{\pi/3}\}_{i=0}^\infty$ is an \textit{eventually constant compatible sequence of periodic orbits}. 
\end{remark}

As a result of Proposition \ref{prop:constantCompSeq} and Remark \ref{rmk:eventuallyConstCompSeq}, we have the following definition.

\begin{definition}
\label{def:eventuallyConstantCompSeq}
If $\xoo\in I$ has a ternary expansion consisting of finitely many $1$'s and infinitely many $0$'s \textit{and} $2$'s, then we call the resulting compatible sequence of periodic orbits $\{\orbitksixang{i}{x}{\pi/3}\}_{i=0}^\infty$ an \textit{eventually constant compatible sequence of periodic orbits}.
\end{definition}

\begin{example}
Consider the element $\xoo \in I$ with the ternary representation given by $0.10\overline{20}_3$, where, as usual, the overbar indicates that the corresponding string $20$ is repeated ad infinitum.  This is the value $1/3+1/12=5/12$, which is not an element of $M(\cantor)$.  Then, choose $\xio{1}$ in $KS_1$ on a side $s_{1,k}$ such that $\xio{1}$ is collinear with $\xoo$ in the direction $\pi/3$, relative to the fixed coordinate system.  This point identified is the point $0.\overline{02}_3$ scaled by $3$, or $1/4$ scaled by $3$,  residing on the side $s_{1,k}$.  Since $1/4$ is a point of the Cantor set (because $1/4=0.\overline{02}_3$), $1/12$ is a point of the Cantor set, and hence, the point $\xio{1}$ is the point $1/12$ on the side $s_{1,k}$ and in the scaled Cantor set $\cantor_{n,k}$.\footnote{This can be more clearly seen if one views the side $s_{1,k}$ as a rotation and translation of the left middle third of the unit interval $I$ (viewed as the base of $KS_0$).}

Consequently, the orbit of $\omegaksi{1}$ is compatible with that of $\omegaksi{0}$ and the corresponding compatible sequence of periodic orbits is an eventually constant compatible sequence of periodic orbits.
\end{example}

If $\xoo$ is an element with a ternary expansion consisting of infinitely many $0$'s and $2$'s and finitely many $1$'s, then, by Definition \ref{def:eventuallyConstantCompSeq}, the compatible sequence of periodic orbits $\compseqang{\pi/3}$ is an eventually constant compatible sequence of periodic orbits.  If $N\geq 0$ is such that $\compseqangi{\pi/3}{N}$ is a constant sequence of compatible orbits, then, for each $n\geq N$, we call the corresponding orbit  $\orbitksiang{n}{\pi/3}$ a  $\cantor$\textit{-orbit} of $\omegaksi{n}$. In the context of the fractal billiard $\omegaks$, we will refer to a $\cantor$-orbit as a \textit{stabilizing periodic orbit} (or, simply, a stabilizing orbit) of the Koch snowflake billiard $\omegaks$; see \S\ref{subsec:stabilizingOrbitsKS}.

\subsection{Compatible sequences of generalized piecewise Fagnano orbits of $\mathbf{\omegaksi{n}}$}
A recurring theme thus far in \S 4 is that if $(\yio{n},\theta(\pi/3))$ is an initial condition of a periodic orbit $\orbitksixNOiang{n}{y}{\theta(\pi/3)}$ in $\omegaksi{n}$, then there exists $[(\xio{n},\pi/3)]\in (\omegaksi{n}\times S^1)/\sim$ such that $(\xio{n},\pi/3) \in \orbitksixNOiang{n}{y}{\theta(\pi/3)}$, $\xio{n}$ is collinear with a point $\xoo$ in the unit interval $I$ and $\compseqang{\pi/3}$ forms a compatible sequence of periodic orbits.  The points $x$ of the billiard table $\omegaksi{n}$ for which we have shown this to be true are 

\begin{enumerate}
\item $x\in \omegaksi{n}$ such that $x\in \cantor_{n,k}' \subset s_{n,k}$, for some $k\leq 3\cdot 4^n$;
\item $x\in \omegaksi{n}$ such that for every $n'\geq n$, there exist $k'\leq 3\cdot 4^{n'}$ and $x' \in \cantor_{n',k'}'\subset s_{n',k'}$ such that $x$ and $x'$ are compatible in the direction of $\theta(\pi/3)$;% in the direction $\theta(\pi/3)$ (measured relative to $s_{n,k})$;
\item $x\in s_{n,k}$ such that $x\in M(\cantor)_{n,k}$, the set $M(\cantor)_{n,k}$ being a scaled (or scaled and rotated) copy of $M(\cantor)$, viewed as a subset of $s_{n,k}$.
\end{enumerate}

A point $\yio{n}$ on a side $s_{n,k}$ not falling into any of these categories (and which, for any $m\geq 0$, is not compatible with a corner of $\omegaksi{m}$ in the direction of $\theta(\pi/3)$) constitutes an initial basepoint for which we call the resulting orbit $\orbitksixang{n}{y}{\theta(\pi/3)}$ a \textit{generalized piecewise Fagnano orbit} of $\omegaksi{n}$.  Because we are measuring angles with respect to the fixed coordinate system, recall that $\theta(\pi/3)$ is our convention for indicating that $\pi/3$ was the measure of the angle when measured relative to the side $s_{n,k}$ on which $\yio{n}$ lies and $\theta(\pi/3)$ is then the measure of the same angle measured relative to our fixed coordinate system.  We note that such an orbit is never a piecewise Fagnano orbit nor is it an orbit $\orbitksixNOiang{n'}{y}{\theta(\pi/3)}$ such that $\yio{n}$ and $\yio{n'}$ are collinear in the direction $\theta(\pi/3)$ and $\orbitksixNOiang{n'}{y}{\theta(\pi/3)}$ is a piecewise Fagnano orbit for any $n'\geq n$.  As was done before, we may find $\xio{n}\in s_{n,k}$ such that $\orbitksiang{n}{\pi/3} = \orbitksixNOiang{n}{y}{\theta(\pi/3)}$ and $\xio{n}$ is collinear with some $\xoo \in I$ so that $\compseqang{\pi/3}$ constitutes a compatible sequence of periodic orbits.  We call such a compatible sequence a \textit{compatible sequence of generalized piecewise Fagnano orbits}.

The element $\xoo$ in the corresponding compatible sequence of basepoints is not an element of $M(\cantor)$; otherwise, there would exist $N\geq 0$ and $\orbitksiang{N}{\pi/3}$  in the compatible sequence such that the basepoints of $\orbitksiang{N}{\pi/3}$ would have to be midpoints of particular sides of some prefractal approximation $\omegaksi{N}$.  In addition, if $\xio{N}$ is the initial basepoint of a generalized piecewise Fagnano orbit, then the initial basepoint $\xoo\in I$ of a compatible orbit $\orbitksiang{0}{\pi/3}$ of $\omegaksi{0}$ has a ternary expansion consisting of infinitely many $0$'s and $1$'s, $2$'s and $1$'s or $0$'s, $1$'s and $2$'s; otherwise, $\orbitksiang{\pi/3}$ would be an orbit in an eventually constant compatible sequence of orbits or a compatible sequence of piecewise Fagnano orbits.

%In fact, the ternary representation of $\xoo$ is given by the following expression:

%\begin{equation}
%\xoo = \sum_{i=1}^\infty \frac{p_{j_i}}{3^i} + \frac{1}{2\cdot 3^k},
%\end{equation}

%\noindent for some $k\geq 0$ and some (possibly aperiodic) sequence $\{p_{j_i}\}_{i=1}^\infty$ (where $j_i\in\{1,2,3\}$ for each $i\geq 1$) such that $\{p_{j_i}\}_{i=1}^\infty$ contains infinitely many 0's, 1's and 2's, guaranteeing that $x_0$ 1) is not an element of $M(\cantor)$ and 2) is not collinear in the direction $\pi/3$ with any element $x_n'$ of any set $\cantor_{n',k'}$ for any $n'\geq n$.  

\subsection{Properties of periodic orbits of $\mathbf{\omegaksi{n}}$ in the direction of $\pi/3$}
\label{subsec:aGeneralDiscussionOfOrbits}
According to the preceding discussion, in each prefractal approximations $\omegaksi{n}$, we may collect orbits according to the nature of the basepoints, or as \textit{piecewise Fagnano orbits}, \textit{$\cantor$-orbits} and \textit{generalized piecewise Fagnano orbits}.  Note that this does not determine equivalence classes of orbits in the direction $\pi/3$.  Indeed, if $\xio{n}$ and $\yio{n}$ are two elements of a side $s_{n,k}$ of $\omegaksi{n}$, then $\orbitksixNOiang{n}{x}{\theta(\pi/3)}$ and $\orbitksixNOiang{n}{y}{\theta(\pi/3)}$ will be equivalent orbits, i.e., they will have exactly the same length.\footnote{See \S\ref{subsec:mathBillandBillMap} for a more thorough discussion of the equivalence relation on orbits.}  It is only when $\xio{n}$ and $\yio{n}$ differ in their nature, e.g., $\xio{n}\in \cantor_{n,k}'$ and $\yio{n}\notin \cantor_{n,k}$, that we see a difference in the orbits in each compatible sequence given to us by $\{\orbitksiang{i}{\pi/3}\}_{i=0}^\infty$ and $\{\orbitksixang{i}{y}{\pi/3}\}_{i=0}^\infty$.  The significance of this fact will become apparent in \S 5 where it will aid us in giving a description of periodic orbits in the direction of $\pi/3$ of the Koch snowflake billiard $\omegaks$.

\begin{remark}
Let $n\geq 1$.  If $\orbitksiang{n}{\pi/3}$ and $\orbitksix{n}{y}$ are two periodic orbits of $\omegaksi{n}$ with the same period (that is, $\#\orbitksiang{n}{\pi/3}=\#\orbitksixang{n}{y}{\pi/3}$) such that, for every $k_n\leq \#\orbitksi{n}$, $x_n^{k_n}$ and $y_n^{k_n}$ lie on the same side of $\omegaksi{n}$, then the corresponding paths traced out by connecting consecutive basepoints have exactly the same length.  Such a fact follows from the known equivalence between the billiard flow on the rational billiard $\omegaksi{n}$ and the geodesic flow on the associated prefractal flat surface $\mathcal{S}(KS_n)$; see the discussion at the very end of \S\ref{subsec:flatStructuresandFlatSurfaces}.
\label{rmk:sameAddressesSameLengths}
\end{remark}

Let $n\geq 1$.  If $\xooi{i}$ denotes the $i$th character in the ternary expansion of $\xoo$ (in terms of $0,1,2$), then define $\omega_n(\xoo)$ to be the cardinality of the set $\{\xooi{i}\,|\, \xooi{i} = 1, \, 1\leq i\leq n\}$. 

\begin{proposition}
\label{prop:nearByPForbits}
Let $n\geq 0$, $k\leq 3\cdot 4^n$ and $\xio{n}\in s_{n,k}$.  Then there exist $m\leq n$, $h\leq 3\cdot 4^m$ and $\yio{m}\in s_{m,h}$ such that $\orbitksixang{m}{y}{\pi/3}$ is a piecewise Fagnano orbit of $\omegaksi{m}$ and $\yio{m}$ compatible with some $\yoo\in M(\cantor)$ with the following being true: there exists $(\yii{m}{k_m},\theta(\pi/3)\in \orbitksixang{m}{y}{\pi/3}$ such that $\xio{n}$ and $\yii{m}{k_m}$ both lie on the same side $s_{m,h'}$ of $\omegaksi{m}$, $h'\leq 3\cdot 4^m$, and the ternary representation of $\xio{n}$ and $\yii{m}{k_m}$ \emph{(}relative to $s_{m,h'}$\emph{)} are identical up to the first $m-1$ characters.
%Let $n\geq 0$ and $k\leq 3\cdot 4^n$.  If $\xio{n}\in s_{n,k}$, then there exist $m\leq n$, $h\leq 3\cdot 4^m$  and $\yio{m}\in s_{m,h}$ such that $\orbitksixang{m}{y}{\pi/3}$ is a piecewise Fagnano orbit of $\omegaksi{m}$ \emph{(}that is, $\orbitksixang{n'}{y}{\pi/3}$ is compatible with an orbit $\orbitksixang{0}{y}{\pi/3}$\emph{)} and, for some $k_{m}\in \mathbb{Z}$ such that $\yii{m}{k_{m}}$ is a basepoint of an element $(\yii{m}{k_{m}},\theta(\pi/3))\in\orbitksixang{m}{y}{\pi/3}$, and $\xio{n}\in s_{m,h}$ such that $\xio{m}$ and $\yii{m}{k_{m}}$ have ternary expansions  \emph{(}relative to $s_{m,h}${)} that are identical up to the $n'-1$th digit.
\end{proposition}

\begin{proof}
We defer the proof of this statement until the very end of  \S\ref{subsubsec:straighteningAddresses}, since much of the machinery introduced in \S5 makes proving this proposition considerably less tedious.
\end{proof}

\begin{theorem}[Computation of the period of $\orbitksiang{n}{\pi/3}$]
\label{thm:periodOfOrbitksi}
Let $\compseqang{\pi/3}$ be a compatible sequence of periodic orbits. Then every orbit $\orbitksiang{n}{\pi/3}$ in the compatible sequence has a period $\#\orbitksiang{n}{\pi/3}$ that is determined by the number $\omega_n(\xoo)$. Specifically, for every $n\geq 0$, the orbit $\#\orbitksiang{n}{\pi/3}$ is given by the following formula:

%The period of a periodic orbit $\orbitksiang{n}{\pi/3}$ of $\omegaksi{n}$ is determined by the number \emph{(}denoted by $\omega_n(\xoo)$\emph{)} of characters $1$ in the ternary expansion \emph{(}in terms of the characters $0,1,2$\emph{)} of the basepoint $\xoo$ of the compatible orbit $\orbitksiang{0}{\pi/3}$.  Specifically, the period of $\orbitksiang{n}{\pi/3}$ \emph{(}denoted by $\#\orbitksiang{n}{\pi/3}$\emph{)} is given by

\begin{equation}
\#\orbitksiang{n}{\pi/3} = 3\cdot 2^{\omega_n(\xoo)}.
\end{equation} 
\end{theorem}

\begin{proof}
Let $\compseqang{\pi/3}$ be a compatible sequence of periodic orbits.  For the purpose of readibility, we consider the following three cases:

\begin{enumerate}
\item{$\compseqang{\pi/3}$ is a compatible sequence of piecewise Fagnano orbits.}
\item{$\compseqang{\pi/3}$ is an eventually constant compatible sequence of periodic orbits.}
\item{$\compseqang{\pi/3}$ is a compatible sequence of generalized piecewise Fagnano orbits.}
\end{enumerate}

\noindent \textbf{Case 1}: Suppose $\compseqang{\pi/3}$ is a compatible sequence of piecewise Fagnano orbits. Then there exists $N\geq 0$ such that for every $n\geq N$, $\orbitksiang{n}{\pi/3}$ is a piecewise Fagnano orbit and for every $n< N$, $\orbitksiang{n}{\pi/3}$ is not a piecewise Fagnano orbit.  Let $n = 0$.  Then $\orbitksixang{n}{x}{\pi/3}$ has period $\#\orbitksiang{n}{\pi/3}=3$ if the orbit is the Fagnano orbit of $\omegaksi{0}$, and $\#\orbitksiang{n}{\pi/3} =6$ otherwise.

Now suppose that $n = 1$.  Without loss of generality, we may assume $\xoo\in (1/3,2/3)$.  Therefore, $\omega_1(\xoo) = 1$.  In addition, we see that $\#\orbitksiang{n}{\pi/3} = 6 = 3\cdot 2^{\omega_1(\xoo)}$. 

Next, let $M\geq 0$.  Suppose that for every $n\leq M$, we have $\#\orbitksiang{n}{\pi/3} = 3\cdot 2^{\omega_n(\xoo)}$.  We may as well assume that $N<M$.  Therefore, $\orbitksiang{M}{\pi/3}$ is a piecewise Fagnano orbit of $\omegaksi{M}$.  By definition, $\pfix{M+1}{\xoo}$ is constructed from $\pfix{M}{\xoo}$ by appending $3\cdot 2^{\omega_M(\xoo)}$ many scaled copies of $\ppfi{0}$ (the Fagnano orbit of the equilateral triangle billiard $\omegaksi{0}$) to each basepoint of $\pfix{M}{\xoo}$.  Hence, 

\begin{align}
\#\pfix{M+1}{\xoo} = 2\cdot \#\pfix{M}{\xoo} &= 2\cdot 3\cdot 2^{\omega_{M}(\xoo)} = 3\cdot 2^{\omega_{M}(\xoo) + 1}.
\end{align}

\noindent Since $(\xoo)_n=1$ for every $n\geq N$, we deduce that $\omega_{M+1}(\xoo) = \omega_M(\xoo)+1$, and the result follows for a compatible sequence of piecewise Fagnano orbits.

\vspace{2 mm}
\noindent\textbf{Case 2}: Suppose $\compseqang{\pi/3}$ is an eventually constant compatible sequence of periodic orbits.  Then there exists $N\geq 0$ such that for every $n\geq N$, $\orbitksiang{n}{\pi/3}$ is identical to $\orbitksiang{N}{\pi/3}$.  Moreover, $\omega_n(\xoo)=\omega_N(\xoo)$ for every $n\geq N$.  Suppose $N=0$.  Then $\#\orbitksiang{N}{\pi/3} = 6$.  Suppose $N>0$.   By Proposition \ref{prop:nearByPForbits}, there exist $M\leq N$ and $\yoo\in I$ such that  $\pfix{M}{\yoo}$ forms a piecewise Fagnano orbit of $\omegaksi{M}$, $\omega_N(\xoo) = \omega_M(\yoo)$ and $\#\orbitksiang{N}{\pi/3} = \#\pfix{M}{\yoo}$.  In Case 1, we saw that the period of a piecewise Fagnano orbit was determined by the formula $\#\pfix{M}{\yoo} = 3\cdot 2^{\omega_M(\yoo)}$.  Therefore,

\begin{align}
\#\orbitksiang{N}{\pi/3} = \#\pfix{M}{\yoo}  &= 3\cdot 2^{\omega_M(\yoo)} = 3\cdot 2^{\omega_N(\xoo)}.
\end{align}

\noindent Since $\orbitksiang{n}{\pi/3}$ is identical to $\orbitksiang{N}{\pi/3}$, the result follows for an eventually compatible sequence of periodic orbits.

\vspace{2 mm}
\noindent\textbf{Case 3}: Suppose $\compseqang{\pi/3}$ is a compatible sequence of generalized piecewise Fagnano orbits.  For each $n\geq 0$, there exists $\yoo\in I$ and $M\leq n$ such that $\pfix{M}{\yoo}$ is a piecewise Fagnano orbit of $\omegaksi{M}$ and $\omega_n(\xoo) = \omega_M(\yoo)$.  

Let $N\geq 0$ and suppose that for every $n \leq N$, $\#\orbitksiang{n}{\pi/3} = 3\cdot 2^{\omega_n(\xoo)}$.  Then, there exist $\yoo\in I$ and $M\leq N+1$ such that $\pfix{M}{\yoo}$ is a piecewise Fagnano orbit of $\omegaksi{M}$ and $\omega_{N+1}(\xoo) = \omega_M(\yoo)$. Therefore, 

\begin{align}
\#\orbitksiang{N+1}{\pi/3} = \#\pfix{M}{\yoo}  &= 3\cdot 2^{\omega_M(\yoo)} = 3\cdot 2^{\omega_{N+1}(\xoo)}.
\end{align}

This concludes the proof of Theorem \ref{thm:periodOfOrbitksi}.
\end{proof}

\begin{notation}
We now explain the notation that is about to be used in the following theorem and proof.  The characteristic function $\chi$ is defined on the space of characters $\{0,1,2\}$ and is given by

\begin{align}
\chi[\alpha] &:= \left\{\begin{array}{ll} 0 & \text{if } \alpha = 0,2\\ 1 & \text{if } \alpha = 1; \end{array}\right.
\label{eqn:charFunctionOnOneC}
\end{align}
\noindent in other words, it is the characteristic function of $\{1\}$.  A ternary expansion of an element $\xoo\in I$ begins with $0.\cdots$.  Consequently, in order to simplify the notation, when we discuss $\chi[(\xoo)_n]$, we are viewing the ternary expansion of $\xoo$ as the sequence of characters occurring to the right of the decimal point \emph{(}and we no longer indicate the subscript $3$ in the expansion\emph{)}.
\end{notation}

\begin{theorem}[Length of the billiard ball path corresponding to $\orbitksiang{n}{\pi/3}$]
\label{thm:lengthBilliardBallPath}
If $\compseqang{\pi/3}$ is a compatible sequence of periodic orbits and $\FL$ is the length of the Fagnano orbit of the equilateral triangle billiard $\omegaksi{0}$ \emph{(}i.e., the shortest inscribed polygon in $\omegaksi{0}$\emph{)}, then, for each $n\geq 0$, the length $|\orbitksiang{n}{\pi/3}|$ of the path traced out by connecting consecutive basepoints of the orbit $\orbitksiang{n}{\pi/3}$ is given by

\begin{equation}
|\orbitksiang{n}{\pi/3}| = 2\mathscr{L} + \sum_{i=2}^n \chi\left[\xooi{i}\right]\#\orbitksiang{i-1}{\pi/3}\frac{\mathscr{L}}{3^i},
\end{equation}

\noindent where, as before, $\xooi{i}$ denotes the $i$th character in the ternary expansion representing $\xoo$.
\end{theorem}

\begin{proof}
Let $n=0$.  If $\overline{1}$ is the ternary expansion (in terms of the characters $0,1,2$) of an element $\xoo\in I$, where $I$ is viewed as the base of $KS_0$,  then $\orbitksiang{0}{\pi/3}$ is the Fagnano orbit of $\omegaksi{0}$ and has length $|\orbitksiang{0}{\pi/3}| = \FL$.  If $\xoo$ has a ternary representation different from $\overline{1}$, but beginning in the character $1$, then $|\orbitksiang{0}{\pi/3}|=2\FL$. In either case, if $x_1^0$ is collinear with $\xoo$, then $|\orbitksiang{1}{\pi/3}| = 2\FL$.

Consider the basic case $n=2$.  Let $x_2^0 \in KS_2$ be collinear with $\xoo\in I$ (viewed as the base of the equilateral triangle $KS_0$).  Then, $x_2^1$ is the basepoint in the compatible orbit $\orbitksiang{2}{\pi/3}$.  Let $\xooi{i}$ be the $i$th character in the ternary expansion of $\xoo$.  We want to show that 

\begin{equation}
|\orbitksiang{2}{\pi/3}| = 2\FL + \sumpieces{2}.
\label{eqn:caseNisTwo}
\end{equation}

\noindent If $\xooi{2} = 1$, then there exists $\yoo\in M(\cantor)$ such that $|\orbitksiang{2}{\pi/3}| = |\pfix{2}{\yoo}|$ and $\xooi{i}=\yooi{i}$ for every $i\leq n=2$.  There exists $y_1^0 \in KS_1$ that is collinear with $\yoo$ such that $|\pfix{2}{\yoo}|=|\orbitksix{1}{y}|+\#\orbitksix{1}{y}\frac{\mathscr{L}}{3^2}$, since $\pfix{2}{\yoo}$ is a piecewise Fagnano orbit of $KS_2$.  Since $\yoo$ and $\xoo$ have the same first two characters in their respective ternary expansions, we have that

\begin{align}
\label{eqn:baseCaseOfLengthFormula}
\notag |\orbitksix{1}{y}| &= |\orbitksix{1}{x}|
\end{align}
\noindent and
\begin{align}
\notag \#\orbitksix{1}{y} &= \#\orbitksix{1}{x}.
\end{align}

Since $\orbitksix{1}{y} = 2\FL$ (which is independent of the choice of basepoint, so long as such a choice is not a corner of the billiard table $\omegaksi{1}$), it follows that 

\begin{equation}
|\orbitksiang{2}{\pi/3}| = 2\FL +\sum_{i=2}^2 \#\orbitksiang{i-1}{\pi/3}\frac{\FL}{3^i}. 
\end{equation}

\noindent If $\xooi{2} \neq 1$, then $\chi[\xooi{2}]=0$ and 

\begin{align}
|\orbitksiang{2}{\pi/3}| &= |\orbitksiang{1}{\pi/3}| = 2\mathscr{L}.
\end{align}

\noindent In either case, we have shown that Equation (\ref{eqn:caseNisTwo}) holds.

Let us now proceed by induction and fix $N\geq 2$.  Suppose that for every $n\leq N$, 

\begin{equation}
|\orbitksiang{n}{\pi/3}| = 2\FL + \sumpieces{n}.
\label{eqn:inductionHypothesisLength}
\end{equation}

\noindent  Then there exists $\yoo\in M(\cantor)$ and $M\leq N+1$ such that for all $i\leq M$, $\xooi{i}=\yooi{i}$ and

\begin{equation}
|\orbitksiang{N+1}{\pi/3}| = |\pfix{M}{\yoo}|.
\end{equation}

\noindent If $M=N+1$, then the nature of $\pfix{N}{\yoo}$ dictates that 

\begin{equation}
|\pfix{N+1}{\yoo}|=|\orbitksix{N}{y}| + \#\orbitksix{N}{y}\frac{\FL}{3^{N+1}}.
\end{equation}

\noindent Applying this fact, the induction hypothesis (\ref{eqn:inductionHypothesisLength}) and Remark \ref{rmk:sameAddressesSameLengths}, we have that 

\begin{eqnarray}
|\orbitksiang{N+1}{\pi/3}| \notag&=& |\pfix{N+1}{\yoo}| \\
								 \notag&=& |\orbitksix{N}{y}| + \#\orbitksix{N}{y}\frac{\FL}{3^{N+1}}\\
								 \notag&=& |\orbitksix{N}{x}| + \#\orbitksix{N}{x}\frac{\FL}{3^{N+1}}\\
								 &=& 2\FL +\sumpieces{N}+ \\
								 \notag& & \quad \#\orbitksix{N}{x}\frac{\FL}{3^{N+1}}\\
								 \notag&=& 2\FL +\sumpieces{N+1}.
\end{eqnarray}

\noindent If $M<N+1$, then, again applying the previously mentioned fact, the induction hypothesis and Remark \ref{rmk:sameAddressesSameLengths}, we deduce that

\begin{eqnarray}
\notag |\orbitksiang{N+1}{\pi/3}| &=& |\pfix{M}{\yoo}|\\
\notag				 &=& |\orbitksix{M-1}{y}| + \#\orbitksix{M-1}{y}\frac{\FL}{3^{M}}\\
\notag				 &=& |\orbitksix{M-1}{x}| + \#\orbitksix{M-1}{x}\frac{\FL}{3^{M}}\\
							 &=& 2\FL + \sumpieces{M-1} + \label{eqn:lengthofOrbitKSn}\\
\notag				 & & \quad  \#\orbitksix{M-1}{x}\frac{\FL}{3^{M}}\\
\notag				 &=& 2\FL + \sumpieces{M}\\
\notag				 &=& 2\FL + \sumpieces{N+1},
\end{eqnarray}

\noindent where the last lines of the calculation in Equation (\ref{eqn:lengthofOrbitKSn}) follow from the fact that the characters $\xooi{i}$, with $M<i\leq N+1$, are necessarily never equal to $1$, meaning that $\chi{\xooi{i}}=0$ for $M<i\leq N+1$.

\end{proof}

\section{Periodic orbits of $\mathbf{\omegaks}$ in the direction of $\pi/3$}
\label{sec:orbitsOfOmegaKS}

In \S 4, we were able to group orbits with initial directions $\theta(\pi/3)$ into particular categories: piecewise Fagnano orbits, $\cantor$-orbits and generalized piecewise Fagnano orbits.  We stressed that this grouping was not equivalent to the classification of orbits described in \S\ref{subsec:mathBillandBillMap}.  Such a grouping is, however, meant to allude to a \textit{description} of orbits in $\omegaks$ with an initial direction of $\pi/3$, where this description is determined by the ternary expansion of an initial basepoint of an orbit $\orbitksiang{0}{\pi/3}$ with an initial direction of $\pi/3$.  In addition, we showed that for every orbit $\orbitksixang{n}{y}{\theta(\pi/3)}$, there existed $\xio{n}\in KS_n$ and $\xoo\in I$ collinear in the direction of $\pi/3$ such that $\orbitksixang{n}{y}{\theta(\pi/3)} = \orbitksiang{n}{\pi/3}$.

We want to somehow overcome the fact that the boundary of $\omegaks$ is nondifferentiable.  When determining orbits of the snowflake, we want to stress the fact that we are primarily interested in the collision points with the boundary $KS$.  In the theory of mathematical billiards, one may consider the billiard orbit to be the path traced out by the billiard ball or just the (equivalence classes of) ordered pairs $(x,\theta)\in (B\times S^1)/\sim$.  Even more simply, one may think of the path traced out by the billiard ball or the basepoints on the boundary $B$ as a set of dynamically ordered points.\footnote{In the theory of dynamical systems, such a set is called the \textit{Poincar\'e section} of the flow; see \S 2.2.} When we want to make the distinction between the path and the ordered pairs in the prefractal billiard tables $\omegaksi{n}$ and in $\omegaks$, we will explicitly refer to such a path as \textit{the billiard ball path} (or, simply, the orbit) corresponding to a particular orbit and to the collection of collision points as the \textit{footprint of the orbit}. In other words, the \textit{footprint} is the intersection of the orbit with the boundary of the billiard table. Strictly speaking, we cannot speak of basepoints here, because this language is indicative of the existence of a well-defined phase space (a unit tangent bundle), which we have yet to rigorously establish.    

As noted in \S\ref{sec:introduction}, in the sequel, we make a slight abuse of language.  A self-similar set is the unique fixed point attractor of a particular iterated function system; see, for example,  Figure \ref{fig:kochcurveconstruction} of \S\ref{sec:introduction}.  By abuse of language, we will also say that a set is self-similar if it is the union of finitely many isometric, abutting copies of a given self-similar set, much as the snowflake curve is the union of three (isometric, abutting) copies of the Koch curve; see Figure \ref{fig:3kochcurvesLabeled} of \S\ref{sec:introduction}.

We next define a ``self-similar orbit.''

\begin{definition}[Self-similar orbit]
\label{def:selfsimorbit}
Let $\gamma$ be a periodic orbit of $\omegaks$.  Then, $\gamma$ is said to be a \textit{self-similar orbit} if its footprint is a self-similar subset of $KS\subseteq \mathbb{R}^2$. 
\end{definition}

%We want to show that it is possible to appropriately recast a compatible sequence of orbits so that a suitable limit of a compatible sequence may describe an orbit of the Koch snowflake billiard.

%Moreover, we will show that certain orbits of the snowflake in the direction $\pi/3$ may be described as a suitable limit of 

We want to show that it may be possible to appropriately recast particular compatible sequences of periodic orbits in such a way that a suitable notion of limit may be applied so as to yield an orbit of the Koch snowflake billiard $\omegaks$.  Specifically, we will provide a plausibility argument as to why a suitable limit of a compatible sequence of piecewise Fagnano orbits constitutes an orbit of the Koch snowflake billiard $\omegaks$.  In addition to this, we will show how an eventually constant compatible sequence of periodic orbits can be thought of as a periodic orbit of $\omegaks$, called a stabilizing orbit (or a $\cantor$-orbit) of $\omegaks$.  We will also discuss a connection between the ternary representation of the initial basepoint of the initial orbit of the eventually constant compatible sequence of orbits and the prefractal approximation $\omegaksi{n}$ at which such a compatible sequence stabilizes.

%\begin{enumerate}
%\item{an eventually constant compatible sequence of periodic orbits with an initial direction $\pi/3$;}
%\item{a compatible sequence of piecewise Fagnano orbits.}
%\item{generalized piecewise Fagnano orbits.}
%\end{enumerate}

%In case (1), the resulting orbit will be called a \textit{stabilizing} periodic orbit (or a $\mathscr{C}$-orbit) of $\omegaks$, while in case (2), it will be called a \textit{piecewise Fagnano} (periodic) orbit of $\omegaks$.  We stress, however, that we do not yet know what a piecewise Fagnano orbit of $\omegaks$ truly is.  Therefore, we will provide a plausible argument as to why a suitable limit of piecewise Fagnano orbits should be called a piecewise Fagnano orbit of  $\omegaks$ and identify what we expect to be the footprints of such orbits.

\subsection{Addressing system for $KS$}
\label{subsec:addressSystemForKS}
We now describe an addressing system for the Koch snowflake $KS$ that will be used in the sequel.  This addressing system was introduced (and used) by M. Pang and the first author in \cite{LaPa} in their study of the boundary behavior of the eigenfunctions of the Koch snowflake drum. 

\label{idx:lcr}
First, we mention that the ternary representation of an element of some side $s_{n,k}$ of $\omegaksi{n}$ will be given in terms of the alphabet $\{l,c,r\}$ (for \textit{left}, \textit{center} and \textit{right}) so as to not confuse the reader later when we discuss the addresses of points of the Koch snowflake curve $KS$ in terms of the alphabet $\{0,1,2,3,4,5\}$.\footnote{In \S4, we used the standard alphabet $\{0,1,2\}$ to describe the ternary representation of an element in $I$.}  In the Koch snowflake curve, there are three types of points, as described below:

\begin{enumerate}
\item{\textit{Corners}: a vertex of the polygon $KS_n$ is a corner.  As expected, such a corner in $KS_n$ is also a corner in $KS_m$ for all $m\geq n$.}
\item{\textit{Non-ternary points}: an element with a ternary representation (relative to a side $s_{n,k}$) not terminating solely in $l$'s or $r$'s.  As such, a \textit{non-ternary point} of $KS$ is an element of some $\cantor_{n,k}'\subseteq s_{n,k}$, where $\cantor_{n,k}'$ was defined in the paragraph preceding Proposition \ref{prop:constantCompSeq}.}
\item{\textit{Elusive limit points}: an element of $KS$ that is never an element of a finite approximation $KS_n$, for any $n\geq 0$.  Though every point of $KS$ is a limit point, we make this abuse of language so as to emphasize the elusive nature of these points which are not vertices of $KS_n$, for any $n$, and are not elements of $\cantor_{n,k}'$, for any $n$ and any $k\leq 3\cdot 4^n$.}
%\item	Corners of some prefractal approximation $\omegaksi{n}$ and every subsequent prefractal $\omegaksi{m}$,  $m \geq n$.
%\item	Elements of $\cantor_{n,k}$, which can be viewed as a subset of a side $s_{n,k}$ of $\omegaksi{n}$, where $k$ is some positive integer such that $k \leq 3\cdot 4^n$.
%\item	Points which are never elements of any prefractal approximation.  At times, we may refer to this third class of points as the \textit{elusive limit points of the Koch snowflake}, which is another abuse of language.
\end{enumerate}

\begin{figure}
\begin{center}
\includegraphics[scale=.5]{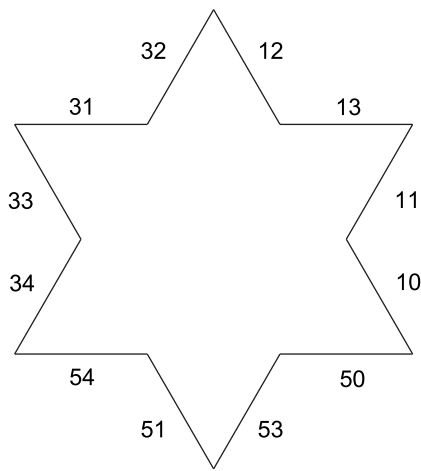}
\end{center}
\caption{The sides of the prefractal billiard $\omegaksi{1}$ labeled according to the particular addressing system described in \S\ref{subsec:addressSystemForKS}.}
\label{fig:sidesLabeled}
\end{figure}

In order to understand how these points are assigned addresses, we need to understand exactly how the addressing system works.  The alphabet from which we will draw upon to construct addresses of points of the snowflake is $\{0,1,2,3,4,5\}$.  We start by noting that every side of a prefractal billiard table may be identified by a finite address.  An address of a side may end in any character 0, 1, 2, 3, or 4 but never 5.  However, an address may start with 5.  If an address of $n$ characters ends in a 1 or a 3, then the side $s_{n,k}$ identified by this address is not a subset of any of the sides of $KS_{n-1}$.  If an address of $n$ characters ends in 0, 2 or 4, then the side $s_{n,k}$ is either the left third or right third interval of a side $s_{n-1,k'}$ of $\omegaksi{n-1}$, and has a length of $1/3^n$.   In Figure \ref{fig:sidesLabeled}, we label all the sides of $\omegaksi{1}$ according to this convention.  

Heuristically, the labeling is such that the next character in an address is decided by `moving forward' or shifting left or right.  When moving forward, there is an option of a `twisting' right or `twisting' left, which results in appending a 1 or 3 to the address, respectively. When making a shift to the left third or the right third of a particular side and the last character of an address is a 1, then the address is appended with a 4 or 2, respectively.  Likewise, having an address end in 3 and performing a shift to the left third or right third of a side results in appending this address with a 2 or 0, respectively.  Note that when `moving forward' into a cell emanating from the side $s_{0,k}$ with the label 5, one must first rotate the directional by $\pi$ so that the proper addresses may be applied to the sides of the cell $C_{n,k}$; see Figure \ref{fig:examplesOfDirectional}.

\begin{figure}
\begin{center}
\includegraphics[scale=.5]{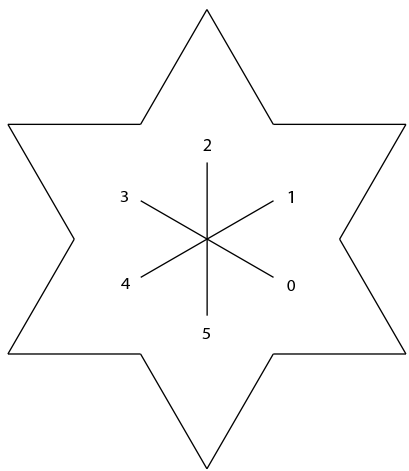}
\end{center}
\caption{The directional used to identify points on the Koch snowflake $KS$.}
\label{fig:directionalInKS1}
\end{figure}

\begin{figure}
\begin{center}
\includegraphics[scale=.4]{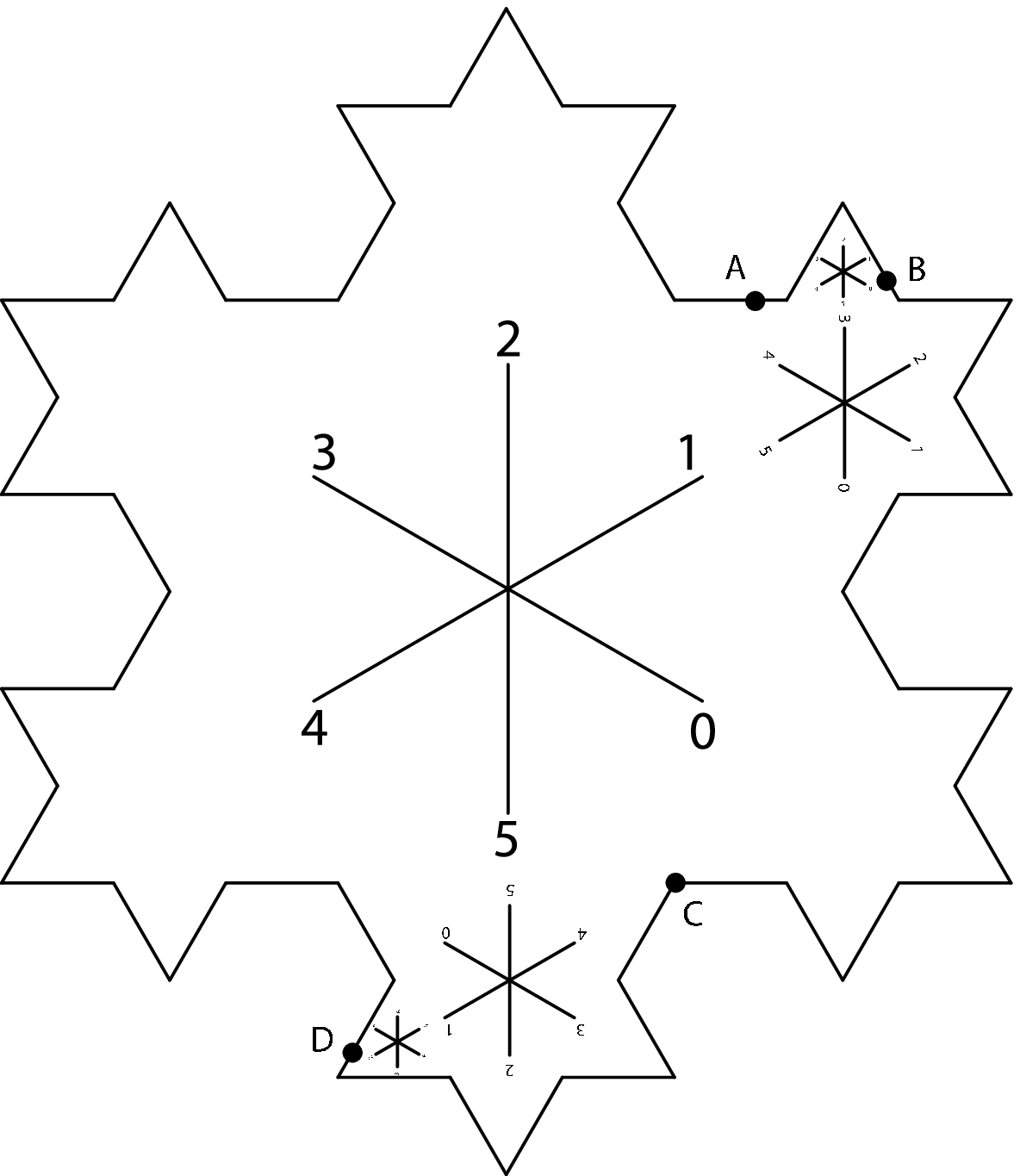}
\end{center}
\caption{Addresses of points of the Koch snowflake are not necessarily unique.  In the context of a finite approximation $KS_n$, the addressing system serves to label sides of the prefractal approximation $KS_n$.  Points $A$, $B$ and $D$ then lie on sides which are labeled by the addresses $134$, $131$, $511$, respectively.  The point $C$, however, does not lie exactly on a side, but at a vertex.  Hence, one could say that $C$ lies on the side with address $504$ or $534$.  Moreover, in the limit, $C$ may be represented by the equivalent addresses $50\overline{4}$ and $53\overline{4}$, where (as before) the overbar indicates that $4$ is repeated ad infinitum.  Since a billiard orbit is a collection of dynamically ordered points, it follows that each basepoint of an element of an orbit can be given an address in terms of this addressing system and still be distinguished from other basepoints of elements of the orbit with the same address.}
\label{fig:examplesOfDirectional}
\end{figure}

In terms of the addressing system described above and in Figures \ref{fig:sidesLabeled}, \ref{fig:directionalInKS1} and \ref{fig:examplesOfDirectional}, the three types of points of $KS$ listed above have very specific symbolic representations:
\begin{enumerate}
\item{Points of the snowflake that are corners (i.e., vertices of some prefractal approximation $KS_n$, with $n\geq 0$) have finite (i.e., terminating) addresses; equivalently, the addresses of these points may be infinite, terminating solely in 0's, 2's or 4's.}
\item{Points that are elements of a set $\cantor_{n,k}'$ when viewed as a subset of some side $s_{n,k}$ of $\omegaksi{n}$  have addresses ending solely in either 0's and 2's or solely in 4's and 2's.}  
\item{Elusive limit points of the Koch snowflake $KS$ have addresses given by (finitely or infinitely many) 0's or 2's and/or 4's or 2's interspersed throughout a word consisting of:
\begin{enumerate}
\item{infinitely many $1$'s and $3$'s,}
\item{infinitely many $1$'s and finitely many $3$'s,}
\item{or infinitely many $3$'s and finitely many $1$'s.}
\end{enumerate}}
\end{enumerate}

\begin{example}
The address $131310002\overline{20}$ identifies an element of $\cantor_{4,k}'$ on the side $13131$ of the prefractal  $\omegaksi{4}$ \emph{(}and hence, a type \emph{(}2\emph{)} point\emph{)}.  Moreover, the address $13111234\overline{3}$ identifies an elusive limit point of the snowflake curve $KS$ \emph{(}and hence, a type \emph{(}3\emph{)} point\emph{)}.  Finally, the address $5\overline{4}$ identifies the point $1/3$ on the base of the equilateral triangle billiard $\omegaksi{0}$; this is a type \emph{(}1\emph{)} point which will end up being an obtuse corner in every subsequent prefractal billiard table. 
\end{example}

\subsection{The symbolic representation of a compatible sequence of initial basepoints}\label{subsec:theSymbolicRepresentationOfACompSequenceOfInitialOrbits}

In \S\ref{sec:orbitsOfOmegaKSn}, we defined what it meant to be a compatible sequence of orbits and a compatible sequence of initial conditions (see Definition \ref{def:compSeqOrbits} and Definition \ref{def:compSeqInitConditions}, respectively, along with Theorem \ref{thm:compOrbitiffCompInitCond}).  Such a construction amounted to picking an initial basepoint $\xoo$ of an orbit $\orbitksiang{0}{\theta_0^0}$ such that the initial basepoint of every subsequent orbit in the compatible sequence  was collinear with $\xoo$ in the same initial direction relative to a fixed coordinate system.  

Now that we have a clearly defined addressing system for points of the Koch snowflake $KS$, we can discuss the symbolic description of a compatible sequence of initial basepoints. In \S\ref{subsec:addressSystemForKS}, we discussed how appending an address with a 1 or 3 essentially amounts to moving forward into a cell.  If $\xoo$ is the midpoint of the base of the equilateral triangle, then the next point in the sequence of compatible initial basepoints is the midpoint of a side of $\omegaksi{1}$ such that $\xio{1}$ and $\xoo$ are collinear.  Symbolically, this element lies on the side $s_{1,k}$ with address $51$.  Since $\xio{1}$ is a midpoint of the side $51$, it follows that $\xio{2}$ lies on the side $513$.  Consequently, the sequence of points in the compatible sequence of initial basepoints must converge to a point on the snowflake given by the address $5\overline{13}$.  If $\xii{0}{1}$ is the basepoint of the next element of the footprint of the orbit $\orbitksi{0}$, then we see that $\xii{0}{1}$ lies on a side with $\omegaksi{0}$ with address 1.  Since $\xii{0}{1}$ is a midpoint of this side, it follows that the basepoint $\xii{1}{1}$ of the orbit $\orbitksi{1}$ lies on the side with address $13$.  Consequently, the sequence of basepoints $\{\xii{i}{1}\}_{i=0}^\infty$ is converging to a point of the snowflake with an address $\overline{13}$. 

In general, if $\xoo$ is an initial basepoint of an orbit in $\omegaksi{0}$, then the address of the side $s_{n,k}$ on which $\xio{n}$ resides is determined in a similar manner.  Similarly, the address of the side $s_{n,k'}$ containing the basepoint of the image $f_n(x_n^0,\theta_n^0)$ under the billiard map $f_n$ associated with $\omegaksi{n}$ is also determined in such a way; recall that $f_n$ is a map from the phase space $(KS_n\times S^1)/\sim$ to itself and see \S\ref{subsec:mathBillandBillMap} for the definition of the billiard map for a general polygonal billiard $\Omega(B)$.  For example, suppose the unit interval base $I$ of $\omegaksi{0}$ was partitioned into intervals of length $1/3^n$, for some $n\geq 0$.  Then suppose $\xoo$ was in one of these intervals.  It follows that $\xii{0}{1}$ would be in a similar interval of length $1/3^n$ if the initial direction were $\pi/3$.  If this interval constitutes a ghost of some side of the prefractal approximation $\omegaksi{n}$, then such a ghost is removed and $\xii{n}{1}$ is on some side with an address ending in either $1$ or $3$.

Now that we have symbolically described a compatible sequence of initial basepoints, we want to show that such a symbolic construction is entirely determined by the ternary representation of $\xoo\in I$. Each element of the unit interval $I$ has a ternary expansion.  We represent such expansions not in terms of $0$, $1$ and $2$ but in terms of $l$, $c$ and $r$, for \textit{left}, \textit{center} and \textit{right}, respectively.\footnote{In \S 4, the use of a sequence $\{p_{j_i}\}_{i=1}^N$, with $N\leq \infty$ (and $j_i\in\{1,2,3\}$ for each $i\geq 1$), consisting of either 0's, 1's and/or 2's in order to express the decimal value of elements in the unit interval, was meant to allude to the ternary expansion of a number and the representation of such a number via an address, given here in terms of $l$, $c$ and $r$, instead of 0's, 1's and 2's.}   Such a coding system will prevent us from confusing addresses of points in $I$ with addresses of points of the snowflake, but we will see that there is a very intimate connection between these points and a certain subset of points of the snowflake.  If a ternary expansion of an element of the unit interval $I$ contains an $l$, $c$ or $r$ in the $(n+1)$th position, then it follows that the element resides in either the left, middle or right third of an interval of length $1/3^{n}$ somewhere along the base of the equilateral triangle $KS_0$, respectively. 

%In \S 4 an element $x_0$ of $I$, viewed as the base of $KS_0$, constituted the initial basepoint of an orbit with an initial direction $\pi/3$.  Such an orbit $\orbitksiang{n}{\pi/3}$ was the first element of a compatible sequence of periodic orbits $\compseq$ so long as $x_0$ was not an element with a finite ternary expansion.  We next use this fact to describe initial conditions of orbits of the Koch snowflake. 

Consider the equilateral triangle with an inscribed similar triangle sharing the right vertex of $KS_0$.  The side of the scaled, inscribed copy of $KS_0$ not common to the base and right edge of $KS_0$ constitutes a segment of a billiard ball path in $\omegaksi{0}$, if the initial direction of the ball was $\pi/3$. Consequently, such a leg of a triangle identifies a point $\xii{0}{1}$ with the same ternary \textit{representation} as $\xoo$ on $s_{0,k}$ such that $\xii{0}{1}$ is the basepoint of $f_0(\xoo,\pi/3)$, where $f_0$ is the billiard map associated with $\omegaksi{0}$, namely, a map from $(KS_0\times S^1)/\sim$ to itself. For example, if $x_0$ is the point $3/4$ in the base of $\omegaksi{0}$, then $x_0'$ will also be the point $3/4$, measured relative to the top vertex of $KS_0$.  

We have already seen (in \S\ref{subsec:compSeqPF}) that the point $\overline{c}$ (or $1/2$) is the first element in the compatible sequence of initial basepoints giving rise to \textit{the} compatible sequence of primary piecewise Fagnano orbits (and, therefore, each orbit has an initial angle of $\pi/3$).  If $\xoo$ is some other point with a ternary expansion not equal to $\overline{c}$, then how does one determine the address of the basepoint $\xii{n}{1}$ of the image $f_n(\xio{n},\theta_n^0)$ from the ternary representation of $\xoo\in I$?  As noted above, shifting along the base of the equilateral triangle results in an identical shift along any side, so long as we are considering a billiard ball moving in the direction $\pi/3$. Let us suppose that $\xoo\in I$ has a ternary representation with  $c$ in the $(n+1)$th position.  Then $\xii{0}{1}$ will also have the same ternary representation and a $c$ in the $(n+1)$th position.  This means that both points live on ghosts $g_{0,k}$ and $g_{0,k'}$, respectively, of $\omegaksi{0}$ and will not be basepoints of orbits in subsequent approximations.  This also implies that the subsequent basepoints will have addresses (in terms of the alphabet $\{0,1,2,3,4,5\}$) such that a $1$ or a $3$ occurs in the $(n+1)$th position.  Likewise, an $l$ or an $r$ in the $(n+1)$th position of the ternary representation for $\xoo$ dictates whether or not the character in the $(n+1)$th position of a subsequent initial basepoint is going to have a $0$, $2$ or $4$.  Thus, we determine the address of sides on which initial basepoints reside by examining the address of $\xoo$ in the unit interval.

%\begin{example}
%Something to appear soon
%\end{example}
\vspace{6 pt}
\subsubsection{Straightening addresses of the Koch snowflake}
\label{subsubsec:straighteningAddresses}
Let $x$ be an elusive limit point of the Koch snowflake.  Then $x$ has an address that is given by a sequence of $1$'s and $3$'s, with $0$'s and $2$'s and/or 2's and $4$'s interspersed, as described above.  Without loss of generality, we may assume the address begins with a $1$.  At some point later in the address, there must be a character that is neither a $2$ nor a $0$.  This character will be either a $3$ or a $1$.  If it is a $3$, we continue on to the next character of the address. If it is a $1$, then we may switch this $1$ for a $3$ and until the next $1$ or $3$, switch every $4$ for a $0$ while keeping the $2$'s unchanged.  Continuing this process, making sure that every subsequent $1$ is followed eventually by a $3$ and not a $1$, we end up with what we call the \textit{straightening of $x$}, which we denote by $s(x)$.  Since $s(x)$ has an address consisting of $1$'s and $3$'s interspersed with $0$'s and $2$'s and/or $2$'s and $4$'s such that the $1$'s and $3$'s alternate, we see that this point is collinear with a point of $I$ along the base of the equilateral triangle and such a line connecting them contains no other points of the snowflake.

\begin{example}
Let $n$ and $k$ be positive integers.  Consider the address of a side $s_{n,k}$ of $\omegaksi{n}$ given by $\alpha = 13123232113133100324$.  The straightening of this address is $13123212313131344120$.  Next, suppose the address of a particular point is given by $\overline{\alpha} =  \overline{13123232113133100324}$.  Then $s(x) = 13123212313131344120$.
\end{example}

Recall from \S\ref{subsec:compSeqPF} that $M(\cantor)$ is dense in $I$ and that $I$ is the unique fixed point attractor of the iterated function system $\{\phi_1,\phi_2,\phi_3\}$  described in Equation (\ref{eqn:IFS}).\footnote{Recall that $M(\cantor)$ is defined in Equation (\ref{eqn:MCantor}).} Consequently, $\overline{M(\cantor)}$ (the closure of $M(\cantor)$ in $\mathbb{R}$) is the unit interval $I$.  Define $T$ as the collection of all points of the unit interval $I$ with finite ternary expansions (or equivalently, expansions terminating solely in $l$'s or $r$'s).  Thus, $\overline{M(\cantor)}\setminus T = I\setminus T$.\footnote{We can then write $T^c$ (the complement of $T$ in $I$) when discussing the elements of $M(\cantor)\setminus T$.  However, so as to place more significance on the fact that such elements are limit points of $M(\cantor)$, we will continue using the more verbose notation.} Then $\overline{M(\cantor)}\setminus T$ consists of points which have infinite ternary expansions not terminating solely in $l$'s or $r$'s. 

%As we did in \S\ref{subsec:compSeqPF},

If $\xoo$ were an element of $T$ and $f_0(\xoo,\pi/3) = (\xii{0}{1},\theta_0^1)$, then $\xii{0}{1}$ would also be a point with a finite ternary expansion (or an infinite ternary expansion terminating solely in either 0's or 2's).  From this we deduce that the line segment connecting $\xoo$ and $\xii{0}{1}$ in the direction $\pi/3$ would form a saddle connection in some prefractal billiard $\omegaksi{n}$, for some $n>0$.  So, when we consider allowable points from which to start in the direction $\pi/3$, we consider $\overline{M(\cantor)}\setminus T$.  If $x_0$ is an element of this set, then the compatible sequence of periodic orbits $\{\orbitksi{i}\}_{i=0}^\infty$ is well defined.  

\begin{remark}
Many of the \textit{elusive limit points} of the Koch snowflake may be reached via a piecewise linear logarithmic spiral.  As the name suggests, such spirals may be \textit{straightened} to straight line segments that begin in the interior of the snowflake and end on the boundary of the snowflake, thus identifying an initial basepoint $\xoo$ of an orbit of the billiard $\omegaksi{0}$ and the basepoint $\xii{0}{1}$ of $f_0(\xoo,\pi/3)$.
\end{remark}

It is clear that the addressing system and the directional exploit the inherent symmetry  of the Koch snowflake and its prefractal approximations.  In \S\ref{subsec:aGeneralDiscussionOfOrbits}, we were implicitly using the underlying symmetry of the prefractal approximation $\omegaksi{n}$ as an aid in proving Theorems \ref{thm:periodOfOrbitksi} and \ref{thm:lengthBilliardBallPath}.  In particular, we deferred the proof of Proposition \ref{prop:nearByPForbits} until now, because we wanted to be able to straighten a finite address and to explain why every basepoint of an orbit must have an identical ternary representation, relative to their respective sides.\footnote{Recall that Proposition \ref{prop:nearByPForbits} (from \S\ref{subsec:aGeneralDiscussionOfOrbits}) was used in the proof of Theorem \ref{thm:periodOfOrbitksi} to calculate the periods of certain periodic orbits of $\omegaksi{n}$.}

\begin{proof}[Proof of Proposition \ref{prop:nearByPForbits}]
Let $n\geq 0$ and $k\leq 3\cdot 4^n$.  The point $\xii{n}{1}$ on a side $s_{n,k}$ has a finite address (where this address, in actuality, is identifying the side $s_{n,k}$).  We know that this address can be straightened so as to identify another side $s_{n,k'}$ of $\omegaksi{n}$ such that this side is visible from the interior of $\omegaksi{0}$.  In addition to this, one may uniquely identify $\xii{n}{1}$ in terms of the address given in terms of $\{0,1,2,3,4,5\}$ and its ternary representation given in terms of $\{l,c,r\}$.  Then, utilizing the local and global symmetry, one can straighten $\xii{n}{1}$ to an element of $s_{n,k'}$.  Without loss of generality, we may then suppose that the point on $s_{n,k'}$ is $\xio{n}$. Since the side $s_{n,k'}$ is visible from the interior of $\omegaksi{0}$, so is the midpoint of $s_{n,k'}$.  Given the nature of the construction of the Koch snowflake, there exist $0\leq m\leq n$ and $h\leq 3\cdot 4^m$ such that $s_{n,k'}$ is a subset of a side $s_{m,h}$ of $\omegaksi{m}$.  Contained in $s_{m,h}$ is a scaled (or scaled and rotated) copy of $M(\cantor)$, which we denote by $M(\cantor)_{m,h}$.  Since every point of $M(\cantor)_{m,h'}$ is visible from the interior of $\omegaksi{0}$, it follows that the midpoint $\yii{m}{1}$ of the side $s_{m,h}$ is visible from the interior of $\omegaksi{0}$, as well.  We can now construct a piecewise Fagnano orbit $\orbitksixang{m}{y}{\pi/3} = \pfix{m}{\yoo}$, where the ternary representation of $\xii{n}{1}$ and $\yii{m}{1}$ have the same ternary representation up to the first $m$ many characters.  This construction is precisely the content of Proposition \ref{prop:nearByPForbits}.
\end{proof}

\subsection{Footprints of piecewise Fagnano orbits of $\mathbf{\omegaks}$}
\label{subsec:footprintsOfPeriodicOrbitsOfTheKochSnowflakeBilliard}
From \S\ref{subsec:theSymbolicRepresentationOfACompSequenceOfInitialOrbits}, we now understand how to identify points of the boundary $KS$ and to determine the type of compatible sequence of periodic orbits from the initial condition $(\xoo, \pi/3)$ of $\orbitksiang{0}{\pi/3}$ in $\omegaksi{0}$. 

Let $X=\{0,1,2,3,4,5\}$ and $X_n:= \prod_{i=1}^n X=X^n$, the space of all words of length $n\geq 1$ with characters (or symbols) in the alphabet $X$; by default, we let $X^0:=\emptyset$.  Define $X_\infty:= \prod_{i=1}^\infty X=X^\mathbb{N}$ to be the space of infinite words expressed in terms of the elements of $X$.  Let $x\in X_\infty$ and $n\in\mathbb{N}$.  We define $\tau_n:X_\infty \to X_n$ by $\tau_n(x):=x\vert_{n}$, where $x\vert_n\in X_n$ is the finite word of length $n$ consisting of the first $n$ characters of the infinite word $x$. (Hence, $x\vert_{n}\in X_n$ is the \textit{truncation at level} $n$ of $x\in X_\infty$.) 

%Essentially, $\tau_n$ truncates $x$ to the first $n$ characters.  

%\begin{remark}
%Throughout this article, we have been indexing the prefractal approximations $KS_n$ by the nonnegative integers $\mathbb{N}_0$.  So as to form inverse limits of compatible sequences of orbits, we form the product $\prod_{i=N}^\infty \orbitksi{i-1}$, for some suitable integer $N\geq 1$, so that an inverse limit sequence of orbits begins with the orbit $\omegaksi{N-1}$.  
%\end{remark}

Restricting $\tau_n$ to $KS$ results in a map that is well defined on the points of $KS$ which are not corners of any finite approximation $KS_i$, for any $i\geq 0$.  For example, $1\overline{2}$ and $3\overline{2}$ identify the same point in $KS$, but the words of length $n$ given by $12...2$ and $32...2$ do not identify the same segments of $KS_n$.

Since the Koch snowflake curve $KS$ may be viewed as the inverse limit of its prefractal approximations $KS_n$, when attempting to define orbits of the Koch snowflake, we take as the suitable notion of limit that we have alluded to throughout this paper and \cite{LaNie1} to be the  inverse limit.\footnote{We note that the notion of inverse limit was used in a different context, but also in relation with fractals, in the recent works \cite{BeDeMiSt, RoSt, St} (among others).}

Recall that the billiard ball path is the orbit and the footprint $\mathcal{F}_n(\xoo)$ of the orbit $\orbitksiang{n}{\pi/3}$ is the intersection of the orbit $\orbitksiang{n}{\pi/3}$ with the boundary $KS_n$.  In \S\ref{subsec:plausibilityArgument}, we will provide a plausibility argument as to why it may be possible to recast a compatible sequence of piecewise Fagnano orbits as an inverse limit sequence of orbits.  For now, we show that the inverse limit of the footprints $\mathcal{F}_i(\xoo)$ exists.  Such a task will require us to define the proper transition maps so that the definition of inverse limit is satisfied.  (See \S\ref{subsec:invLimSeqInvLim} above for a brief discussion of inverse limits, and \cite{Bo,HoYo} for further information.) In other words, we recast the sequence of footprints $\{\mathcal{F}_i(\xoo)\}_{i=0}^\infty$ corresponding to a compatible sequence of orbits as an inverse limit sequence of footprints of periodic orbits, where each basepoint is given in terms of its address; see \S\ref{subsec:addressSystemForKS} (including Figures \ref{fig:sidesLabeled}, \ref{fig:directionalInKS1} and \ref{fig:examplesOfDirectional}) for a description of the addressing system. Once (and if) an appropriate billiard flow can be defined, we will seek in future work to show that the inverse limit of footprints is in some way an analogue of a Poincar\'e section of the suitably defined billiard flow.

Given $m,n\in \mathbb{N}^*$ with $m\leq n$, consider the map $\tau_{nm}:X_n\to X_m$, where $\tau_{nm}(x_n) = x_{n}\vert_{m}$ is the truncation of the word $x_n$ (of length $n$) by $n-m$ characters (thereby producing a word of length $m$).  Then $\tau_{nm}\vert_{KS_n}$ is a map that truncates the addresses of segments of $KS_n$ of length $n+1$ to produce addresses of the segments of $KS_m$ of length $m+1$.  Let the map $\iota = \iota_{mn}$ be the `identity map' acting on the unit circle $S^1$.  The map $\iota = \iota_{mn}$ will serve to preserve the compatibility of the two initial conditions $(\xio{m},\theta_{m}^0)$ and $(\xio{n},\theta_n^0)$. In actuality, $\iota(\theta_n^0) = \theta_m^0$, meaning that $(\xio{m},\theta_m^0)$ is identified with $(\xio{m},\theta_n^0)$.\footnote{The map $\iota$ maintains the compatibility of the initial basepoints $\xio{n}$ and $\xio{m}$, by relying on the fact that angles are measured with respect to a fixed coordinate system.}

Consider $\xoo\in I$ with a ternary representation consisting of infinitely many $c$'s and finitely many $l$'s and $r$'s.  Recall from \S\ref{subsec:compSeqPF} that this is a necessary and sufficient condition for $\compseq$ forming a compatible sequence of piecewise Fagnano orbits.  Specifically, there exists a unique integer $N\geq 1$ such that $\orbitksiang{N}{\pi/3}=\pfix{N}{\xoo}$ is a piecewise Fagnano orbit of $\omegaksi{N}$, yet for every $m<N$, $\orbitksiang{m}{\pi/3}$ is not.  Moreover, for every $n\geq N$, $\orbitksiang{n}{\pi/3}$ is a piecewise Fagnano orbit.  It follows that the first element $\xii{N}{1}$ in the footprint of the orbit $\orbitksi{N}$ has a finite address ending in either $13$ or $31$.  Likewise, for each $n\geq N$, $\xii{n}{1}$  has a finite address ending in either $13$ or $31$.  Specifically, for every $i$ such that $N\leq i\leq n$, $(\xii{n}{0})_i\neq 0,2,4$.

Each $\xii{n}{k_n}$ has an address that can be determined from $\xii{n}{0}$ (or $\xii{n}{1}$, if one prefers).  By utilizing the local and global symmetry of the snowflake (see Figure \ref{fig:kochSymmetries}), one may determine the footprint of $\orbitksiang{n}{\pi/3}$ as a non-dynamically ordered set of points.  One may then consider the inverse limit of footprints of the orbits $\pfix{n}{\xii{n}{0}}$, where one forms the inverse limit by considering as our transition maps  $\tau_{mn}$ (for $m\leq n$) the truncation of finite addresses, as defined above in this section.  Letting $\mathcal{F}_n(\xii{n}{0})$ be the footprint of the orbit $\orbitksiang{n}{\pi/3}$, we write the inverse limit of the footprints as

\begin{align}
\varprojlim \mathcal{F}_n(\xii{n}{0}) &= \left\{(\xii{i}{k_i})_{i=N}^\infty \in\prod_{i=N}^\infty  \mathcal{F}_i(\xii{i}{0}) | \tau_{mn}(\xii{n}{k_n}) = \xii{m}{k_m}\text{ for all } N\leq m\leq n \right\}. 
\end{align}

Due to the fact that this set lacks a dynamical ordering, we cannot say that it represents a section of any billiard flow.  Therefore, we next examine how to formulate the footprint of an orbit of $\omegaks$ as an inverse limit of footprints of orbits.

Consider $\zeta_n$, the number of pairs of points preceding $\xii{n}{k_n}$ in the footprint $\mathcal{F}_n(\xii{n}{0})$ of the orbit $\orbitksiang{n}{\pi/3}$. One can easily see that $2\zeta_n+\beta_n = k_n$, for some $\beta_n\in\{1,2\}$.  Then, for $(\xii{i}{k_i})\in \varprojlim \mathcal{F}_i(\xii{i}{0})$, we have that

\begin{align}
\xii{n-1}{k_{n-1}} &=\tau_{n-1,n}(\xii{n}{k_n}) = \tau_{n-1,n}(\xii{n}{2\zeta_n+\beta_n})
\end{align}

\noindent and $\xii{n-1}{\zeta_n + 1}$ is the basepoint of $f_{n-1}^{\zeta_n +1}(\xii{n-1}{0},\pi/3)$.  By the definition of a piecewise Fagnano orbit (see Definition \ref{def:pfn}), we see that $\xii{n-1}{k_{n-1}} =  \xii{n-1}{\zeta_n+1}$, and hence that $k_{n-1} = \zeta_n+1$. Consider the map $F_{n-1,n}$ given by 

\begin{align}
F_{n-1,n}(\xii{n}{k_n},\theta_n^{k_n}) &= f_{n-1}^{\zeta_n+1}\circ \tau_{n-1,n}\times \iota_{n-1,n}\circ f_n^{-k_n}(\xii{n}{k_n},\theta_n^{k_n}).
\end{align}

%We define $F_n:\prod_{i=1}^\infty \orbitksiang{i-1}{\theta_{i-1}^{0}}\to \orbitksiang{n}{\theta_n^{0}}$ by $F_n((x_i^{k_i},\theta_i^{k_i})_{i=0}^\infty) := (x_n^{k_n},\theta_{n}^{k_n})$. Then $F_{nm}: \orbitksiang{n}{\theta_n^0}\to \orbitksiang{m}{\theta_m^0}$ is defined by 

%\begin{eqnarray}
%F_{nm}(x_n^{k_n},\theta_n^{k_n}) &:=& f_m^{k_n\text{(mod } \#\orbitksi{n})}\circ \tau_{nm}\times \iota \circ f_n^{-k_n}(x_n^{k_n},\theta_n^{k_n}) \hstr[2] = \hstr[2] (x_m^{k_n\text{(mod } \#\orbitksi{n})},\theta_m^{k_n \text{(mod } \#\orbitksi{n})}).
%\end{eqnarray}

The map $F_{n-1,n}$ is well defined.  In general, if $m < n$, then the map

\begin{align}
F_{m,n}(\xii{n}{k_n},\theta_n^{k_n}) &= f_{m}^{\zeta_{m+1}+1}\circ \tau_{n-1,n}\times \iota_{n-1,n}\circ f_n^{-k_n}(\xii{n}{k_n},\theta_n^{k_n})
\end{align}

\noindent constitutes the proper transition map needed to construct the inverse limit of footprints of a compatible sequence of piecewise Fagnano orbits $\{\orbitksiang{i}{\theta_i^0}\}_{i=N}^\infty$. We denote the resulting inverse limit as follows:

\begin{align}
\mathcal{F}(\xoo) := &\varprojlim \mathcal{F}_i(\xio{i}) = \\
				 	\notag \left\{(x_i^{k_i},\theta_i^{k_i})_{i=N}^\infty \in\right.& \left.  \prod_{i=N}^\infty \orbitksiang{i}{\theta_i^0}|F_{nm}(x_n^{k_n},\theta_n^{k_n})=(x_m^{k_m},\theta_m^{k_m})\text{, for all } N\leq m \leq n\right\}.
%				 	\notag &= \left\{(x_i^{k_i},\theta_i^{k_i})_{i=1}^\infty \in \prod_{i=1}^\infty \orbitksi{i}{\theta_i^0}\vert k_m = k_n \text{mod }|\mathscr{O_m}|\right\},
\end{align}

\noindent We define $x^0\in KS$ to be $\lim_{i\to\infty} \xio{i}$, which is the limit (in the plane) of the compatible sequence of initial basepoints $\{\xio{i}\}_{i=N}^\infty$.

As will be further discussed later on, we would like to think of $\mathcal{F}(\xoo)$ as the `\textit{footprint}' (or `\textit{Poincar\'e section}') of a `\textit{piecewise Fagnano periodic orbit}' of $\omegaks$, namely, the piecewise Fagnano periodic orbit with initial basepoint $x^0:=\lim_{i\to\infty} \xio{i}$ in the direction of $\pi/3$.
\vspace{4 mm}
\subsubsection{Topological and geometric properties of the footprint $\mathcal{F}(\xoo)$}
\label{subsubsec:topologicalAndGeometricPropertiesOfFootprint}
%We have been able to determine various properties of periodic orbits, each with an initial direction of angle $\pi/3$.  In particular, we will prove in this section that every periodic orbit of the Koch snowflake with an initial direction of angle  $\pi/3$ has a finite length. In addition to this, there is a particular subcollection of the orbits with finite length for which (i) the union of all of their footprints constitutes the \textit{natural} Cantor set in the Koch snowflake $KS$ and (ii) the footprint of each individual periodic orbit is a topological Cantor set.  Such facts will rely heavily on the machinery developed in \S\S\ref{subsec:compSeqPF}--\ref{subsec:aGeneralDiscussionOfOrbits}, the addressing system discussed in \S\ref{subsec:addressSystemForKS} and the local symmetry discussed in the caption of Figure \ref{fig:kochSymmetries}.

Before discussing our main results regarding the period, length and topological nature of the orbits, we establish a 1--1 correspondence between elements of $M(\cantor)$ and the collection of footprints $\{\mathcal{F}(\xoo)\}_{\xoo\in M(\cantor)}$, as well as a similar statement for the geodesic flow (and the closed geodesic paths corresponding to the \textit{unfolded} billiard ball paths) on the associated flat surface.

%Though no closed formula may be given for the length of $\orbitks$ if $\xoo$ has an aperiodic ternary representation, such an element $\xoo$ \textit{is} a limit point of $M(\cantor)$.  The question that we may ask is whether or not the length of $\orbitks$ can be shown to be the limit of a sequence of explicit lengths of orbits of the snowflake given in terms of lengths of piecewise Fagnano orbits and $\cantor$-orbits.

%\begin{theorem}
%If $\gpfx{\yoo}$ is a generalized piecewise Fagnano orbit of the Koch snowflake, then there exists a sequence of piecewise Fagnano orbits $\pfx{\xoo}$ of the Koch snowflake billiard $\omegaks$ such that $\lim_{\xoo\to\yoo} |\pfx{\xoo}| = |\gpfx{\yoo}|$.
%\end{theorem}

%\begin{proof}
%Let $\yoo$ be an element of $\overline{M(\cantor)}\setminus T$ with a ternary representation (in terms of $l,c,r$) consisting of infinitely many $l$'s, $c$'s and $r$'s.  For every $n\geq 1$, one may form a word representing some element $\xoo$ that is $n$ characters in length such that $\xooi{i}=\yooi{i}$ for every $i\leq n$.  Then, we are allowed to stipulate that the ternary representation of $\xoo$ is such that $\xooi{j} = c$ for every $j> n$.  Continuing in this fashion, we construct a sequence of initial basepoints of orbits of the equilateral triangle billiard $\omegaksi{0}$ such that the lengths of the corresponding orbits of $\omegaks$ converge to the length of the generalized piecewise Fagnano orbit $\gpfx{\yoo}$.
%\end{proof}

\begin{theorem}
\label{thm:orbitsAreCantorSets}
If $\{\orbitksiang{i}{\theta^0_i}\}_{i=0}^\infty$ is a compatible sequence of piecewise Fagnano orbits, then the `footprint' $\mathcal{F}(\xoo)$, given by the inverse limit $\varprojlim \orbitksi{i}$, is a topological Cantor set \emph{(}i.e., a totally disconnected and perfect compact set\emph{)}. Moreover, the above inverse limit of the footprints is a self-similar Cantor set.
\end{theorem}

\begin{proof}
Let $\{\orbitksiang{i}{\theta^0_i}\}_{i=0}^\infty$ be a compatible sequence of piecewise Fagnano orbits. By Theorem \ref{thm:invLimitTotDisc}, $\mathcal{F}(\xoo)$ is totally disconnected and compact, being the inverse limit of finite sets.  What remains to be shown is that $\mathcal{F}(\xoo)$ is a perfect set.

Consider an element $(x_i^{k_i},\theta_i^{k_i})_{i=0}^\infty \in \mathcal{F}(\xoo)$.  Fix $n\geq 0$.  Then there exists $(y_i^{j_i},\phi_i^{j_i})_{i=0}^\infty \in \mathcal{F}(\xoo)$ such that 1) $(x_i^{k_i},\theta_i^{k_i})=(y_i^{j_i},\phi_i^{j_i})$ for all $i\leq n$ and 2) $x_{n+1}^{k_{n+1}} \neq y_{n+1}^{j_{n+1}}$. Concretely, this element $(y_i^{j_{i}},\phi_i^{j_{i}})_{i=0}^\infty$ is determined from $(x_i^{k_i},\theta_i^{k_i})_{i=0}^\infty$ by way of the local symmetry.    As such, we can continue to construct a sequence of elements in $\mathcal{F}(\xoo)$ that converges (with respect to either the Euclidean metric or a metric defined on the space of addresses) to $(x_i^{k_i},\theta_i^{k_i})_{i=0}^\infty \in \mathcal{F}(\xoo)$.

To see that such a topological Cantor set is a self-similar Cantor set, one simply recognizes the fact that for every footprint $\mathcal{F}(\xoo)$, there is a finite collection of IFS's $\{\Phi_{\xoo,i}\}_{i=1}^\Xi$, each IFS $\Phi_{\xoo,i}$ consisting of two contraction mappings, and each giving rise to a unique fixed point attractor that is a self-similar set in its own right.  The union of these fixed point attractors is then the self-similar footprint and one can easily reproduce a graphical representation of the footprint.  We defer further explanation of this to \S\ref{subsec:plausibilityArgument}.\footnote{More specifically, for more detailed information concerning this IFS, we refer to the discussion immediately following the proof of Theorem \ref{thm:limitOfLengths} in \S\ref{subsec:plausibilityArgument} below.}
\end{proof}

\begin{proposition} 
If $\{\mathcal{F}(\xoo)\}_{\xoo\in M(\cantor)}$ is the collection of all footprints of piecewise Fagnano orbits of the Koch snowflake billiard $\omegaks$ and $\xoo,\yoo\in M(\cantor)$, then  $\xoo\neq\yoo$ if and only if  $\mathcal{F}(\xoo)$ and $\mathcal{F}(\yoo)$ have no elements in common.  
\end{proposition}

\begin{proof}
Suppose $\mathcal{F}(\xoo) \neq \mathcal{F}(\yoo)$.  Recall that an address of a point can be straightened so as to identify a point $x^0\in KS$ that is compatible with an element $\xoo$ (in the direction of $\pi/3$ in the interior of $KS$) that is also an element of the unit interval base of $KS_0$, which we have been denoting by $I$ throughout the paper.  This element may then be taken to be $\lim_{i\to\infty} \xio{i}$, with $\{\xio{i}\}_{i=0}^\infty$ being a compatible sequence of initial basepoints (again, compatible in the direction of $\pi/3$).  That is, $x^0 = \lim_{i\to\infty} \xio{i}$ and $y^0 = \lim_{i\to\infty} \yio{i}$.  Suppose $x^0\neq y^0$.  Then there exists  $\nu>0$ such that $||x^0-y^0||>\nu$ (where we have taken $||\cdot||$ to be the Euclidean norm in $\mathbb{R}^2$).  Therefore, there exists $n\geq 0$ such that
\begin{align}
||\xio{n} - \yio{n}|| &> \nu
\end{align}                                                                                                                                                                                                                                                                                                                                                                                                             
\noindent or
\begin{align}
\xio{n} &\neq  \yio{n} 
\end{align}                                                                                                                                                                                                                                                                                                                                                                                                             
Since $\xio{n}$ and $\xoo$ are compatible in the direction of $\pi/3$ and $\yio{n}$ and $\yoo$ are compatible in the direction of $\pi/3$, it follows that $\xoo\neq \yoo$. The converse holds since $\xoo\neq \yoo$ implies $x^0\neq y^0$.  
  %Consider the flow $T_\theta$, where $\theta=\pi/3$ on the flat surface $\mathcal{S}(KS_n)$.  Fixing $\xoo$ and $\yoo$, two distinct points in $M(\cantor)$, we see that $\xoo$ and $\yoo$ cannot have the same ternary expansions nor have ternary expansions representing the same point in the unit interval.  This follows from the fact that the points without unique representations were the points in $T$. The points $\xoo$ and $\yoo$ determine distinct compatible sequences of piecewise Fagnano orbits.  More precisely, there exist $N_1$ and $N_2$ (without loss of generality, $N_1\leq N_2$) such that $\pfix{N_1}{\xoo}$ and $\pfix{N_2}{\yoo}$ are piecewise Fagnano orbits and $\pfix{n_1}{\xoo}\neq\pfix{n_2}{\yoo}$, for every $n_1$, $n_2\geq N_1$.  Therefore, $\compseqang{\pi/3} \neq \compseqangx{\pi/3}{y}$.
 
%If the flow lines $T^t_{(\xoo,\theta)}$ and $T^t_{(\yoo,\theta)}$ have any points in common, it must be the case that they have all points in common. In other words, $x_1=x_2$, which contradicts our assumption.
\end{proof}

\subsection{Stabilizing orbits (or $\cantor$-orbits) of $\mathbf{\omegaks}$}
\label{subsec:stabilizingOrbitsKS}
\begin{figure}
\begin{center}
\includegraphics{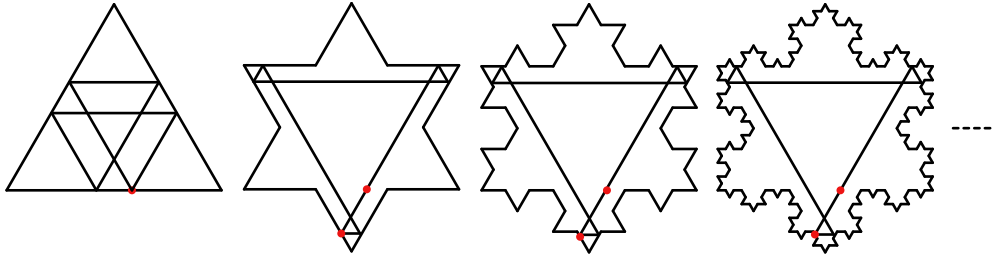}
\end{center}
\caption{A stabilizing periodic orbit of the Koch snowflake billiard $\omegaks$.  The initial basepoint $x^0$ of the stabilizing orbit is $\xio{1}$, which has the address $\alpha(\xio) = 51\overline{24}$.  The initial basepoint $\xio{1}$ is compatible with the element $\xoo = 7/2\in I$, which has a ternary representation of $c\overline{rl}$.}
\label{fig:cantorOrbit7-12}
\end{figure}

In \S\ref{subsec:eventConstComSeq}, we described what we call an \textit{eventually constant compatible sequence} of periodic orbits.  Such a compatible sequence of orbits was determined by the initial basepoint $\xoo$ of the first orbit $\orbitksiang{0}{\pi/3}$ in the compatible sequence $\compseqang{\pi/3}$.  The ternary representation of $\xoo$ is comprised of finitely many $c$'s and infinitely many $l$'s and $r$'s. Let $N$ be a positive integer such that for every $n\geq N$ we have that $\orbitksiang{n}{\pi/3}=\orbitksiang{N}{\pi/3}$.\footnote{Throughout \S4, we chose to write the ternary expansion of an element of $M(\cantor)$ in terms of the characters $\{0,1,2\}$.  Because of the obvious confusion with the alphabet $\{0,1,2,3,4,5\}$ that we use when addressing points of the Koch snowflake $KS$, we use instead the characters $l$ $c$ and $r$ to express the ternary representation of elements of $M(C)$ and, in general, of the unit interval $I$, viewed as the base of $KS_0$.} To be clear, the footprint of the orbit $\orbitksiang{N}{\pi/3}$ consists of finitely many points with addresses in $KS$ that contain finitely many $1$'s and $3$'s.  We then say that $\orbitksiang{N}{\pi/3}$ is a \textit{stabilizing} periodic orbit (or a $\cantor$-orbit) of $\omegaks$; see Figure \ref{fig:cantorOrbit7-12}.

An example of a $\cantor$-orbit of $\omegaks$ is one for which $\xoo = 7/12$.  This element is not in $M(\cantor)$ and the corresponding eventually constant compatible sequence of periodic orbits $\compseqang{\pi/3}$ is constant after the second element of the compatible sequence of periodic orbits since the next element of the compatible sequence of initial basepoints corresponds to an element of the Cantor set on the side with address $51$.    

In \S\ref{subsec:aGeneralDiscussionOfOrbits}, Theorems  \ref{thm:periodOfOrbitksi} and \ref{thm:lengthBilliardBallPath} highlighted a relationship between the ternary representation of the initial basepoint $\xoo$ of $\orbitksiang{0}{\pi/3}$ and the length and period of an orbit $\orbitksiang{n}{\pi/3}$.  The ternary representation of the element $\xoo$ also provides further information about the behavior of a compatible sequence of periodic orbits when such a compatible sequence is eventually constant. 

\begin{theorem}
Any element $\xoo\in I$ with a ternary representation consisting of finitely many $c$'s and infinitely many $l$'s and $r$'s constitutes the first element in a compatible sequence of initial basepoints giving rise to a compatible sequence of periodic orbits that is eventually constant.
\end{theorem}

\begin{example}
The element $\xoo = c\overline{rl}$ is such an element.  The eventually constant compatible sequence of periodic orbits stabilizes after the zeroth level approximation $\orbitksiang{0}{\pi/3}$, as depicted in Figure \ref{fig:cantorOrbit7-12}.
\end{example}

Further, one can determine the approximation at which the compatible sequence of periodic orbits \textit{stabilizes}.

\begin{theorem}
Let $\xoo\in I$ have a ternary representation consisting of finitely many $c$'s and an infinite number of $l$'s and $r$'s.  Furthermore, let $N$ be the number of left shifts required to eliminate all characters $c$ from the ternary representation of the initial basepoint $\xoo$ of the orbit $\orbitksiang{0}{\pi/3}$. Then, for every $n\geq N$, we have that $\orbitksiang{n}{\pi/3} = \orbitksiang{N}{\pi/3}$.
\end{theorem}

Interestingly enough, if $\xio{N}$ is the \textit{initial compatible basepoint} of a stabilizing periodic orbit, then we observe that the ternary representation $\xio{N}$ relative to the side $s_{N,k}$ containing $\xio{N}$ can be determined from the ternary representation of $\xoo$.

\begin{proposition}
Let $\xio{N}$ be the initial basepoint of a stabilizing orbit of $\omegaks$. Then the ternary representation of $\xio{N}$ \emph{(}relative to a segment of length $1/3^N$\emph{)} is given by the $N$th left-shift of the ternary representation of the point $\xoo$.  
\end{proposition}

\begin{example}
If $x_0^0 = 5/12$, then the ternary expansion of this number is $c\overline{lr}$.  The ternary expansion of the initial basepoint $\xio{1}$  of the stabilizing orbit $\orbitksiang{1}{\pi/3}$ is given by $l\overline{lr}$.  Relative to the side of length $1/3$ on which this point can be found, the ternary expansion is $\overline{lr}$.
\end{example}

%One may then say that a stabilizing orbit of $\Omega(KS_n)$ is an orbit of $\Omega(KS)$, because for each approximation $\Omega(KS_n)$ and each point with a ternary representation containing finitely many $c$'s and infinitely many $l$'s and $r$'s on some side of $\Omega(KS_n)$, there exists an open connected neighborhood about $x_n^0$ of positive volume.  

%Hence, reflection can be defined in each approximation and an orbit can be determined for every $n$ that is \textit{the same} orbit for every $n$.  

%we demonstrated the existence of a particular type of compatible sequence of periodic orbits: one for which $\compseq$ was eventually constant.  The fact that the compatible sequence of periodic orbits was eventually constant indicates that there exists $N\geq 0$ such that an initial basepoint $x_N^0$ of an orbit $\orbitksiang{N}{\theta_N^0}$ residing on a side $s_{N,k}$ corresponds to an element of $\cantor_{N,k}'\subseteq s_{N,k}$.  More importantly, $\{x_i^0\}_{i=0}^\infty$ is a compatible sequence of basepoints converging to a point on the boundary $KS$ with an address consisting of finitely many 1's and 3's and infinitely many 0's and 2's or 2's and 4's, corresponding to a shift to the left third or right third along the side $s_{N,k}$, depending on whether or not the finite address of $s_{N,k}$ ended in a 1 or 3, respectively.

%Therefore, we call $\orbitks$ a $\cantor$-orbit of $\omegaks$ when $\theta^0_i = \pi/3$, and $\xoo\in I$ has a ternary representation containing at most finitely many characters $c$.

\subsection{Piecewise Fagnano `orbits' of $\mathbf{\omegaks}$}
\label{subsec:plausibilityArgument}
In \S\ref{subsec:aGeneralDiscussionOfOrbits}, we saw that for every $n\geq 0$, there was an explicit formula for the period and length of a piecewise Fagnano orbit $\pfix{n}{\xoo}$ of $\omegaksi{n}$.  Specifically, the length and period were determined by the ternary expansion of $\xoo$.  Moreover, in \S\ref{subsec:footprintsOfPeriodicOrbitsOfTheKochSnowflakeBilliard}, we constructed the inverse limit of the Poincar\'e sections of the orbits comprising a compatible sequence of piecewise Fagnano orbits $\compseqang{\pi/3}$.  Based on the explicit formulas for the period and length (see Theorems \ref{thm:periodOfOrbitksi} and \ref{thm:lengthBilliardBallPath}), one sees that 

\begin{align}
\lim_{n\to\infty} \#\orbitksiang{n}{\pi/3} &= \infty 
\end{align}

\noindent and, as we next explain, that

\begin{align}
\lim_{n\to\infty} &|\orbitksiang{n}{\pi/3}| 
\end{align}

\noindent is finite.

\begin{theorem}
\label{thm:limitOfLengths}
Let $\compseqang{\pi/3}$ be a compatible sequence of piecewise Fagnano orbits.  Then

\begin{align}
\lim_{n\to\infty} &|\orbitksiang{n}{\pi/3}|
\end{align}

\noindent exists and is finite.
\end{theorem}

\begin{proof}
Let $\mathscr{L}$ be the length of the Fagnano orbit of the equilateral triangle billiard $\omegaksi{0}$. Let $\compseqang{\pi/3}$ be a compatible sequence of piecewise Fagnano orbits.  Then, for each $n$, Theorem \ref{thm:lengthBilliardBallPath} implies (with the notation introduced in Equation (\ref{eqn:charFunctionOnOneC})) that 

\begin{equation}
|\orbitksiang{n}{\pi/3}| = 2\mathscr{L} + \sum_{i=2}^n \chi\left[\xooi{i}\right]\#\orbitksiang{i-1}{\pi/3}\frac{\mathscr{L}}{3^i}.
\end{equation}

Since, for each $n\geq 1$, $|\orbitksiang{n}{\pi/3}|$ is bounded by $|\ppfi{n}|$ and $\lim_{n\to\infty} |\ppfi{n}| = 4\mathscr{L}$, it follows that 
\begin{align}
\notag \lim_{n\to\infty} |\orbitksiang{n}{\pi/3}| &= 2\mathscr{L} + \lim_{n\to\infty}\sum_{i=2}^n \chi\left[\xooi{i}\right]\#\orbitksiang{i-1}{\pi/3}\frac{\mathscr{L}}{3^i}\\
\label{eqn:finitelimt}
																 &= 2\mathscr{L} + \sum_{i=2}^{\infty} \chi\left[\xooi{i}\right]\#\orbitksiang{i-1}{\pi/3}\frac{\mathscr{L}}{3^i}\\
\notag 														 &\leq 2\mathscr{L} + \sum_{i=2}^{\infty} 2^i\frac{\mathscr{L}}{3^i}\hstr[2] = \hstr[2] \frac{10}{3}\mathscr{L}< \infty.
\end{align}

We note that since $\{|\orbitksiang{i}{\pi/3}|\}_{i=0}^\infty$ is a monotonically increasing sequence (which follows from the fact that the sequence of partial sums representing this sequence is monotonically increasing), we deduce that the limit exists and (by Equation (\ref{eqn:finitelimt})) is finite.  
%If $\xoo \in M(\cantor)$ has an (eventually) periodic ternary representation, then one may be able to calculate a specific value for the length of the orbit $\orbitks$, since such a calculation will essentially reduce to a finite sum of infinite geometric series plus a finite number of terms not contributing to any of the geometric series.  If $\xoo$ is an element with an aperiodic ternary representation, then no convenient closed formula may be given, but the length will be finite, as demonstrated above.
\end{proof}
Such facts support the following plausibility argument, according to which some suitable limiting set can be considered as an orbit of the Koch snowflake billiard.

Let $\compseqang{\pi/3}$ be a compatible sequence of piecewise Fagnano orbits.  Then there exists a least integer $n\geq 0$ such that $\orbitksixang{n}{x}{\pi/3}$ is a piecewise Fagnano orbit and $\orbitksixang{n-1}{x}{\pi/3}$ is not.  Denote by $\Xi$ the period of $\orbitksixang{n-1}{x}{\pi/3}$.

Recall from the proof of Theorem \ref{thm:orbitsAreCantorSets} that we claimed that the footprint $\mathcal{F}(\xoo)$ of a piecewise Fagnano orbit was a self-similar subset of $KS$.  Concretely, $\mathcal{F}(\xoo)$ is the finite union of (abutting) self-similar sets.  Specifically, there are $\Xi$ many IFS's involved in the union describing $\mathcal{F}(\xoo)$.  That is, if for each $i\leq\Xi$, we have that $\Phi_{\xoo,i}$ is the operator associated with an IFS (which we denote by $\{\Phi_{\xoo,i,j}\}_{j=1}^2$) and $\mathscr{F}(\xoo)_i$ is the unique fixed point attractor of $\Phi_{\xoo,i}$, then 

\begin{align}
\Phi_{\xoo,i}(\mathcal{F}(\xoo)_i) := \bigcup_{j=1}^2 &\Phi_{\xoo,i,j}(\mathcal{F}(\xoo)_i) = \mathcal{F}(\xoo)_i.
\end{align}

\noindent As a result, we have that

\begin{align}
\notag \mathcal{F}(\xoo) &= \bigcup_{i=1}^\Xi \Phi_{\xoo,i}(\mathcal{F}(\xoo)_i) = \bigcup_{i=1}^\Xi \mathcal{F}(\xoo)_i.
\end{align}

%Each self-similar set $\mathcal{F}(\xoo)_i$ involved in the union is the unique fixed point attractor of a particular IFS, which we denote by .  Then,

For each $i\leq \Xi$, there exist a suitably defined non-expansive map $\rho_i$ that, together with $\Phi_{\xoo,i,1}$ and $\Phi_{\xoo,i,2}$, constitutes a new IFS $\{\rho_i,\Phi_{\xoo,i,1},\Phi_{\xoo,i,2}\}$ (in a non-standard sense).  Such an IFS has an associated operator that can be used to produce the piecewise Fagnano orbit $\pfix{n}{\xoo}$ and every subsequent piecewise Fagnano orbit in the compatible sequence of piecewise Fagnano orbits.  We denote the associated operator by 

\begin{align}
\widetilde{\Phi}_{\xoo,i}:= \rho_{i}\cup\Phi_{\xoo,i,1}\cup\Phi_{\xoo,i,2}.  
\end{align}

Before showing exactly how such an operator may be utilized to produce a piecewise Fagnano orbit, we introduce a key notion.  Let $\mathscr{F}$ be the Fagnano orbit of $\omegaksi{0}$.  Then, for each $i\leq \Xi$, there is a contraction map $\Psi_i$ acting on $\mathscr{F}$ in such a way that 

\begin{align}
\label{eqn:pfnAsUnionOfleaves}
\pfix{n}{\xoo} &= \bigcup_{i=1}^\Xi \Psi_{i}(\mathscr{F})\cup\orbitksixang{n-1}{x}{\pi/3}.
\end{align}

\noindent That is, $\pfix{n}{n}$ is given by Equation (\ref{eqn:pfnAsUnionOfleaves}) once a proper dynamical ordering is reestablished for $\bigcup_{i=1}^\Xi \Psi_{i}(\mathscr{F})\cup\orbitksixang{n-1}{x}{\pi/3}.$ Let $\mathscr{F}_k^i := \Psi_i(\mathscr{F})$, where $k$ is an integer meant to indicate that $\Psi_i(\mathscr{F})$ is a scaled copy of $\mathscr{F}$, in the scaling ratio of $1/3^k$.  

Finally, we have that

\begin{align}
\notag \pfix{n}{\xoo} &= \bigcup_{i=1}^\Xi \mathscr{F}_k^i\cup\orbitksixang{n-1}{x}{\pi/3}
\end{align}

\noindent and, in general, for $p\geq 0$,

\begin{align}
\notag \pfix{n+p}{\xoo} &= \bigcup_{i=1}^\Xi \widetilde{\Phi}^p_{\xoo,i}(\mathscr{F}_k^i)\cup\orbitksixang{n-1}{x}{\pi/3}.
\end{align}

\begin{figure}
\begin{center}
\begin{overpic}[scale=.5,unit=1mm]{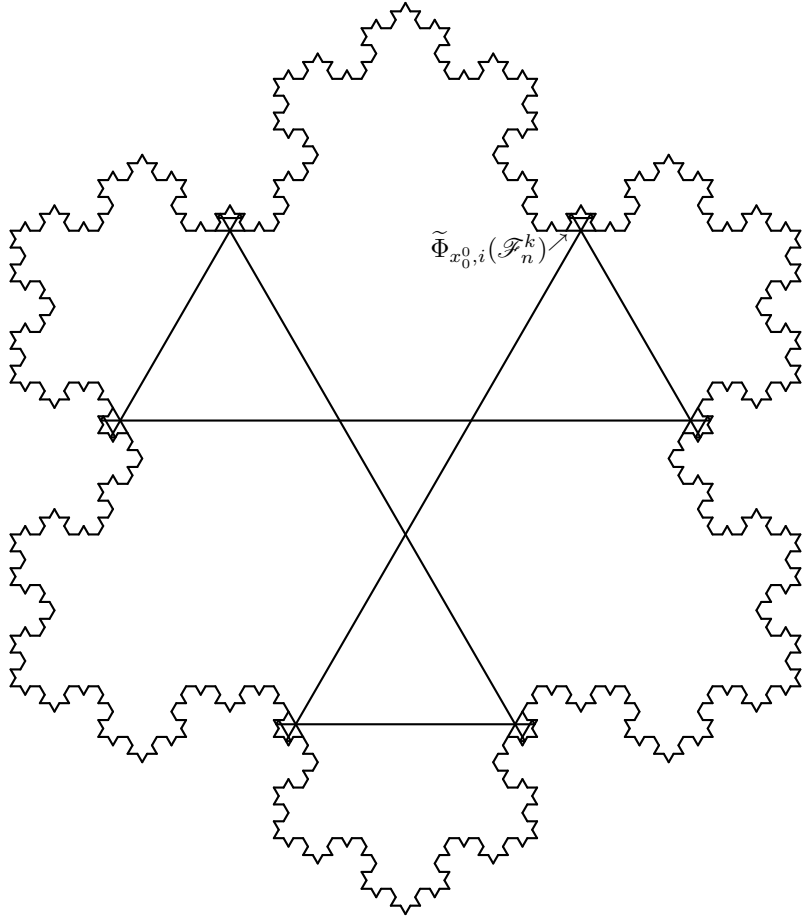}
    \put(46,72){$\widetilde{\Phi}_{\xoo,i}(\mathscr{F}_n^k)^{\huge{\nearrow}}$}
    %\put(34,5){\includegraphics[scale=.5]{plausibility.eps}}
\end{overpic}  
\end{center}
\caption{Depicted in this figure is a piecewise Fagnano orbit $\pfix{4}{\xoo}$ of $\omegaksi{4}$. In addition, this orbit is an orbit in a compatible sequence of piecewise Fagnano orbits.  We have drawn in a subset of the ghosts of $\omegaksi{4}$ so as to see how it is one orbit is compatible with another. The least $n$ such that $\orbitksixang{n}{x}{\pi/3}$ is a piecewise Fagnano orbit is $n=3$.  The compatible orbit $\orbitksixang{2}{x}{\pi/3}$ has period $6$.  The iterated function system producing $\pfix{4}{\xoo}$ from $\pfix{3}{\xoo}$ is shown.  In Figure \ref{fig:plausibilityZoomed}, we enlarge that particular piece of the billiard table $\omegaksi{4}$ and examine the effect of one of the IFS's $\widetilde{\Phi}_{\xoo,i}$.}
\label{fig:plausibility}
\end{figure}

\noindent Again, a dynamical ordering must be reestablished for this equality to be true.  Given the fact that $\pfix{n+p}{\xoo}$ is a finite orbit with finite period, imposing such an ordering is not a concern.  See Figure \ref{fig:plausibility} and the caption therein for a detailed example illustrating the above discussion.

\begin{figure}
\begin{center}
\begin{overpic}[scale=4,unit=1mm]{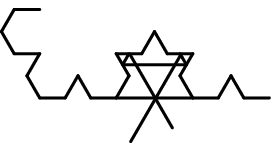}
    \put(37.5,12){$\widetilde{\Phi}_{\xoo,i}(\mathscr{F}_k^i)^{\huge{\nearrow}}$}
	\put(55.5,25){$\mathscr{F}_k^i$}
	\put(29,36){$\Phi_{\xoo,i,1}(\mathscr{F}_k^i)_{\huge{\searrow}}$}
	\put(67.25,36){${}_\swarrow\Phi_{\xoo,i,2}(\mathscr{F}_k^i)$}
    %\put(34,5){\includegraphics[scale=.5]{plausibility.eps}}
\end{overpic} 
\end{center}
\caption{The effect of the IFS $\widetilde{\Phi}_{\xoo,i}$, $i\leq 6 = \#\orbitksi{2}$.}
\label{fig:plausibilityZoomed}
\end{figure} 
%  Positioning the \LaTeX\ output with \verb#\put# commands:  
%  (default for \texttt{unit} is \verb#\unitlength#)  

%\begin{figure}[htbp]
	%\centering
	%	\includegraphics{plausibility.eps}
	%\caption{test}
	%\label{fig:plausibility}
%\end{figure}

%Recall from the proof of  that we claimed that the footprint $\mathcal{F}(\xoo)$ of a piecewise Fagnano orbit was a self-similar subset of $KS$.  We denote the IFS for which $\mathcal{F}(\xoo)$ is the unique fixed point attractor by $\Phi_{\xoo}$.  Then, including in the collection of contraction mappings a suitably defined nonexpansive map, we recover an IFS $\widetilde{\Phi}_{\xoo}$.  If $\mathscr{F}_k$ is the Fagnano orbit of $\omegaksi{0}$ at scale $k$, then such a map $\widetilde{\Phi}_{\xoo}$ can be used to recover the graphical representation of the orbit $\orbitksiang{n}{\pi/3}$ by evaluating $\widetilde{\Phi}_{\xoo}^n$ at $\mathscr{F}_k$. 

\vspace{3 mm}
\subsubsection*{Plausibility argument}
\label{pls:piecewiseFagnanoOrbitsOfKS}
%Let $\compseqang{\pi/3}$ be a compatible sequence of piecewise Fagnano orbits.  The footprint associated with the compatible sequence is then denoted by $\mathcal{F}(\xoo)$. 

%There exists an IFS consisting of contraction mappings that is associated with the footprint $\mathcal{F}(\xoo)$.  Specifically, $\mathcal{F}(\xoo)$ is the unique fixed point attractor of the operator $\Phi_{\xoo}$ determined by the IFS.   As mentioned above, if $\mathscr{F}$ is the Fagnano orbit of $\omegaksi{0}$, then there exists a suitably defined non-expansive mapping that,  when combined with the contraction mappings comprising the IFS, determines another iterated function system with an associated operator  $\widetilde{\Phi}_{\xoo}$  such that, for some $k\geq 0$, $\pfix{n}{\xoo} = \widetilde{\Phi}(\mathscr{F}_k)\cup \orbitksixang{n-1}{x}{\pi/3}$.

Let $\compseqang{\pi/3}$ be a compatible sequence of piecewise Fagnano orbits and 

\begin{align}
\pfx{\xoo}:= \lim_{p\to\infty} \bigcup_{i=1}^\Xi \widetilde{\Phi}_{\xoo,i}^{n+p}(\mathscr{F}_k^i) \cup \orbitksixang{n-1}{x}{\pi/3}\label{eqn:pfxoo},
\end{align} 

\noindent  where $n$ is the least integer such that $\orbitksixang{n}{x}{\pi/3}$ is a piecewise Fagnano orbit of $\omegaksi{n}$ in the compatible sequence of piecewise Fagnano orbits $\compseqang{\pi/3}$.  Since $\mathcal{F}(\xoo)$ is an uncountable set and a Poincar\'e section is at most countably infinite, it is plausible that there exists a suitably defined map $f$ on a countably dense subset of $\mathcal{F}(\xoo)$ that constitutes the Poincar\'e section of an orbit of $\omegaks$ whose graphical representation is given by $\pfx{\xoo}$.  This may then allow us to reconstitute $\pfx{\xoo}$ as a true orbit of the Koch snowflake billiard $\omegaks$.

%\begin{align}
%\pf^{\xoo}:=\lim_{N\to\infty} \Phi_{\xoo}^N (\Delta),
%\end{align}

%\noindent then there exists a suitably defined map $f$ defined on a dense subset of $\mathcal{F}(\xoo)$ that constitutes the Poincar\'e section of an orbit whose graphical representation is given by $\pf^{\xoo}$.  Furthermore, this may then allow us to reconstitute $\pf^{\xoo}$ as an elusive orbit of the Koch snowflake fractal billiard $\omegaks$. 

In light of Theorem \ref{thm:orbitsAreCantorSets} and Theorem \ref{thm:limitOfLengths}, one would then expect this periodic orbit to be `self-similar' and to have finite length.

\vspace{3 mm}
\subsubsection*{Support for the plausibility argument}
%Consider a compatible sequence of piecewise Fagnano orbits $\compseqang{\pi/3}$.  There exists a least positive integer $N$ such that $\orbitksiang{N}{\pi/3}$ is a piecewise Fagnano orbit of $\omegaksi{N}$.  

%There exists an iterated function system $\Phi_{\xoo}$ which can be used in the construction of the orbit $\orbitksiang{n}{\pi/3}$, for every $n>N$.  Concretely, if $\mathscr{F}_k$ is the Fagnano orbit at scale $k$, then there exists an iterated function system $\widetilde{\Phi}_{\xoo}$ such that

%\begin{align}
%\pfix{n}{\xoo} &= \widetilde{\Phi}_{\xoo}^{N+n}(\mathscr{F}_k)\cup\orbitksixang{n-1}{x}{\pi/3} 
%\label{eqn:pfixAsIFS}
%\end{align}

%\noindent where $\widetilde{\Phi}_{\xoo}$ is comprised of contraction mappings and a non-expansive mapping.  Note that $\Phi_{\xoo}$ may not necessarily be comprised entirely of contraction mappings.

The Poincar\'e section of $\pfix{n+p}{\xoo}$ is finite and one can reestablish a dynamical ordering on $\bigcup_{i=1}^\Xi \widetilde{\Phi}_{\xoo,i}^{n+p}(\mathscr{F}_k^i) \cup \orbitksixang{n-1}{x}{\pi/3}$. However, the limiting set $\pfx{\xoo}$ defined by Equation (\ref{eqn:pfxoo}) does not naturally come equipped with such an ordering.  Without additional knowledge, one cannot determine any dynamical ordering on $\pfx{\xoo}$.  

Consider again the footprint $\mathcal{F}(\xoo)$ of the piecewise Fagnano orbit with initial basepoint $x^0 = \lim_{i\to\infty} \xio{i}$ (and in the initial direction of $\pi/3$).  Given how the transition maps used in the construction of $\mathcal{F}(\xoo)$ are defined in terms of the billiard maps $f_n$ defined on the prefractal billiard approximations $\omegaksi{n}$, it is reasonable to expect that one can recover a suitable map that would force a dynamical ordering on $\mathcal{F}(\xoo)$.  This is, unfortunately, problematic, because in Theorem \ref{thm:orbitsAreCantorSets} we established that $\mathcal{F}(\xoo)$ was a topological Cantor set, i.e., uncountable.  Any attempt to recover a discrete set reminiscient of a Poincar\'e section will necessarily yield an at most countably infinite set.  Consequently, we expect to eventually be able to construct a suitable map defined on a countable dense subset of $\mathcal{F}(\xoo)$ that then determines an appropriate analog of the Poincar\'e section of a billiard flow for the Koch snowflake billiard $\omegaks$.  With a natural dynamical ordering imposed on a countable dense subset of the footprint $\mathcal{F}(\xoo)$, one may then be able to recover a continuous path that can be reasonably called the billiard orbit of $\omegaks$ with an initial condition $(x^0,\pi/3)$.

As mentioned, this countably dense subset would be determined by a suitably chosen map.  The orbit itself would then be determined by how this map acts, that is, on what point of $\mathcal{F}(\xoo)$ such a map initially acts.    As such, one may dictate that the initial point of departure on $KS$ must be a point compatible with some element $\xoo\in M(\cantor)$ in the direction of $\pi/3$ that then earns the \textit{name} of a \textit{piecewise Fagnano orbit} of $\omegaks$. Thus, perhaps, the closure of the continuous orbit determined from a countable dense subset may correspond to the actual footprint of a piecewise Fagnano orbit.  Moreover, if one can show that the closure of this set is the closure of the countable union of piecewise Fagnano orbits $\bigcup_{p=0}^\infty \pfix{n+p}{\xoo}$, then one has obtained not only a dynamical interpretation of the orbit, but also a natural geometric interpretation in terms of a collection of interated function systems $\{\widetilde{\Phi}_{\xoo,i}\}_{i=1}^\Xi$.

%%or consider some way of collecting all of the possible countably dense subsets into a family of orbits 

%We then define a piecewise Fagnano orbit of the Koch snowflake billiard $\Omega(KS)$ to be one for which $x^0$ has an address ending in an infinite number of alternating 1's and 3's such that $x^0$ is collinear with an element $\xoo \in M(\cantor)$ in the direction $\pi/3$. Consequently, $\xoo$ has a ternary representation terminating solely in infinitely many $c$'s. Thus, a piecewise Fagnano orbit of the Koch snowflake billiard is the inverse limit of a particular compatible sequence of piecewise Fagnano orbits.

%\begin{theorem}
%Every compatible sequence of piecewise Fagnano orbits is an inverse limit sequence of periodic orbits, and the inverse limit of such a sequence is a periodic orbit of the Koch snowflake billiard $\omegaks$. 
%\end{theorem} 

\subsection{Generalized piecewise Fagnano `orbits' of $\mathbf{\omegaks}$}
%If $\{O_i(x_i,\pi/3)\}_{i=0}^\infty$ is a compatible sequence of periodic orbits with $\xoo \in M(\cantor)$, then we may form an inverse limit sequence from those orbits.

Up until now, there are two types of periodic orbits of the snowflake.  Those for which the initial basepoint of the first element of the compatible sequence is an element of $M(\cantor)$, and those for which the initial basepoint is eventually collinear with some element of a ternary Cantor set on some side $s_{n,k}$ of some prefractal billiard $\omegaksi{n}$ not corresponding to a point with a finite ternary expansion.  

So a question naturally arises: what `periodic orbits' (in the direction of $\pi/3$) have not yet been described?  Let $\xoo \in I$ be an element with a ternary representation comprised of an infinite number of $c$'s, $l$'s and $r$'s (specifically, $c$'s and $r$'s, $c$'s and $l$'s or $c$'s, $l$'s and $r$'s).  For any $n\geq 1$, the compatible orbit $\orbitksiang{n}{\pi/3}$ is not a piecewise Fagnano orbit.  Nor is the compatible orbit ever an element of an eventually constant compatible sequence of periodic orbits.  However, the existence of a well-defined compatible sequence of periodic orbits is not a sufficient condition for such a compatible sequence being recasted as an inverse limit sequence of orbits.  Concretely, within the current framework, one cannot (at least, without further modifications) recast such a compatible sequence as an inverse limit sequence, let alone provide a plausibility argument for the existence of such an orbit, because of the nature of the initial basepoint $\xoo\in I$ of the orbit $\orbitksiang{0}{\pi/3}$.  But this does not preclude us from coming to terms with what exactly the \textit{limit} of a compatible sequence of such orbits would look like.  

\begin{conjecture}[Footprint of a generalized piecewise Fagnano orbit]\label{conj:footprintOfGeneralizedPiecewiseFagnanoOrbit}
Let $\xoo\in I$.   If $\compseqang{\pi/3}$ is a compatible sequence of generalized piecewise Fagnano orbits, then there exists a suitable family of transitions maps that one can use in order to properly formulate an inverse limit from the Poincar\'e sections of the generalized piecewise Fagnano orbits.
\end{conjecture}

\begin{conjecture}[A generalized piecewise Fagnano orbit]
If one can properly determine the footprint of a generalized piecewise Fagnano orbit, as stated in Conjecture \ref{conj:footprintOfGeneralizedPiecewiseFagnanoOrbit}, then, in analogy with the support for the plausibility argument given in \S\ref{subsec:plausibilityArgument}, one should be able to recover an associated periodic orbit of the Koch snowflake fractal billiard $\omegaks$.
\end{conjecture}

%The type of orbit that we have yet to discuss is the inverse limit of a particular compatible sequence of periodic orbits for which $\xoo$ has a particular ternary representation and the point $x^0\in KS$ has a unique address given by the limit the sequence of compatible points $(\xio{i})_{i=0}^\infty$.  If $\xoo$ is an element of $\overline{M(\cantor)}\setminus T$ such that the ternary representation consists of an infinite number of $l$'s and $r$'s \textit{and} an infinite number of $c$'s, then $\xoo$ will neither be an element of $M(\cantor)$ nor be an element of the Cantor set $\cantor\subseteq I$ (nor a scaled copy of $\cantor$ viewed as a subset of a ghost $g_{n,0}$ of the base of the equilateral triangle $KS_0$).

%Let $\xoo$ be such a value.  Then the basepoints of any periodic orbit never collide with any corners nor with any points that are midpoints of ghosts of some sides of some future approximation.  The resulting inverse limit is then a periodic orbit of the Koch snowflake.

%\begin{example}
%something to come soon
%\end{example}

\section{Open Problems and Conjectures}
\label{sec:openProblemsConj}
We close this paper by stating various additional questions, open problems and conjectures pointing to future research in this new field of `fractal billiards'.  We also offer some evidence (either theoretical or experimental) towards several of these conjectures.  Throughout much of this discussion, we will use the present important and prototypical example of the Koch snowflake billiard, but naturally, similar questions could be asked about more general fractal billiards.

\subsection{$\mathbf{\Omega(KS)}$ eventually realized as a billiard (i.e., there may be a well-defined phase space)}

As mentioned in the introduction, the Koch snowflake $KS$ is a closed, non-rectifiable and everywhere nondifferentiable curve.  This means that there is a priori no well-defined phase space.  At the heart of the description of the billiard flow is the fact that we can always obtain an inward pointing vector at the point of collision in the boundary $KS_n$ (or any other rational billiard, for that matter), modulo the conical singularities of the table.  We would like to eventually show that a suitable analog of such a notion makes sense for the proposed Koch snowflake billiard $\Omega(KS)$. Certainly, the absence of an inward pointing vector at a point $x$ of $KS$ would preclude us from describing an initial condition $(x^0,\theta^0)$ of an orbit.  Ideally, there is an infinite collection $\Theta$ of vectors $\theta$ of inward pointing directions such that at every point $x\in KS$ and for each $\theta\in \Theta$, there exists a compatible sequence of initial conditions $(x_i^0,\theta_i^0)_{i=0}^\infty$ such that $(x,\theta)\in \mathcal{F}(\xoo)$.

We therefore ask whether it is possible to find at every point of $KS$ a collection of inward pointing vectors.  Of course, one may then simplify the question to be one regarding the existence of \textit{an} inward pointing vector.

%The orbit construction that we have detailed in \S\ref{sec:orbitsOfOmegaKS} does not depend on the existence of an inward pointing direction at any of the points of the boundary $KS$.  So long as we only consider compatible sequences of orbits, we can express an orbit of the Koch snowflake billiard $\omegaks$ as a suitable inverse limit sequence of orbits.  However, there are three points of concern: 1) each basepoint $x_i^0$ of an initial condition in a compatible sequence of initial conditions must not be a corner of $KS_j$ for any $j\geq i$, 2) each basepoint $x_i^{k_i}$ of every orbit $\orbitksi{i}$ in a compatible sequence of orbits must never be a corner of any finite approximation $\omegaksi{j}$, and 3) considering an arbitrary point $x$ of $KS$, it must be the case that there is at least one unit vector (with angle different from $\theta(\pi/3)$) pointing inward (that is not $\theta(\pi/3)$).  

\subsection{Fractal flat surfaces}

In \S\ref{sec:SKSnAsABranchedCover}, we have shown how to construct the flat surfaces $\mathcal{S}(KS_n)$ of the associated rational billiards $\Omega(KS_n)$.  The genus of $\mathcal{S}(KS_n)$ is $g_n = 3\cdot 4^n - 2$; see \cite[Lemma 1.2, p. 1022]{MasTa} for a specific formula for calculating the genus of surfaces constructed from rational billiards.  This implies that the sequence $\{g_n\}_{n=0}^\infty$ diverges to $+\infty$.  We conjecture that there exists a surface $\mathcal{S}(KS)$ that can be appropriately viewed as a suitable limit of the flat surfaces $\mathcal{S}(KS_n)$ and such that the conjectured billiard flow on $\Omega(KS)$ is dynamically equivalent to the geodesic flow on $\mathcal{S}(KS)$.  If we can indeed view $\mathcal{S}(KS)$ as a suitable limit of the prefractal flat surfaces $\mathcal{S}(KS_n)$, then it is natural to expect that the \textit{fractal flat surface} $\mathcal{S}(KS)$ should have infinite genus.

We further conjecture that the self-similarity of the Koch snowflake curve is somehow reflected in the geometry of $\mathcal{S}(KS)$.  We have gathered some experimental evidence towards this conjecture (see Figures 9 and 11 in \cite{LaNie1}), but are not yet ready to propose a precise formulation.

In summary, we conjecture that $\mathcal{S}(KS)$ is a well-defined flat surface of infinite genus, that it has a fractal, self-similar structure, and that the geodesic flow on $\mathcal{S}(KS)$ should help us define (and, in fact, should be dynamically equivalent to) the billiard flow on the Koch snowflake table $\omegaks$.  This is clearly a difficult and long-term  problem that should be the focus of significant research in the future.

\begin{remark}
When considering suitable limits of various geometric objects (such as flat surfaces associated with prefractal approximations to the given fractal billiards), the notion of Hausdorff--Gromov limit \cite{Gr} of metric spaces may also be useful, in addition to that of inverse limit.
\end{remark}

\begin{remark}
We note that it would be interesting to investigate the possible connections between the types of `fractal flat surfaces' with infinite genus conjectured (here and in \cite[Conjecture 4.7]{LaNie1}) to be associated with fractal billiards, and the (noncommutative, adelic) `fractal membranes' introduced by the first author in \cite{La3} and viewed as (noncommutative) Riemann surfaces with infinite genus. 
\end{remark}

\subsection{Unique ergodicity of the flow}
We saw in \S\ref{sec:SKSnAsABranchedCover} that $\Delta(KS_n)$  could be tiled by equilateral triangles $\Delta_{n+1}$.  In the field of mathematical billiards, there is an important property of a particular family of flat surfaces, called \textit{the Veech dichotomy}. (See [\textbf{Ve1--3}], along with, e.g., [\textbf{GuStVo,Gu1,GuJu1--2,HuSc,KaHa2, MasTa,Sm,Vo,Zo}].)

\begin{statement}[The Veech dichotomy]
For each direction $\theta$, the geodesic flow in the direction $\theta$ is either closed or uniquely ergodic.
\end{statement} 

Since the geodesic flow and the billiard flow on a rational billiard are equivalent, one may rephrase the Veech dichotomy in terms of the billiard flow on a rational billiard $\Omega(B)$ as follows.

\begin{statement}
For each inward pointing direction $\theta$, the billiard flow in the direction $\theta$ is either closed or uniquely ergodic.  Specifically, independent of the choice of the initial basepoint $x^0$, the path traced out by the orbit  $\mathscr{O}(x^0,\theta^0)$ of a billiard $\Omega(B)$ is either closed or uniquely ergodic in the billiard table $\Omega(B)$ \emph{(}as opposed to the billiard flow being uniquely ergodic in the phase space $(B\times S^1/)\sim$\emph{)}.\footnote{We discuss the phase space of the billiard dynamics in \S\ref{subsec:mathBillandBillMap}.}
\end{statement}

Heuristically, a uniquely ergodic orbit of a billiard table $B$ is one for which the orbit fills the table.\footnote{It is therefore \textit{not} periodic.}  Likewise, a uniquely ergodic flow fills the corresponding flat surface.

\begin{theorem}[Veech's Theorem \cite{Ve1,Ve2,Ve3}]
If the stabilizer of the flat surface $\mathcal{S}(B)$ is a uniform lattice in $SL(2,\mathbb{R})$, then the Veech dichotomy holds for the flat surface.
\end{theorem}

It follows from a result of Gutkin and Judge in [\textbf{GuJu1--2}], that, for every $n\geq 1$, the Veech group $\Gamma(KS_n)$ (the stabilizer of the flat surface $\mathcal{S}(KS_n)$ in $SL(2,\mathbb{R})$) is commensurate with the Veech group $\Gamma(KS_0)$.  Since $\Gamma(KS_0)$ is commensurate with $SL(2,\mathbb{Z})$, it follows from [\textbf{GuJu}] that, for every $n\geq 0$, the Veech dichotomy holds for $\Omega(KS_n)$.  

If a suitable surface $\mathcal{S}(KS)$ can be defined so that the conjectured billiard flow on $\omegaks$ is dynamically equivalent to the geodesic flow on $\mathcal{S}(KS)$, we conjecture that the Veech dichotomy will hold for the Koch snowflake billiard table $\Omega(KS)$. To such end, we have the following theorem and definition.

\begin{theorem} 
If the flow in a direction $\theta$ is uniquely ergodic in $\mathcal{S}(KS_0)$, then, for every $n\geq 1$, an orbit in the same direction and starting at a point in $\mathcal{S}(KS_n)$ that is collinear in the direction of $\theta$ intersects with (and stops at) at most one conic singularity of the surface $\mathcal{S}(KS_n)$, in both the past and future, thus remaining a uniquely ergodic direction in the surface $\mathcal{S}(KS_n)$.
\end{theorem} 

\begin{proof} 
This easily follows from the discussion in \S3 and \S4.  The precise justification, however, is somewhat lengthy, and hence the details will be provided in \cite{LaNie2}.
\end{proof}

As a result, we propose the following definition.

\begin{definition}[Compatible sequence of uniquely ergodic orbits]
If a compatible sequence of orbits is entirely comprised of uniquely ergodic orbits, then we say that it is \textit{a compatible sequence of uniquely ergodic orbits}.
\end{definition}

Once a suitable notion of \textit{fractal flat surface} has been defined and the geodesic flow on such a surface has been shown to be dynamically equivalent to the conjectured billiard flow, a natural question to ask is whether or not the inverse limit (with suitably defined transition maps) of a compatible sequence of uniquely ergodic orbits is in fact uniquely ergodic in the billiard table $\Omega(KS)$.  Some work in progress towards this goal is being carried out in \cite{LaNie2}.

In summary, in the long term, one of the over-arching goals of this billiards project is to establish some analogue of the Veech dichotomy for the proposed fractal flat surface $\mathcal{S}(KS)$ and the associated geodesic flow.  In general, we will attempt to establish some sufficient condition for the fractal analogue of the Veech dichotomy to hold for other fractal flat surfaces.

Another important and long-term goal for the continuation of this project is the study of the connections between the length spectrum (i.e., the lengths of the periodic orbits of the given fractal billiard) and the spectrum of the Dirichlet or Neumann Laplacian (or of another Hamiltonian) on the associated fractal drum (or `drum with fractal boundary', as in, e.g., \cite{La1}, \cite{La2}, \cite{LaPa}, \cite{LaNRG}, \cite[\S12.5]{La-vF}, and the relevant references therein).  We refer the interested reader to \cite[\S4]{LaNie1} (particularly, Open Problem 4.8) for a brief discussion of the existence of possible Gutzwiller-type formulas [\textbf{Gz1,2}] (as well as [\textbf{Ch,Co1-2,DuGn}]) in this context.

In conclusion, it is hoped that the present work on the Koch snowflake billiard and its prefractal, rational billiard approximations, along wih its forthcoming sequels (including \cite{LaNie2}), will eventually help lay the foundations for a general theory of fractal billiards and their associated geodesic flows on fractal flat surfaces.

\vspace{5 pt}
\noindent \textit{Acknowledgements}.  We would like to thank Benjamin Steinhurst for helpful comments on a preliminary version of this paper, as well as for discussions on inverse limits and inverse limit sequences.  We would also like to thank Sergei Troubetzkoy for additional helpful comments during a workshop on mathematical billiards. In addition, the second author would like to thank Sergei Troubetzkoy, Sergei Tabachnikov, Robert Strichartz and Alexander Teplyaev for their positive feedback and encouragement at various conferences at which the material in this paper was presented.


\begin{thebibliography}{99}

\bibitem[Ba]{Ba} M. Barnsley, \textit{SuperFractals}, Cambridge University Press, New York, 2006.

\bibitem[BaxUm]{BaxUm} A. Baxter and R. Umble, Periodic orbits of billiards on an equilateral triangle, \textit{Amer. Math. Monthly} No. 8, \textbf{115} (2008), 479--491.

\bibitem[BeDeMiSt]{BeDeMiSt} M. Begue, L. DeValve, D. Miller and B. Steinhurst, Spectrum and heat kernel asymptotics on general Laakso spaces, preprint, arXiv:0912.2176v2, 2010.

\bibitem[Bo]{Bo} N. Bourbaki, \textit{General Topology}, Springer-Verlag, Berlin, 1989 (English translation).

\bibitem[Ch]{Ch} J. Chazarain, Formule de Poisson pour les vari\'{e}t\'{e}s riemanniennes, \textit{Invent. Math.} \textbf{24} (1974), 65--82.

\bibitem[Co1]{Co1} Y. Colin de Verdi\`{e}re, Spectre du laplacien et longueur des g\'{e}od\'{e}siques p\'{e}riodiques, I et II, \textit{Compositio  Math.} \textbf{27} (1973), 83--106 and 159--184.

\bibitem[Co2]{Co2} Y. Colin de Verdi\`ere, Spectrum of the Laplace operator and periodic geodesics: thirty years after, \textit{Ann. Inst. Fourier} No. 7, \textbf{57} (2008), 2429--2463.

\bibitem[DuGn]{DuGn} J. J. Duistermaat and V. Guillemin, The spectrum of positive elliptic operators and periodic bicharacteristics, \textit{Invent. Math.} \textbf{29} (1975), 39--79. 

\bibitem[Ed]{Ed} G. Edgar, \textit{Measure, Topology and Fractal Geometry}, Springer, New York, 2000.

\bibitem[Fa]{Fa} K. J. Falconer, \textit{Fractal Geometry}: \textit{Mathematical foundations and applications}, 2nd ed., John Wiley \& Sons, Chichester, 2003.

\bibitem[GaStVo]{GaStVo} G. Galperin, Ya. B. Vorobets and A. M. Stepin, Periodic billiard trajectories in polygons, \textit{Russian Math. Surveys} No. 3, \textbf{47} (1992), 5--80.

\bibitem[Gr]{Gr} M. Gromov, \textit{Metric Structures for Riemannian and Non-Riemannian Spaces}, Modern Birkh\"auser Classics, Birkh\"auser, Basel and Boston, 2001.

\bibitem[Gu1]{Gu1} E. Gutkin, Billiards in polygons.  Survey of recent results, \textit{J. Stat. Phys.} \textbf{83} (1996), 7--26.

\bibitem[Gu2]{Gu2} E. Gutkin, Billiards on almost integrable polyhedral surfaces, \textit{Erg. Th. and Dyn. Syst.} \textbf{4} (1984), 569--584.
  	
\bibitem[GuJu1]{GuJu1} E. Gutkin and C. Judge, The geometry and arithmetic of translation surfaces with applications to polygonal billiards, \textit{Math. Res. Lett.} \textbf{3} (1996), 391--403.

\bibitem[GuJu2]{GuJu2} E. Gutkin and C. Judge, Affine mappings of translation surfaces: Geometry and arithmetic, \textit{Duke Math. J.} \textbf{103} (2000), 191--213.

\bibitem[Gz1]{Gz1} M. C. Gutzwiller, Periodic orbits and classical quantization conditions, \textit{J. Math. Phys.} \textbf{12} (1971), 343--358.

\bibitem[Gz2]{Gz2} M. C. Gutzwiller, \textit{Chaos in Classical and Quantum Mechanics}, Interdisciplinary Applied Mathematics, vol. 1, Springer-Verlag, New York, 1990.


\bibitem[HoYo]{HoYo} J. G. Hocking and G. S. Young, \textit{Topology}, Dover Publ., Mineola, 1988.

\bibitem[HuSc]{HuSc} P. Hubert and T. Schmidt, An introduction to Veech surfaces, in: \textit{Handbook of Dynamical Systems}, vol. 1B (A. Katok and B. Hasselblatt, eds.), Elsevier, Amsterdam, 2005, pp. 501--526.

\bibitem[Hut]{Hut} J. E. Hutchinson, Fractals and self-similarity, \textit{Indiana Univ. Math. J.} \textbf{30} (1981), 713--747.

\bibitem[KaHa1]{KaHa1} A. Katok and B. Hasselblatt, \textit{A First Course in Dynamics}: \textit{With a panorama of recent developments}, Cambridge Univ. Press, Cambridge, 2003.

\bibitem[KaHa2]{KaHa2} A. Katok and B. Hasselblatt, \textit{Introduction to the Modern Theory of Dynamical Systems}, Cambridge Univ. Press, Cambridge, 1995.


\bibitem[KaZe]{KaZe} A. Katok and A. Zemlyakov, Topological transitivity of billiards in polygons, \textit{Math. Notes} 
\textbf{18} (1975), 760--764.

\bibitem[KeMaSm]{KeMaSm} S. Kerckhoff, H. Masur, J. Smillie, Ergodicity of billiard flows and quadratic differentials, \textit{Annals of Math}. \textbf{124} (1986), 293--311.

\bibitem[La1]{La1} M. L. Lapidus, Fractal drum, inverse spectral problems for elliptic operators and a partial resolution of the Weyl--Berry conjecture, \textit{Trans. Amer. Math. Soc.} \textbf{325} (1991), 465--529.

\bibitem[La2]{La2} M. L. Lapidus, Vibrations of fractal drums, the Riemann hypothesis, waves in fractal media, and the Weyl--Berry conjecture, in: \textit{Ordinary and Partial Differential Equations} (B. D. Sleeman and R. J. Jarvis, eds.), vol. IV, Proc. Twelfth Internat. Conf. (Dundee, Scotland, UK, June 1992), Pitman Research Notes in Math. Series, vol. 289, Longman, Scientific and Technical, London, 1993, pp. 126--209.

\bibitem[La3]{La3} M. L. Lapidus, \textit{In Search of the Riemann Zeros}: \textit{Strings, fractal membranes and noncommutative spacetimes}, Amer. Math. Soc, Providence, R. I., 2008.

\bibitem[LaNie1]{LaNie1} M. L. Lapidus and R. G. Niemeyer, Towards the Koch snowflake fractal billiard---Computer experiments and
mathematical conjectures, in: ``\textit{Gems in Experimental Mathematics}'' (T. Amdeberhan, L. A. Medina and V. H. Moll, eds.), Contemporary Mathematics, Amer. Math. Soc., Providence, R. I., \textbf{517} (2010), pp. 231--263. [E-print: arXiv:math.DS:0912.3948v1, 2009.]


\bibitem[LaNie2]{LaNie2} M. L. Lapidus and R. G. Niemeyer, Veech groups $\Gamma_n$ of the Koch snowflake prefractal flat surfaces $\mathcal{S}(KS_n)$, work in progress, 2011.


\bibitem[LaNRG]{LaNRG} M. L. Lapidus, J. W. Neuberger, R. J. Renka and C. A. Griffith, Snowflake harmonics and computer graphics: Numerical computation of spectra on fractal domains, \textit{Internat. J. Bifurcation \& Chaos} \textbf{6} (1996), 1185--1210.

\bibitem[LaPa]{LaPa} M. L. Lapidus and M. M. H. Pang, Eigenfunctions of the Koch snowflake drum, \textit{Commun. Math. Phys.} \textbf{172} (1995), 359--376.

\bibitem[La-vF]{La-vF} M. L. Lapidus and M. van Frankenhuijsen, \textit{Fractal Geometry, Complex Dimensions and Zeta Functions: Geometry and spectra of fractal strings}, Springer Monographs in Mathematics, Springer-Verlag, New York, 2006. (Second revised and enlarged edition to appear in 2011.)

\bibitem[Ma]{Ma} W. S. Massey, \textit{Algebraic Topology: An Introduction}, Springer-Verlag, New York, 1977.

\bibitem[Mas]{Mas} H. Masur, Closed trajectories for quadratic differentials with an applications to billiards, \textit{Duke Math. J.} \textbf{53} (1986), 307--314.

\bibitem[MasTa]{MasTa} H. Masur and S. Tabachnikov, Rational billiards and flat structures, in: \textit{Handbook of Dynamical Systems}, vol. 1A (A. Katok and B. Hasselblatt, eds.), Elsevier, Amsterdam, 2002, pp. 1015--1090.

\bibitem[McL]{McL} S. Mac Lane, \textit{Categories for the Working Mathematician}, 2nd ed., Graduate Text in Mathematics, vol. 5, Springer-Verlag, New York, 1989.

\bibitem[RoSt]{RoSt} K. Romeo and B. Steinhurst, Eigenmodes of a Laplacian on some Laakso spaces, \textit{Complex Variables and Elliptic Equations} {\bf 54} (2009), 623--637. 


\bibitem[Sm]{Sm} J. Smillie, Dynamics of billiard flow in rational polygons, in: \textit{Dynamical Systems}, Encyclopedia of Math. Sciences, vol. 100, Math. Physics 1 (Ya. G. Sinai, ed.), Springer-Verlag, New York, 2000, pp. 360--382.

\bibitem[Ta]{Ta} S. Tabachnikov, \textit{Billiards}, Panoramas et Synth\`{e}ses, vol. 1, Soc. Math. France, Paris, 1995.

\bibitem[St]{St} B. Steinhurst, ``\textit{Diffusions and Laplacians on Laakso, Barlow--Evans, and other fractals}'', Ph.D. Dissertation, Univ. of Connecticut, May 2010.

\bibitem[Ve1]{Ve1} W. Veech, Teichm\"{u}ller geodesic flow, \textit{Annals of Math}. \textbf{124} (1986), 441--530.

\bibitem[Ve2]{Ve2} W. Veech, Teichm\"{u}ller curves in modular space, Eisenstein series, and an application to triangular billiards, \textit{Invent. Math.} \textbf{97} (1989), 553--583.

\bibitem[Ve3]{Ve3} W. Veech, Flat surfaces, \textit{Amer. J. Math.} \textbf{115} (1993), 589--689.

\bibitem[Vo]{Vo} Ya. B. Vorobets, Plane structures and billiards in rational polygons: The Veech alternative, \textit{Russian Math. Surveys} \textbf{51} (1996), 779--817.

\bibitem[Zo]{Zo} A. Zorich, Flat surfaces, in: \textit{Frontiers in Number Theory, Physics and Geometry} I (P. Cartier, \textit{et al}., eds.), Springer-Verlag, Berlin, 2002, pp. 439--585.

\end{thebibliography}
\end{document}